\documentclass[a4paper,12pt]{amsart}

\newcommand{\version}{20}
\title[D-minimal structures v. \version]{D-minimal structures\\
\small{Version~\version}}

\author[A. Fornasiero]{Antongiulio Fornasiero}
\address{Universit\`a di Parma (Italy)}
\email{antongiulio.fornasiero@gmail.com} 

\date{12 April 2016}

\usepackage{ifpdf}
\usepackage{amsmath, amsthm, amssymb}
\usepackage{mathtools}
\usepackage{xspace}
\usepackage[neveradjust]{paralist}
\AtBeginDocument{\setdefaultleftmargin{0pt}{}{}{}{}{}}
\setdefaultenum{\upshape(1)}{\upshape i)}{}{}
\usepackage[T1]{fontenc}
\ifpdf
\usepackage[pdftex,% linktocpage=true, 
pdfborder={0 0 0}, plainpages=false,%
pdfpagelabels, pdfdisplaydoctitle=true,%
pdfstartview={XYZ null null null},%
colorlinks=false
]{hyperref}
\hypersetup{%
  pdfauthor={A. Fornasiero},
  pdftitle={D-minimal structures},
  pdfkeywords={D-minimal, o-minimal, open core, ordered field, definably complete, 
    dense pair, locally o-minimal, Baire},
  pdfsubject={D-minality, generalizations of o-minimality}
  }
\else
\usepackage[plainpages=false,linktocpage=true,pdfpagelabels, colorlinks=false
]{hyperref}
\fi
\usepackage[backrefs, alphabetic, msc-links]{amsrefs}
\makeatletter
\@namedef{subjclassname@2010}{%
 \textup{2010} Mathematics Subject Classification}
\makeatother
\newcommand{\Sp}{\Sfam^p}

\DeclareMathOperator{\Graph}{Graph}
\newcommand*{\intro}[1]{\textbf{#1}}
\newcommand*{\Pa}[1]{\bigl( #1 \bigr)}
\newcommand*{\set}[1]{\{\,#1\,\}}
\newcommand*{\mset}[1]{\{#1\}}
\newcommand*{\abs}[1]{\lvert#1\rvert}
\newcommand*{\norm}[1]{\lVert#1\rVert}
\newcommand*{\card}[1]{\lvert#1\rvert}
\newcommand{\Nat}{\mathbb{N}}
\newcommand{\N}{\mathbb{N}}
\newcommand{\Real}{\mathbb{R}}
\newcommand{\R}{\mathbb{R}}
\newcommand{\Ralg}{\Real^{\mathrm{alg}}}
\newcommand{\Rbar}{\bar\Real}
\DeclareMathOperator{\GL}{GL}
\newcommand*{\inter}[1]{\mathring{#1}}
\DeclareMathOperator{\interior}{int}
\newcommand*{\cl}[1]{\overline{#1}}
\DeclareMathOperator{\cll}{cl}
\newcommand{\fr}{\partial}
\DeclareMathOperator{\bd}{bd}
\DeclareMathOperator{\fdim}{fdim}
\DeclareMathOperator{\lc}{lc}
\newcommand*{\nlc}[1]{\ulcorner #1 \urcorner}
\newcommand*{\inlc}[2]{{#1}^{\ulcorner #2 \urcorner}}
\newcommand{\Dis}{\mathcal{D}}
\newcommand*{\rest}[1]{\upharpoonright_{#1}}
\DeclareMathOperator{\id}{id}
\newcommand{\Zclosed}{Z\hyph closed\xspace}
\newcommand{\Zclosure}{Z\hyph closure\xspace}
\newcommand{\Zapplication}{Z\hyph application\xspace}
\DeclareMathOperator{\dcl}{dcl}
\DeclareMathOperator{\acl}{acl}
\DeclareMathOperator{\zcl}{Zcl}
\newcommand{\zclK}{\zcl^{\K}}
\newcommand{\zclM}{\zcl^{\monster}}
\DeclareMathOperator{\scl}{Scl}
\DeclareMathOperator{\sdim}{Sdim}
\DeclareMathOperator{\lexmin}{lex\, min}
\DeclareMathOperator{\lexinf}{lex\, inf}
\newcommand*{\gen}[1]{\langle #1 \rangle}
\newcommand{\app}{\leadsto}%{\hookrightarrow}
\newcommand*{\pair}[1]{\langle #1 \rangle}
\newcommand*{\tuple}[1]{\langle #1 \rangle}
\newcommand{\Am}{\mathbb A}
\newcommand{\Bm}{\mathbb B}
\newcommand{\Cm}{\mathbb C}
\newcommand{\Dm}{\mathbb D}
\newcommand{\Lm}{\mathbb L}
\newcommand{\monster}{\mathbb M}
\newcommand{\av}{{\bar a}}

\newcommand{\bv}{{\bar b}}

\newcommand{\cv}{{\bar c}}
\newcommand{\dv}{{\bar d}}

\newcommand{\fv}{{\bar f}}

\newcommand{\x}{{\bar x}}
\newcommand{\y}{{\bar y}}
\newcommand{\z}{{\bar z}}

\newcommand{\eps}{\varepsilon}
\def\Ind#1#2{#1\setbox0=\hbox{$#1x$}\kern\wd0\hbox to
  0pt{\hss$#1\mid$\hss}\lower.9\ht0\hbox to 0pt{\hss$#1\smile$\hss}\kern\wd0}
\newcommand{\ind}[1][]{\mathop{\mathpalette\Ind{}^{\!\!\!\!\rlap{$\scriptscriptstyle\textnormal{#1}$}\,\,\,\,}}}
\def\mind{\ind[M]}
\def\zind{\ind[Z]}
\def\tind{\ind[\th]}

\newcommand{\notind}[1][]{\mathrel{\not\mkern-7mu{\ind[#1]}}}
\newcommand*{\enumref}[1]{\rm{(\ref{#1})}}
\newcommand{\et}{\ \&\ }
\newcommand{\K}{\mathbb{K}}
\newcommand{\Z}{\mathbb{Z}}
\newcommand{\Rc}{\mathcal R}
\newcommand{\F}{\mathbb{F}}
\newcommand{\G}{\mathbb{G}}

\newcommand{\KC}{\K^C}
\newcommand{\KbC}{{\Kb}^C}

\newcommand{\Kb}{\overline{\K}}

\newcommand{\B}{\mathfrak B}
\newcommand{\Afam}{\mathcal A}
\newcommand{\Bfam}{\mathcal B}
\newcommand{\Cfam}{\mathcal C}
\newcommand{\Dfam}{\mathcal D}
\newcommand{\Efam}{\mathcal E}
\newcommand{\Sfam}{\mathcal S}

\newcommand{\Vfam}{\mathcal V}
\newcommand{\Part}{\mathcal P}
\newcommand{\elem}{\equiv}
\DeclareMathOperator{\Th}{Th}
\newcommand{\Lang}{\mathcal L}
\newcommand{\Langf}{\mathcal L_{OF}}
\newcommand{\Td}{T^d}
\newcommand{\Cp}{\mathcal C^p}
\newcommand{\Ck}{\mathcal C^k}
\newcommand{\Cone}{\mathcal C^1}
\newcommand{\Ctwo}{\mathcal C^2}
\newcommand{\identity}{\mathrm 1}
\newcommand{\idmatrix}{\mathbf 1}
\DeclareMathOperator{\Ker}{Ker}

\DeclareMathOperator{\reg}{reg}
\DeclareMathOperator{\isol}{isol}
\newcommand{\iminimal}{i\hyph minimal\xspace}
\newcommand{\iminimality}{i\hyph minimality\xspace}

\newcommand{\Iminimality}{I\hyph minimality\xspace}
\newcommand{\ipminimal}{constructible\xspace}
\newcommand{\ipminimality}{constructibility\xspace}

\newcommand{\dminimal}{d\hyph minimal\xspace}

\newcommand{\dminimality}{d\hyph minimality\xspace}
\newcommand{\pN}{pseudo\hyph$\Nat$\xspace}
\newcommand{\ominimal}{o\hyph minimal\xspace}

\newcommand{\DSF}{DSF\xspace}
\DeclareMathOperator{\RK}{rk}
\newcommand{\rkCB}{\RK^{CB}}
\newcommand*{\CBd}[2]{#1^{(#2)}}
\newcommand{\rkP}{\RK^{P}}
\newcommand{\rkM}{\RK}
\DeclareMathOperator{\Zrk}{Zrk}
\newcommand{\sdiff}{\mathbin{\Delta}}
\newcommand{\dcompact}{d-compact\xspace}
\newcommand{\Fs}{\mathcal F_\sigma}
\newcommand{\Gd}{\mathcal G_\delta}
\newcommand{\ao}{a.o\mbox{.}\xspace}
\def\hyph{\nobreakdash-\hspace{0pt}\relax}
\providecommand{\rom}{\textup}

\newcommand{\Wlog}{W.l.o.g\mbox{.}\xspace}
\newcommand{\wloG}{w.l.o.g\mbox{.}\xspace}
\newcommand{\eg}{e.g\mbox{.}\xspace}
\newcommand{\ie}{i.e\mbox{.}\xspace}

\newcommand{\wrt}{w.r.t\mbox{.}\xspace}
\newcommand{\st}{s.t\mbox{.}\xspace}
\newcommand{\cf}{cf\mbox{.}\xspace}
\newcommand{\Cf}{Cf\mbox{.}\xspace}
\newcommand{\tfae}{t.f.a.e\mbox{.}\xspace}
\newcommand{\Tfae}{T.f.a.e\mbox{.}\xspace}

\newcommand{\Case}[1]{\par\medskip\noindent\textsc{Case #1}:}

\newtheorem{lemma}{Lemma}[section]
\newtheorem{thm}[lemma]{Theorem}
\newtheorem{corollary}[lemma]{Corollary}
\newtheorem{conjecture}[lemma]{Conjecture}
\newtheorem{proposition}[lemma]{Proposition}
\newtheorem{open problem}[lemma]{Open problem}

\newtheorem*{proviso}{Proviso}
\newtheorem*{fact*}{Fact}
\newtheorem{fact}[lemma]{Fact}

\theoremstyle{remark}

\newtheorem{claim}{Claim}
\newtheorem*{claim*}{Claim}
\newtheorem{exercise}[lemma]{Exercise}

\theoremstyle{definition}
\newtheorem{definizione}[lemma]{Definition}

\newtheorem{remark}[lemma]{Remark}
\newtheorem{final remark}[lemma]{Final remark}
\newtheorem{example}[lemma]{Example}
\newtheorem{examples}[lemma]{Examples}
\newtheorem{question}[lemma]{Question}

\newtheorem{warning}[lemma]{Warning}

\newenvironment{sentence}[1][]{%
  \begin{list}{}{%
    \setlength\topsep{0.5ex}%
    \setlength\leftmargin{\parindent}%
  }%
  \item[#1]
 }
 {\end{list}}

\newenvironment{sentences}{%
  \begin{list}{}{%
   \setlength\parsep{0pt}%
    \setlength\topsep{0.5ex}%
    \setlength\itemsep{0.5ex}%
    \setlength\leftmargin{\parindent}%
    \setlength\labelwidth{\parindent}%
  }%
 }{\end{list}}
\begin{document}

\begin{abstract}
We study d-minimal expansions of ordered fields, and dense pairs thereof. We also consider other generalizations of o-minimality.
\end{abstract}

\maketitle

{\small
%\setcounter{tocdepth}{1}
%% only sections in the TOC
\tableofcontents
}

\makeatletter
\renewcommand\@makefnmark%
   {\normalfont(\@textsuperscript{\normalfont\@thefnmark})}
\renewcommand\@makefntext[1]%
   {\noindent\makebox[1.8em][r]{\@makefnmark\ }#1}
% change style of footnotes: (1) instead of 1 
% done here instead than in the preamble to avoid messing with amsart formatting of keywords et similia
\makeatother

\section{Introduction}

Let $\K$ be a first-order expansion of an ordered field.
Recall that $\K$ is said to be \intro{definably complete} (DC) if every definable
subset of $\K$ has a supremum in $\K \cup \set {\pm \infty}$
(see \eg \cite{FS:1} and its bibliography).
In this article we study the following generalization of o-minimality:

\begin{definizione}\label{def:dmin-intro}
$\K$ is \intro{\dminimal{}} if it is definably complete, and every definable
set $X \subset \K$ is the union of an open set and finitely many discrete sets,
where the number of discrete sets does not depend on the parameters of
definition of~$X$.
\end{definizione}

\cite{dries85} gives the first known example of a \dminimal non o-minimal
expansion of~$\Real$, 
\cites{FM05, MT06} give more examples of d-minimal expansions of~$\Real$ (and
introduce the notion of d-minimality), and
\cite{miller05} studies general properties of \dminimal expansions of~$\Real$
(and other such ``tameness'' notions).
Here we focus on the general case, when $\K$ is a \dminimal expansion of a
field, not necessarily the reals.

In~\cite{fornasiero-lomin} we studied a notion which is in between
o-minimality and \dminimality: locally o-minimal structures (see also
\cite{schoutens}).
In~\cite{F-Hier} we studied DC structures (expansing a field) in general, and
proved a dichotomy theorem: a DC structure either defines a discrete subring,
or it is ``restrained'' (see \S\ref{sec:restrained})
(restrained structures are a generalization of \dminimal ones).
Here we continue the study in \cites{fornasiero-lomin, F-Hier}
and apply their main results.

For \dminimal structures, we have the following fact: every definable set
can be partitioned into finitely many (definable and embedded) manifolds
(Proposition~\ref{prop:d-special}, which generalizes a result in
\cite{miller05} for expansions of~$\R$); we also have stronger property that a
definable set
admits a \emph{stratification} into finitely many manifolds  satisfying a
suitable version of Whitney condition (a)
(see Thm.~\ref{thm:stratification} and Prop.~\ref{prop:Whitney}); for locally
o-minimal structures we can even find a Whitney stratification
(see Corollary~\ref{cor:lomin-Whitney}, which generalizes \cite{Loi} for
o-minimal structures).

%It is trivial to see that a locally o-minimal non o-minimal structure cannot
%be an expansion of~$\Real$, and that

Other important properties of a \dminimal structure$\K$ are:
\begin{sentence}[(*)]
Every definable subset of $\K^n$
is constructible (\ie, a Boolean combination of finitely many open sets), 
definable sets have a well-behaved dimension function, 
and $\K$ has definable Skolem functions (Remark~\ref{rem:dmin-constructible}, 
Theorem~\ref{thm:imin}, and \S\ref{sec:DSF}).
\end{sentence}

We also study (\S\ref{sec:i0min} and \S\ref{sec:DSF}) the even more general notions of ``\iminimal'' structures
(\ie structures such that every definable set with empty interior is nowhere
dense) 
and ``\ipminimal'' structures (\ie structures such that every definable set is constructible), before focusing on
\dminimal structures (\S\S\ref{sec:dmin}--\ref{sec:dense}).
The reasons are: on the one hand, (*) holds in greater generality than for
\dminimal structures, and proving the result in this additional generality
does not involve extra difficulties (notice however that \cite{miller-definable-choice}
proved already that \dminimal structures have definable Skolem functions); 
on the other hand, we need to show that a
\dminimal structure is \iminimal and \ipminimal, and hence we need some
criteria for this.

Moreover, \iminimal structures are of independent interest, as shown in
\S\ref{sec:iminimal} and \cite{fornasiero-restrained}; in particular, many
restrained expansions of the real field are \iminimal
(see Fact~\ref{fact:imin-real-unary}), and \iminimal structures enjoy strong
``regularity'' properties (see Fact~\ref{fact:imin-real-Hausdorff} and
\cite{miller05} for expansion of the reals, and Theorem~\ref{thm:imin} and
Lemma~\ref{lem:imin-C1} for the general case).

A useful tools in the study of \dminimal structure is the
Pillay rank of a
definable set (see \S\ref{subsec:rank} and \S\ref{subsec:rank-2}): to every
definable set $X$ we associate %a natural number $\dim(X)$ and
an ordinal number~$\rkP(X)$.
The Pillay rank is a refinement of the dimension function, in the sense that
if $X$ is a nowhere dense subset of $Y$ then $\rkP(X) < \rkP(Y)$
(while in general we only have $\dim(X) \leq \dim(Y)$).

However, some properties from o-minimality do not generalize well to 
\dminimality: 
for instance, in general 
a set definable in \dminimal structure doesn't have Whitney
stratification (see Example~\ref{ex:no-Verdier});
moreover, in a \dminimal non o-minimal theory 
the algebraic closure does not have the exchange property
(\cite{DMS}*{2.12}),%
\footnote{In a previous version of this article, it was erroneously claimed
that a \dminimal non o-minimal theory is never rosy: we don't know if this is
true  or not; \cf Lemma~\ref{lem:asymmetry}.} 
and it is easy to show that there exists
a \dminimal structure with the Independence Property 
(\cf \cite{fornasiero-matroids}*{Example~12.7}).

\S\ref{sec:dense} is devoted to the study of dense (elementary) pairs of
\dminimal structures, and augments the results in \cite{vdd-dense} 
about dense pairs of o-minimal structures (and more generally the literature
on lovely pairs of geometric structures: see \cite{boxall} for an introduction
to the topic and a bibliography).
In \cite{fornasiero-matroids}*{\S9 and \S13} we studied another notion: \dminimal topological
structures, and proved some results about dense pairs of such structures.
Here we show that if $\K$ is a definably complete expansion of a field, and
$\K$ is \dminimal in the sense of Definition~\ref{def:dmin-intro},
then it is a \dminimal topological structure (see \S\ref{subsec:zclosure}):
hence, we can apply the results in \cite{fornasiero-matroids} to our
situation.
In particular, we have the following theorem (which is new even for expansions
of the reals).
\begin{thm}
Let $\K$ be a \dminimal structure, and $T$ be its theory.
Fix $n \in \Nat$; let $T^{n d}$ be the theory of tuples
$\Am_0 \prec \Am_1 \prec \dots \Am_n \models T$,
such that each $\Am_{i+1}$ is a proper elementary extension of $\Am_i$,
and $\Am_0$ is topologically dense in~$\Am_n$.
Then:
\begin{enumerate}
\item 
$T^{n d}$ is consistent and complete.  
\item
$\Am_n$ is the open core of
$\pair{\Am_0 \prec \Am_1 \prec \dots \Am_n} \models T^{n d}$
\rom(see Definition~\ref{def:core}\rom).
\item
Any model of $T^{nd}$ is definably complete.
\end{enumerate}
\end{thm}

In \S\ref{sec:cauchy} we show that if $\K$ is \dminimal, then $\KC$, the
Cauchy completion of~$\K$, has exactly one structure such that 
$\K \preceq \KC$.
This is new even in the case when $\K$ is o-minimal (see \cites{LS,Frecon}).

In \S\ref{sec:preliminary} we continue the study from 
\cites{DMS, FS:1, fornasiero-lomin, F-Hier} 
of definably complete %and Baire 
structures, introducing
some concepts and proving the results we need in the remainder of the paper.

Many of the results proved in this article are adaptions of results (and
proofs) either in o-minimal structures, or in expansions of~$\Real$ (and in
particular from \cites{miller05, DMS, Loi}).

\subsection*{Acknowledgments}
Thanks to L.~Kramer, and to an anonymous referee for his many suggestions on
how to improve this article.
%H.~Adler, A.~Berarducci, P.~Hieronymi , D.~Ikegami, C.~Miller, T.~Servi, K.~Tent, M.~Ziegler.

% -----------------------------

\section{Preliminaries}
\label{sec:preliminary}

\subsection{Conventions, basic definitions, and notation}

See \cite{fornasiero-lomin}*{\S2} for our main conventions and notations; 
in particular,
$\K$ will always be a definably complete expansion of a field, and
``definable'' will always mean ``definable with parameters''.
$\Lang$ is the language of~$\K$.
Moreover, $\cl X$ or $\cll(X)$ denote the topological closure of~$X$,
$\partial X$ denotes $\cl X \setminus X$,
while $\bd(X)$ denotes $\cl X \setminus \inter X$.

% -----------------------------

\subsection{Previous results}

See \cite{fornasiero-lomin}*{\S5} for general topology facts, definably
complete structures, dimension and full dimension, and 
pseudo-finite, locally closed, or constructible  sets.

\begin{definizione}[\cite{FS:1}*{\S2}]
Let $X \subseteq Y \subseteq \K^n$ be definable sets.
\begin{itemize}
\item 
$X$ is definably \intro{meager} in $Y$ (or, as we will almost always say, $X$
is meager in~$Y$) if there is a definable increasing family $\Pa{C_t: t \in
  \K}$ of subsets of~$Y$, such that each $C_t$ is a nowhere dense subset of
$Y$ (that is, the closure of $C_t$ in $Y$ has empty interior as a subset
of~$Y$), and $X \subseteq \bigcup_t C_t$.  If $Y = \K^n$, we simply say that
$X$ is meager.
\item
$Y$~is \intro{definably Baire} if, for every definable nonempty subset~$U
\subseteq Y$ which is open in~$Y$, $U$~is not meager in~$Y$.  
\item
$X$ is almost open (or \intro{\ao{}} for short) if there exists a definable open set $U \subseteq \K^n$, such that
$X \sdiff U$ is meager (in~$\K^n$);
notice that every meager set is \ao.
\end{itemize}

\end{definizione}

\begin{fact}[{\cite{hieronymi-Baire} see also \cite{F-Hier}*{Theorem~36}}]\label{fact:Baire}
$\K$ is definably Baire.
\end{fact}

See~\cite{FS:1}*{\S1--5} for more on the definably Baire property; 
a fundamental result, besides Fact~\ref{fact:Baire}, is the following definable
version of Kuratowski-Ulam Theorem, 
which implies that $\K^n$ is also definably Baire.

\begin{fact}[Kuratowski-Ulam Theorem \cite{FS:1}*{Proposition~5.4}]
Let $D$ be an \ao subset of $\K^{m + n}$.
Then, $D$ is meager iff the set
$\set{\x \in \K^m: D_\x \text{ is not meager}}$ is meager.
\end{fact}

Moreover, we will sometimes use the following result.
\begin{fact}[\cite{FS:1}*{Proposition~2.11}]\label{fact:open-meager}
Let $Y$ be definable and $U \subseteq Y$ be definable, open \rom(in~$Y$\rom)
and nonempty.
Then, $U$ is meager in $Y$ iff $U$ is meager in itself.
\end{fact}

Let $X \subseteq \K^n$ be definable.  Remember that call a definable set
$X \subseteq \K^n$
\intro{\dcompact{}} if it is closed and bounded, \intro{pseudo-finite} if it
is \dcompact and discrete, \intro{at most pseudo-enumerable} if there exist a
definable function $f$ whose domain is a discrete subset of $\K^m$ and whose
image is~$X$, and \intro{pseudo-enumerable} if it is at most pseudo-enumerable
and not pseudo-finite (see~\cite{fornasiero-enumerable}).

%One important result is the following.

Since every pseudo-enumerable set is meager,
Fact~\ref{fact:Baire} implies the following result
(see also~\cite{fornasiero-enumerable} for a proof that does not use
Fact~\ref{fact:Baire}).
\begin{fact}%[{\cite[Main Theorem]{fornasiero-enumerable}}]
\label{fact:enumerable}
$\K$ is not pseudo-enumerable.
\end{fact}

Given a definable set~$X$, let
\[
\delta(X) \coloneqq \inf\set{d(c, X \setminus \mset c): c \in X}.
\]
%Notice that 
\begin{remark}
Let $X$ be definable and bounded.
$X$ is pseudo-finite iff $\delta(X) > 0$.
\end{remark}

\begin{lemma}\label{lem:endpoints}
Let $A \subseteq \K$ be definable, closed, with empty interior.
Let $D$ be the set of endpoints of the complement of~$A$
(see \cite{fornasiero-enumerable}*{\S2}).
Then, $D \subseteq A$, $D$ is dense in~$A$,
and is at most pseudo-enumerable.
\end{lemma}
\begin{proof}
By \cite{fornasiero-enumerable}*{Corollary~4.18}.
\end{proof}

\begin{definizione}\label{def:core}
Let $\F = \pair{F; < ,\dots}$ be a structure expanding a linear order.
The \intro{open core} of $\F$
is he structure on $F$ given by a predicate for every $\F$-definable open subset
of $F^n$, as $n$ varies: see \cite{DMS}. 
\end{definizione}

We can always distinguish two cases: either there exists a closed definable
discrete unbounded subset of~$\K$, or $\K$ has locally o-minimal open core
(see \cite{fornasiero-lomin}).

\begin{lemma}\label{lem:increasing-union-discrete}
Let $(D_t: t \in \K)$ be a definable increasing family of discrete subsets of $\K^n$.
Then, $\bigcup_t D_t$ is also at most pseudo-enumerable.
\end{lemma}
\begin{proof}
Let $X := \bigcup_t D_t$.
If $\K$ has locally o-minimal open core, then each $D_t$ is pseudo-finite;
therefore, by \cite{fornasiero-lomin}*{Theorem~3.3(9)}, $X$ is also
pseudo-finite.
Otherwise, 
let $I$ be a closed definable discrete unbounded subset
 of~$\K_{\geq 1}$.
Then,
\[
X = \bigcup_{t \in I} D_t.
\]
Each set in the above union is discrete, the index set is
pseudo-enumerable, and the family of sets is definable;
thus, by \cite{fornasiero-enumerable}*{Corollary~4.16}, 
$X$ is at most pseudo-enumerable.
\end{proof}

%Here are some more results on definably complete structures.

\subsection{Functions}
%Let $\K$ be definably complete.
In this subsection, we will prove that certain subsets of $\K^n$ are meager.
%Thus, if $\K$ is Baire, those subsets cannot be all of~$\K^n$.
%By inspecting the proofs, one will see that in some cases (marked by
%(*)), we actually prove that those sets are at most pseudo-enumerable.
% and
%thus, even without assuming that $\K$ is Baire, they cannot be all of~$\K^n$ 
%(Fact~\ref{fact:enumerable}).

Let $f: A \to \K^n$ be a definable function.
For every $s > 0$, define
\begin{multline*}
\Dis(f, s) := \bigl\{x \in A: \forall t > 0,\\
f\Pa{A \cap B(x,t)} \text{ is not contained in any open ball of radius } s
\bigr\}.%
\footnotemark
\end{multline*}
\footnotetext{Note that there is a misprint in the definition of $\Dis(f, s)$
  in \cite{DMS}*{1.8}: with their definition, $\Dis(f,s)$ is not closed in~$A$.}
Then, each $\Dis (f, s)$ is closed in~$A$.
Moreover, $\Dis(f)$, the set of discontinuity
points of~$f$, is equal to $\bigcup_s \Dis(f, s)$, and, therefore, it is an
$\Fs$ subset of~$A$.

\begin{fact}[\cite{F-Hier}*{Lemma 39 and Theorem B}]
\label{fact:monotone-discontinuous}
Let $f: \K \to \K$ be definable and monotone.
Then, $\Dis(f)$ is at most pseudo-enumerable.
Moreover, $f$ is differentiable on a dense subset of $\K$.
\end{fact}

% \begin{lemma}\label{lem:monotone-discontinuous}
% Let $f: \K \to \K$ be definable and monotone.
% Then, $\Dis(f)$ is at most pseudo-enumerable.
% \end{lemma}
% \begin{proof}
% If $\K$ has locally o-minimal open core, then the result is clear (see
% \cite{fornasiero-lomin}).
% Otherwise, let $N$ be a closed definable discrete unbounded subset
% of~$\K_{\geq 1}$.
% Fix $ \varepsilon > 0$; first, we prove that
% $\Dis(f,\varepsilon)$ is closed and discrete.
% Assume not; then, \wloG, $f$ is increasing,
% and there exists $b \in \K$ such that $b$ is an accumulation point $b$ for 
% Define $m := \sup_{t < b} f(t) = \lim_{t \to b^-} f(t) \in \K$.
% Let $c \in (a,b)$ such that $f(c) > m - \varepsilon/2$.
% The fact that $f$ is increasing and our assumption on $b$ imply that
% there  exists $d \in (c, b)$ such that $f(d) \geq f(c) + \varepsilon > m +
% \varepsilon/2$, contradicting the definition of~$m$.

% Finally,
% \[
% \Dis(f) = \bigcup_{t \in N} \Dis(f, 1/t).
% \]
% Each set in the above union is discrete, the index set is
% pseudo-enumerable, and the family of sets is definable;
% thus, by \cite{fornasiero-enumerable}*{Corollary~4.16}, $\Dis(f)$ is at most pseudo-enumerable.
% \end{proof}

% The open question (see \cite{miller}) is whether, with the same hypothesis as in the
% above lemma, there exists a point $x \in \K$ such that $f$ is differentiable
% at~$x$. 

\begin{lemma}\label{lem:min-image}
Let $f: \K^n \to \K$ be definable.
Define $M_f := \set{x \in \K^n: x \text{ is a local }\\
\text{minimum for } f}$.
Then, $f(M_f)$ is at most pseudo-enumerable.
%meager in~$\K$.
\end{lemma}
\begin{proof}
For every $r > 0$, let
\[
M(r) := \set{x \in \K^n: \abs{x} \leq r \ \&\ \abs{f(x)} \leq r \ \&\
x \text{ is a minimum for } f \text{ in the ball } B(x, 1/r)}.
\]
Note that
$f(M_f) = \bigcup_r f\Pa{M(r)}$. 
Fix $r > 0$.
\begin{claim}
$Y := f\Pa{M(r)}$ is pseudo-finite.
\end{claim}
Assume, for a contradiction, that $Y$ has an accumulation point $y \in \K$.
For every $\delta > 0$, let 
\[
U(\delta) := \set{x \in M(r): f(x) \in B(y, \delta) \setminus \mset y}.
\]

Let $C(\delta)$ be the closure of $U(\delta)$ in~$\K^n$.
Since each $C(\delta)$ is a nonempty \dcompact subset of $\cl{B(0,r)}$, the
intersection of the $C(\delta)$ is nonempty; let $x \in \bigcap_\delta
C(\delta)$.
Choose $x_1$, $x_2$, and $\delta$ such that
$x_1, x_2 \in B(x,1/(2r)) \cap U(\delta)$ and $f(x_1) < f(x_2)$ (they
exist by definitions).
However, $x_1 \in B(x_2, 1/r)$, and $x_2 \in M(r)$;
therefore, $f(x_1) \geq f(x_2)$, absurd.

The fact that $f(M_f)$ is at most pseudo-enumerable follows from
the claim and Lemma~\ref{lem:increasing-union-discrete}.
%
% Now we want to prove that $f(M_f)$ is at most pseudo-enumerable.
% If $\K$ has locally o-minimal open core, then the conclusion is clear 
% (by \cite{fornasiero-lomin}*{Theorem~3.3(9)}, $f(M_f)$ is pseudo-finite).
% Otherwise, let $N$ be a closed definable discrete unbounded subset
%  of~$\K_{\geq 1}$.
% Then,
% \[
% f(M_f) = \bigcup_{t \in N} f(M(1/t)).
% \]
% Each set in the above union is pseudo-finite, the index set is
% pseudo-enumerable, and the family of sets is definable;
% thus, by \cite{fornasiero-enumerable}*{Corollary~4.16}, 
% $f(M_f)$ is at most pseudo-enumerable.
\end{proof}

%%%% More results

\begin{definizione}\label{def:first class}
A definable function $f: X \to Y$ is \intro{of first class} if there exists
a definable function $F: \K \times X \to Y$, such that, for every $t \in \K$ 
and $x \in X$,
\begin{enumerate}
\item  $f_t(x) := F(t, x) : X \to Y$ is a continuous function
of~$x$,
\item $\lim_{t \to + \infty} f_t(x) = f(x)$;
\end{enumerate}
that is, $f$ is a point-wise limit of a definable family of continuous
functions $(f_t)_{t \in \K}$.
\end{definizione}

\begin{lemma}\label{lem:first-class}
Let $X$ be definably Baire, and  $f: X \to \K^m$ be of first class.
Then, $\Dis(f)$ is meager in~$X$.
\end{lemma}
\begin{proof}
Minor variation of \cite{oxtoby}*{Thm.~7.3}.
It suffices to show that, for every $\varepsilon > 0$,
$F_\varepsilon := \set{x \in X: \omega(x) \geq 5 \varepsilon}$ is nowhere
dense, where
\[
\omega(x) := \lim_{\delta \to 0^+} 
\sup \set{ \abs{f(x') -   f(x'')}: x', x'' \in B(x; \delta)}.
\]
Fix an open definable subset $X' \subseteq X$.
For every $i \in \K$, let
\[
E_i := \bigcap_{s, t \geq i}
\set{x \in X': \abs{f_s(x) - f_t(x)} \leq \varepsilon}.
\]
Then, $(E_i)_{i \in \K}$ is an increasing family of closed subsets of~$X'$,
and $\bigcup_i E_i = X'$.
%If $X$ is meager in itself, $\Dis f$ is automatically meager.
Since $X$ is definably Baire and $X'$ is open in~$X$, $X'$ is not meager 
in itself (see Fact~\ref{fact:open-meager}),
and therefore $E_{i_0}$ must have nonempty interior (in~$X'$), 
for some $i_0 \in \K$.
Let $V$ be a definable nonempty open subset of~$E_{i_0}$.
We have $\abs{f_t(x) - f_s(x)} \leq \varepsilon$ for all $x \in V$
and $s, t  \geq i_0$.
Putting $t = i_0$ and letting $s \to \infty$, it follows that
$\abs{f_{i_0}(x) - f(x)} \leq \varepsilon$ for all $x \in V$.
For all $x_0 \in V$ there exists a neighborhood $U(x_0) \subseteq V$,
such that $\abs{f_{i_0}(x) - f_{i_0}(x_0)} \leq \varepsilon$ for all 
$x \in U(x_0)$.
Hence $\abs{f(x) - f_{i_0}(x_0)} \leq 2 \varepsilon$ for all $x \in U(x_0)$.
Therefore $\omega(x_0) \leq 4 \varepsilon$, 
and so no point of $V$ belongs to~$F_\varepsilon$.
Thus, every nonempty open set $X'$ contains a nonempty open subset~$V$
disjoint from~$F_\varepsilon$.
This shows that $F_\varepsilon$ is nowhere dense.
\end{proof}

See also \cite{F-Hier}*{Lemma~46} for a similar result (with a similar proof),
and \cite{FS:1}*{Lemma~3.10} for a related one.

%% ----

%\subsection{Bad and locally closed sets}

% \subsection{Differential equations}

% \begin{proposition}\label{prop:differential-equation}
% Let $U \coloneqq I_1 \times I_2 \subseteq \K^n \times \K^m$ be an open rectangular box,
% and $F: U \to M_n(\K)$ be a $\Cone$ definable function (where $M_n(\K)$ is the
% set of $n \times n$ matrices over~$\K$).
% For $\pair{a, b} \in\inter U$, consider the system of differential equations
% \begin{equation}\label{eq:differential}
% \left\{
% \begin{aligned}
% \phi(a) &= b,\\
% d_x(\phi) &= F(x, \phi(x)).
% \end{aligned}
% \right.
% \end{equation}
% Then, there exists \emph{at most} one function $\phi: I_1 \to I_2$, which is
% definable, $\Cone$, and satisfies 
% \eqref{eq:differential}.
% \end{proposition}
% \begin{proof}
% The same as~\cite{OPP}*{Theorem~2.3}.
% \end{proof}

\subsection{Bad set}

\begin{definizione}\label{def:bad}
Let $A \subseteq \K^{n + m}$.
The set of \intro{bad points} for $A$ is 
\[
\B_n(A) := \set{x \in \K^n: \cll(A)_x \setminus \cll(A_x) \neq \emptyset}.
\]
\end{definizione}
Notice that $\B_n(A) = \set{x \in \K^n: \cll(A)_x \neq \cll(A_x)}$.

In the following, it will often be necessary to prove that $\B_n(A)$ 
%the set of bad points of~$A$, 
is ``small'' (in some suitable sense, usually meaning ``meager'').

\begin{remark}\label{rem:bad-3}
Assume that $A \subseteq C \subseteq \cl A \subseteq \K^{n+m}$.
Then, $\B_n(A) \supseteq \B_n(C)$.
\end{remark}
\begin{proof}
\[
\cll(A)_x \setminus \cll(A_x) = \cll(C)_x \setminus \cll(A_x)
\supseteq \cll(C)_x \setminus \cll(C_x). 
\qedhere
\]
\end{proof}

\begin{lemma}\label{lem:bad-2}
If $A$ is an $\Fs$, then $\B_n(A)$ is %analytic (\ie, 
the projection of a $\Gd$ set.
If $A$ is open, then $\B_n(A)$ is a meager~$\Fs$.
\end{lemma}
\begin{proof}
% Let $\pi := \Pi^{n + m}_n$ and
% $A := \bigcup_{t \in \K} A(t)$, where $\Pa{A(t)}_{t \in \K}$ is a
% definable increasing family of \dcompact subsets of $\K^{n + m}$.

% Let $D := \set{\pair{x, y}: y \in \cll(A_x)}$.
% Notice that $\B_n(A) = \pi(\cl A \setminus D)$.
% Hence, to prove that $\B_n(A)$ is an $\Fs$ set, it suffices to show that $D$
% is a $\Gd$ set.
% For every $r > 0$ and $t \in \K$, let 
% $D(t,r) := \set{\pair{x, y}: d(y, A(t)_x) \leq r}$.
% Notice that $D = \bigcap_{r > 0} \bigcup_{t} D(t,r)$.

% \begin{claim}
% For every $t$ every $r > 0$, $D(t,r)$ is \dcompact.
% \end{claim}
% It suffices to show that each $D(t,r)$ is closed.
% Let $p_0 := \pair{x_0, y_0} \notin D(t,r)$; thus, $d(y_0, A(t)_{x_0}) > r$.
% Assume, for a contradiction, that $p_0 \in \cl{D(t,r)}$.
% Then, for every $s > 0$, there exist $x \in B(x_0: s)$, $y \in B(y_0, s)$, 
% such that $d(y, A_x) \leq r$.
% Hence, for any such~$x$, there exists $z'$ such that $d(z', B(y_0, s)) \leq r$ 
% and $\pair{x, z'} \in A$.
% Let $z$ be an accumulation point, as $s \to 0$, of the sets 
% \[
% Z_s := \set{z' \in \K^m: d(z', B(y_0; s)) \leq r \et \exists x \in B(x; s)\
%   (\pair{x, z'} \in A)}.
% \]
% Notice that $d(z, y_0) \leq r$.
% However, since $A$ is closed, $\pair{x_0, z} \in A$, and therefore 
% $d(y_0, A_x) \leq r$, absurd.

% \textbf{Discard the previous ``proof''.}

Define
\begin{gather*}
F := \set{ \pair{x, y, r, y'} \in \K^n \times \K^m \times \K \times \K^m: 
r > 0 \et \pair{x, y'} \in A \et \abs{y - y'} < r};\\
\begin{aligned}
\pi(x, y) &:= x;\\
\pi'(x, y, r) &:= \pair{x, y};\\
\pi''(x, y, r, y') &:= \pair{x, y, r}.
\end{aligned}
\end{gather*}
%\begin{claim}
Notice that
\[
\B_n(A) = \pi\Pa{\cl A \cap \pi'(\set{r > 0} \setminus \pi''(F))}.
\]

%\end{claim}
If $A$ is an $\Fs$, then $F$ is also an~$\Fs$, and therefore $\B_n(A)$ is
the projection of a $\Gd$ set.

Assume now that $A$ is open.
Then, $F$ is also open, and therefore $\B_n(A)$ is an~$\Fs$.
For every $r > 0$, define 
$C(r) := \set{\pair{x,y} \in \cl A: \abs{x,y} \leq 1/r \et  d(y, A_x) \geq r}$.
Notice that $\B_n(A) = \bigcup_{r > 0} \pi(C(r))$ and
%Assume that $A$ is open.
each $C(r)$ is \dcompact.
Hence, to prove that $\B_n(A)$ is meager, it suffices to prove that, for every
$r > 0$, $\pi(C(r))$ has empty interior.
%Wlog, $\K$~is Baire.
Assume, for a contradiction, that $\pi(C(r_0))$ contains a nonempty open
box~$W$, for some $r_0 > 0$. To simplify the notation, assume that $m = 1$.
Define $f: W \to \K$, $f(x) := \min\Pa{C(r_0)_x}$.
By \cite{DMS}*{Lemma~2.8},
%\cite{FS}, 
there exists a nonempty open box $W' \subseteq W$, such that
$f \rest {W'}$ is continuous; \wloG, $W = W'$.
Fix $x_0 \in W$, call $y_0 := f(x_0)$, and let $V_{x_0}$ be an open box
around~$x_0$ contained in~$W$, and such that, for every $x \in V_{x_0}$,
$d(f(x), y_0) < r_0/4$. 
Since $\pair{x_0, y_0} \in \cl A$, there exists $\pair{x, y'} \in A$,
such that $x \in V_{x_0}$, and $d(y_0, y') < r_0/4$.
Let $y := f(x)$.
Since $\pair{x, y} \in C(r)$, $d(y, A_x) \geq r$; in particular, 
$d(y, y') \geq r$.
However, this contradicts $d(y, y_0) < r_0/4$ and $d(y_0, y') < r_0/4$.
\end{proof}

\begin{remark}
$\B_n(A \cup B) \subseteq \B_n(A) \cup \B_n(B)$.
\end{remark}

\subsection[Locally closed sets]{Locally closed and constructible sets}
%\acrobat

Let $X$ be a topological space $X$, and $Y \subseteq X$.
$Y$~is locally closed (in $X$) if it is of the form $U \setminus V$, where $U$ and
$V$ are open subsets of~$X$.
$Y$~is constructible if it is a finite Boolean combination of open subsets
of~$X$.
If we don't specify, we take $X = \K^n$ (for some $n \in \N$).

See \cite{fornasiero-lomin}*{\S5.6} for basic results on definable locally closed  and constructible sets; we recall here the definition and our notation.
See also \cite{robinson74} and \cite{pillay87}.

\begin{definizione}[\cite{fornasiero-lomin}]\label{def:lc}
Let $X \subseteq \K^n$.
Define
\[
\lc(X) := \set{x \in X : X \text{ is locally  closed at } x} 
\]
(that is,
$x \in \lc(X)$ iff there exists an open ball $B$ of center~$x$, such that $X
\cap B = \cl X \cap B$), and $\nlc X := X \setminus \lc(X)$.
% \end{definizione}

% \begin{definizione}
Define $\inlc X 0 := X$, and, for each $k \in \Nat$,
$\inlc X {k + 1} := \nlc{\inlc X {k}}$.
\end{definizione}
%Notice that $\lc(X)$ is locally closed, and therefore constructible.
%Notice also that, if $X$ is definable, then also $\nlc X$ and each $\inlc X {k}$
%are definable.
One can easily see that $\lc(X)$ is locally closed, and that
if $X$ is an~$\Fs$, then $\nlc X$ is also an~$\Fs$
(see \cite{fornasiero-lomin}).

Here is the basic result on definable constructible sets.
\begin{fact}
$\nlc A  = A \cap \partial(\partial A)$.
$A$ is the union of $m$ locally closed sets if and only if $\inlc A {m + 1}$
is empty.
\end{fact}
\begin{proof}
See~\cite{allouche96}, where $\partial A$ is denoted by $\check A$,
and $\nlc A$ by either $\mathcal B(A)$ or~$H(A)$; see also~\cite{DM} for another proof.
\end{proof}

\begin{fact}%\cite{robinson74}%
\footnote{See the proof of \cite{pillay87}*{Lemma~2.3}}
%\footnote{Give precise reference. \Cf \cite{pillay87}*{Lemma~2.3}}
Let $X$ be a topological space and $Y \subseteq X$ be a constructible subset.
Then, $Y$~is nowhere dense iff it has empty interior (inside~$X$).
\end{fact}

\begin{proposition}\label{prop:lc-fiber}
Let $A \subseteq \K^n$ be definable and locally closed.
Let $U \subseteq \K^n$ be open, such that $A = \cl A \cap U$.
Let $d \leq n$.
Then, for every $x \in \K^d$,
\[
A_x = \cll(A_x) \cap U_x,
\]
and in particular $A_x$ is locally closed.
Moreover, $\B_d(A) \subseteq \B_d(U)$, and therefore $\B_d(A)$ is meager.
\end{proposition}
\begin{proof}
$A_x \subseteq \cll(A_x) \cap U_x$ is obvious.
For the opposite inclusion, 
\[
\cll(A_x) \cap U_x \subseteq (\cl A)_x \cap U_x =
(\cl A \cap U)_x = A_x.
\]
Assume, for a contradiction, that $x \in \B_d(A) \setminus \B_d(U)$.
Let $E := \cl A$; notice that $A = E \cap U$.
Since $x \notin \B_d(U)$, $\cll(U_x) = (\cl U)_x$.
%By the above Remark, applied to~$E_x \cap U_x$ and to $E \cap U$,
Notice that $\cll(E \cap U) = E \cap \cl U = E$, and therefore
\[
\cll(A_x) = \cll(E_x  \cap U_x) = E_x \cap \cll(U_x) =
E_x \cap (\cl U)_x = (E \cap \cl U)_x = E_x = (\cl A)_x,
\]
contradicting $x \notin \B_d(A)$.

By Lemma~\ref{lem:bad-2}, $\B_d(U)$ is meager, and we are done.
\end{proof}

\begin{corollary}\label{cor:small-bad}
Let $A \subseteq \K^n$ be definable and constructible, and $d \leq n$.
Then, $\B_d(A)$ is meager.
\end{corollary}

% --------------------------------------------------------------------

%\acrobat

\subsection{Pettis theorem}
%Let $\K$ be definably complete and \textbf{Baire}.

If $X$ and $Y$ are subsets of $\K^n$,
then $X - Y := \set{x - y: x \in X, y \in Y}$.

\begin{lemma}[Pettis Theorem]
Let $A \subseteq \K^n$ be definable and  \ao.
If $A$ is nonmeager, then $A - A$ contains a nonempty open neighborhood
of~$0$.
Conversely, if $\K$ is \iminimal \rom(see \S\ref{sec:i0min}\rom),
%with Definable Skolem Functions \rom(\DSF for short:see \S\ref{sec:DSF}\rom), 
$n = 1$,
and $A - A$ is nonmeager, then $A$ is nonmeager.%
\footnote{The original version of this and the following lemma required that
  $\K$ has Definable Skolem Functions. 
Thanks to an anonymous referee for pointing out that the
  additional assumption was not necessary.}
\end{lemma}
\begin{proof}
Minor variation of \cite{oxtoby}*{Thm.~4.8}.
Let $A = U \sdiff P$, where $U$ is open and definable, and $P$ is meager.
$A$~is nonmeager iff $U$ is nonempty.
If $A$ is nonmeager, let $B \subseteq U$ be a nonempty open ball, of radius $\delta > 0$.
For any $x \in \K^n$, we have
%\[
\begin{multline*}
(x + A) \cap A = \Pa{(x + U) \sdiff (x + P)} \cap (U \sdiff P) = \\
= \Pa{(x + U) \cap U} \sdiff \Pa{(x + U) \cap P} \sdiff \Pa{(x + P) \cap U}
\sdiff \Pa{(x + P) \cap P} \supseteq \\
\supseteq [(x + B) \cap B] \sdiff [P \cup (x + P)].
\end{multline*}
%\]
If $\abs x < \delta$, the right member represents a nonempty open set,
minus a meager set; it is therefore nonempty.
Thus, for every $x \in B(0;\delta)$, $(x + A) \cap A$ is nonempty,
and therefore $x \in A - A$.

Conversely, assume, for a contradiction, that $\K$ is \iminimal,
$A \subset \K$ is meager, but $A - A$ is nonmeager.
Then, by Theorem~\ref{thm:imin},
by replacing $A$ with $\cl A$, \wloG we can assume that $A$
is closed and nowhere dense, while $A - A$ contains a nonempty interval~$I$.
Let $D \subseteq A$ be the set of endpoints of $\K \setminus A$.
By Lemma~\ref{lem:endpoints}, $D$ is a pseudo-enumerable subset of~$A$, and
$D$ is dense in~$A$.
Since the function $\pair{x,y}\mapsto x- y$ is continuous, $D - D$ is dense in~$I$.
By Theorem~\ref{thm:imin} again, $D - D$ has nonempty interior; however, $D -
D$ is pseudo-enumerable, a contradiction.
%
% if \iminimal and $A - A$ is nonmeager, then, by
% Thm.~\ref{thm:imin}, $\dim (A - A) = 1$; therefore, if $\K$ has DSF,
% Lemma~\ref{lem:dim-function} implies that  $\dim (A) = 1$,  
% and thus $A$ is nonmeager.
\end{proof}

Define $F: \K^4 \to \K$ as $F(x_1, x_2, y_1, y_2) := (x_1 - x_2) / (y_1 - y_2)$
if $y_1 \neq y_2$, and $0$ otherwise.

\begin{corollary}
Let $A \subseteq \K$ be definable and \ao.
If $A$~is nonmeager, then $F(A^4) = \K$.
If $\K$ is \iminimal 
% with \DSF \rom(see \S\ref{sec:DSF}\rom) 
and $F(A^4)$ is nonmeager, then $A$ is nonmeager.
\end{corollary}
\begin{proof}
If $A$ is nonmeager, then $A - A$ contains an open neighborhood of~$0$, 
and therefore $F(A^4) = \K$.
The converse is proved as in the previous lemma.
\end{proof}

See also \cite{fornasiero-matroids}*{Lemma~3.47}.

% ---------------------------------------------------------------

%\acrobat

\section{Restrained, i-minimal and constructible structures}\label{sec:i0min}
As usual, $\K$ is a definably complete structure, expanding a field.

\subsection{Restrained structures}\label{sec:restrained}

\begin{definizione}[\cite{F-Hier}*{Definition~23 and Theorem~A}]
$\K$ is \intro{restrained} if, for every definable discrete subset $D
\subseteq \K^n$, and every definable function $f: D \to \K$,
$f(D)$ is nowhere dense in~$\K$.
\end{definizione}

%In \cite{F-Hier}*{Theorem~A}, we proved a dichotomy theorem: either $\K$ is
% restrained, or it defines a discrete subring.
The structures of interest to us in this article are restrained.

\begin{fact}[\cite{F-Hier}*{Lemma~49 and Theorem~A}]
\label{fact:restrained-meager}
\Tfae:
\begin{enumerate}
\item $\K$ is restrained;
\item $\K$ does not define a discrete subring; 
\item for every $m \in \N$, every meager subset of
$\K^m$ is nowhere dense.
\end{enumerate}
\end{fact}

\begin{lemma}\label{lem:restrained-meager}
$\K$ is restrained iff every meager subset of $\K$ is
nowhere dense.
\end{lemma}
\begin{proof}
The ``only if'' direction is clear from Fact~\ref{fact:restrained-meager}.

For the opposite direction, assume for a contradiction that $\K$ is not restrained.
Then, there exists $Z \subset K$ definable discrete subring.
Let $X \subset \K$ be the field of fractions of~$Z$.
Then, $X$ is pseudo-enumerable, and therefore it is meager;
however, $X$ is dense in $\K$.
\end{proof}

\begin{fact}[\cite{F-Hier}*{Lemma~51}]\label{fact:restrained-Cp}
Let $\K$ be restrained, $U \subseteq \K^n$  be open and definable,
$f: U \to \K$ be a definable continuous function,
and $p \in \N$.
Then, $f$ is $\Cp$ on a dense open subset of~$U$.
\end{fact}

\subsection{I-minimal structures}\label{sec:iminimal}

\begin{definizione}
$\K$ is \intro{\iminimal{}} if, for every unary definable set~$X$, 
if $X$ has empty interior, then $X$ is nowhere dense.
\end{definizione}

See also~\cite{miller05} for the case when $\K$ is an expansion of~$\Real$.

\begin{remark}
If $\K$ is \iminimal, then it is restrained.
\end{remark}
\begin{proof}
By Lemma~\ref{lem:restrained-meager}.
\end{proof}

In \cite{fornasiero-restrained} we prove the following two facts about
expansions of the real field.
\begin{fact}\label{fact:imin-real-Hausdorff}
Let $\Rc$ be an expansion of the real field.
Then, \tfae:
\begin{enumerate}
\item $\Rc$ is \iminimal;
\item every unary definable set with empty interior has Lebesgue measure~$0$;
\item every unary definable set with empty interior has Hausdorff
dimension~$0$.
\end{enumerate}
\end{fact}
\begin{fact}\label{fact:imin-real-unary}
Let $\Rc$ be an o-minimal expansion of the real field, and $C \subseteq \Real$
be a closed set.
The, either $\pair{\Rc, C}$ is \iminimal, or 
%$\pair{\Rc, C}$ is not restrained \rom(that is, 
the set of natural numbers is definable in $\pair{\Rc, C}$.
\end{fact}
Thus, several of the ``tameness'' conditions introduced in
\cite{miller05}*{\S3.1} are equivalent to \iminimality,
and, for expansion of the real field by a unary closed set, \iminimality is
equivalent to being restrained.
%not defining~$\Nat$.

% \begin{lemma}
% If $\K$ is \iminimal, then it is Baire
% \end{lemma}
% \begin{proof}
% By \cite[Lemma~7.1]{fornasiero-enumerable},
% if $\K$ were not Baire, it would contain a definable
% dense and co-dense subset~$Q$.
% However, $Q$ would have empty interior, but it would not be nowhere dense.
% \end{proof}

\begin{examples}
If $\K$ is locally o-minimal, then it is \iminimal (trivial).

If $\K$ is d-minimal, then it is \iminimal.
In fact, by definition of d-minimality if $X \subseteq \K$ is definable and
with empty interior, then $X$ is a finite union of discrete sets~$X_1, \dotsc,
X_n$.
%By changing the $X_i$ if necessary, we can assume that the $X_i$ are also
%definable.
Every discrete subset of a definably complete structure is nowhere dense.
Thus, $X$~is a finite union of nowhere dense sets,
and thus it is nowhere dense.
\end{examples}

\begin{thm}\label{thm:imin}
The following are equivalent:
\begin{enumerate}
\item\label{en:i-imin-1}
$\K$ is \iminimal;
\item\label{en:i-imin-n} 
for every $n\in \Nat$, if $X$ is a definable subset of
$\K^n$ with empty interior, then $X$ is nowhere dense;
\item\label{en:i-bd}
for every definable set~$X$, $\bd(X)$ has empty interior;
\item\label{en:i-dim-1}
for every definable $X \subseteq \K$, $\dim X = \dim \cl X$;
\item\label{en:i-dim-n}
for every definable~$X$, $\dim X = \dim \cl X$;
\item\label{en:i-f} 
if $U \subseteq \K^n$ is definable and open, and $f : U \to \K$ is definable, then $\Dis(f)$ is nowhere dense;
\item\label{en:i-B} 
for every $n, m \in \Nat$, if $A \subseteq \K^{n+m}$ is definable,
then $\B_n(A)$~is nowhere dense (see Definition~\ref{def:bad});
\item\label{en:i-fr}
for every $n, m \in \Nat$, if $A \subseteq \K^{n+m}$ is definable, then the
set
\[
\set{x \in \K^n: (\fr A)_x \neq \fr(A_x)}
\]
is nowhere dense;
\item\label{en:i-meager-1} 
%$\K$ is Baire, and, 
for every definable $X \subseteq \K$, either $X$ has interior, or it is meager;
\item\label{en:i-meager-n}
%$\K$ is Baire, and 
every definable set either has interior, or it is meager;
\item\label{en:i-dim-union-1}
for all definable~$A, B \subseteq \K$,
$\dim (A \cup B) = \max \set{\dim A,  \dim B}$;
\item\label{en:i-dim-union-n}
for all definable~$A, B \subseteq \K^n$,
$\dim (A \cup B) = \max \set{\dim A,  \dim B}$;
\item\label{en:i-dim-fiber}
for all definable~$A \subseteq \K^n$, if $\dim A = d$,
then $\set{x \in \K^d: \dim A_x > 0}$ is nowhere dense;
\item\label{en:i-dim-fiber-2}
let $d, k, m, n \in \Nat$, with $k \leq n$ and $d \leq m$;
let $A \subseteq \K^{n + m}$ be definable, and $\dim A \leq d + k$;
define $C:= \set{x \in \K^n: \dim A_x \geq d}$; then, $\dim(C) \leq k$;
\item\label{en:i-countable-1}
any at most pseudo-enumerable union of subsets of $\K$ with empty interior has
empty interior;
\item\label{en:i-countable-n}
for every $d \leq n \in \Nat$, any at most pseudo-enumerable union of subsets
of $\K^n$ of dimension less or equal to $d$ has dimension less or equal to~$d$;
\end{enumerate}
Moreover, if $\K$ is \iminimal, then:
\begin{enumerate}[\rm(I)]
\item\label{en:i-meager}
every meager set is nowhere dense;
\item\label{en:i-increasing}
for every $d \leq n \in \Nat$, any increasing definable union of subsets of
$\K^n$ of dimension less or equal to $d$ has dimension less or equal to~$d$;
\item\label{en:i-monotonicity}
if $U \subseteq \K$ is open and definable, and $f: U \to \K$ is definable,
then there exists $D \subseteq \K$ definable, closed and with empty
interior, such that, for every definably connected component $I$ of
$U \setminus D$,  $f \rest I$ is continuous, and either constant or strictly
monotone; 
\item\label{en:i-ao}
every definable set is \ao.
%
% \item\label{en:i-continuous}
% if $U \subseteq \K^n$ is open and definable, and $f: U \to \K$ is definable,
% then there exists $D \subseteq \K^n$ definable, closed and with empty
% interior, such that the restriction $f \rest U \setminus D$ is continuous.
\end{enumerate}
\end{thm}
The proof is postponed to Section~\ref{sec:imin-proof};
\cf \cite{miller05}*{Main Lemma and Thm.~3.3}.
We record now some consequences of the above theorem.

\begin{example}
Let $\pair{M', M}$ be o-minimal structures (expanding a field), such that $M$
is a proper
elementary substructure of $M'$ and it is dense in~$M'$.
Thus, the structure $N := \pair{M', M}$ has o-minimal open core 
(see~\cite{DMS}).
%is \bminimal.
Therefore, if $X \subseteq N$ is meager, then $X$ is nowhere dense
(see \cite{fornasiero-lomin}*{\S4}).
However, $N$ is not \iminimal, because $M$ is a definable dense subset of $N$
%such that $\cl M = N$ 
with empty interior
(thus, clause~\eqref{en:i-meager} in Thm.~\ref{thm:imin}
does not imply \iminimality).
\end{example}

\begin{lemma}\label{lem:i-l-min}
$\K$ is locally o-minimal iff it is \iminimal and its open core is locally o-minimal.
\end{lemma}
\begin{proof}
The ``only if'' direction is clear.
Let us prove the ``if'' direction.
Let $X \subseteq \K$ be definable and with empty interior.
By \iminimality, $X$ is nowhere dense.
Since the open core of $\K$ is locally o-minimal, $X$~is pseudo-finite
(by \cite{fornasiero-lomin}*{Lemma~3.2}, applied to the open core of~$\K$ and
the closure of~$X$).
Thus, $\K$~is locally o-minimal.
\end{proof}

\begin{corollary}\label{cor:i-o-min}
Assume that $\K$ is \iminimal.
Then, $\K$~is o-minimal iff every definable discrete set is finite.
\end{corollary}
\begin{proof}
The ``only if'' direction is clear.
Assume now that $\K$ is \iminimal and every definable discrete set is finite.
By \cite{fornasiero-lomin}*{Corollary~4.6}, $\K$~has o-minimal open core.
Thus, by Lemma~\ref{lem:i-l-min}, $\K$~is locally o-minimal, and therefore it
coincides with its own open core (\cite{fornasiero-lomin}*{\S7}).
\end{proof}

 Since there do exist \iminimal structures that are not locally o-minimal
(\eg, \dminimal not o-minimal expansions of the real field), we have that for
some \iminimal structure there is some definable nonempty set $X$ such that
$\dim(\fr X) = \dim X$ (notice that if $\K$ is \iminimal, then $\dim(\fr X)
\leq \dim X$, because $\dim X = \dim(\cl X) = \max\set{\dim(X), \dim(\fr X)}$).

\begin{proviso}
For the remainder of this subsection, we will assume that $\K$ is \iminimal.
\end{proviso}

\begin{lemma}\label{lem:imin-C1}
Let $f : U \to \K$ be definable, where $U \subseteq \K^n$ is open and definable.
Then, for every $p \in \Nat$, there exists $D \subset U$ closed, definable and
nowhere dense, such that $f$ is $\Cp$ on $U \setminus D$.
\end{lemma}
\begin{proof}
The case $p = 0$, %we have already proved it in
            %Thm.~\ref{thm:imin}(\ref{en:i-continuous}).
is Thm.~\ref{thm:imin}(\ref{en:i-f}).

The case $p > 0$ follows from the case $p = 0$ and Fact~\ref{fact:restrained-Cp}.
\end{proof}

\begin{lemma}\label{lem:locally-nd}
Let $d \leq n$, $A \subseteq \K^{n}$ be definable, $\pi := \Pi^{n}_d$, and
\[
Z \coloneqq Z(A) := \set{a \in A: \exists U \text{ neighborhood of } a:
\pi(A \cap U) \text{ is nowhere dense}}.
\]
Then, $Z$ is a definable open \rom(in~$A$\rom) subset of~$A$, and $\pi(Z)$ 
is nowhere dense \rom(in~$\K^d$\rom).
\end{lemma}
\begin{proof}
Follows immediately from %Lemma~\ref{lem:bad-1}.
\cite{fornasiero-lomin}*{Lemma~5.25}.
\end{proof}

\begin{definizione}
We define \intro{$\Pi$-good} sets as in~\cite{miller05}*{\S7}.
That is, we say that a definable set $A \subseteq \K^{n+m}$ is $\pi$-good
(where $\pi := \Pi^{n + m}_m$) if:
\begin{itemize}
\item $\dim A = m$;
\item $\pi A$ is open;
\item $\pi(A \cap U)$ has interior for every $a \in A$ and open neighborhood
$U$ of~$a$;
\item for all $x \in \pi A$, $\dim(A_x) = 0$ and $\cll(A_x) = \cll(A)_x$.
\end{itemize}
More generally, $A$ is $\mu$-good (where $\mu$ is a projection from
$\K^{n + m}$ to an $m$\hyph dimensional coordinate space) if there is a
permutation of coordinates $\sigma$ such that $\mu = \pi \circ \sigma$, and
$\sigma A$ is $\pi$-good.
Finally, $A$ is $\Pi$-good if it is $\mu$-good for some $\mu$ as above,
and a collection of sets is $\Pi$-good if each of its elements is.
\end{definizione}

\begin{lemma}[Partition Lemma]\label{lem:partition}
Let $\Afam$ be a finite collection of definable subsets of~$\K^n$.
Then, there exists $\Bfam$, a finite $\Pi$-good partition of $\K^n$ compatible
with~$\Afam$ 
\rom(that is, $\Bfam$ is a finite collection of sets, 
and every set in $\Afam$ is a union of sets in~$\Bfam$\rom).
\end{lemma}
\begin{proof}
The proof proceeds as in \cite{miller05}*{\S7, Partition Lemma}, using  
Lemma~\ref{lem:locally-nd}.
\end{proof}

\begin{lemma}\label{lem:i-lc}
Let $A \subseteq \K^{n + m}$ be definable.
Then,
\begin{enumerate}
\item $\set{x \in \K^m: \lc(A_x) \neq (\lc A)_x}$ is nowhere dense;
\item for each $k \in \Nat$, 
$\set{x \in \K^m: (\inlc A k)_x \neq \inlc{(A_x)}{k}}$
is nowhere dense \rom(see Definition~\ref{def:lc}\rom);
\item if $\set{x \in \K^m: \lc(A_x) \neq \emptyset}$ is somewhere dense, then
$\lc(A) \neq \emptyset$.
\end{enumerate}
\end{lemma}
\begin{proof}
The same as \cite{miller05}*{Lemma~8.1}.
\end{proof}

% \begin{conjecture}[Monotonicity theorem]
% Let $f: (a, b) \to \K$ be definable.
% Then, there exists a closed nowhere dense set $D \subseteq (a,b)$,
% such that, on each definably connected component $I$ of $U$, $f$ is
% continuous, and either constant, or strictly monotone,
% where $U := \K \setminus D$.
% \end{conjecture}
% \begin{proof}
%\end{proof}

\subsection{Proof of Thm.~\ref{thm:imin}} \label{sec:imin-proof}
%\mbox{}\indent
($\ref{en:i-imin-1} \Leftrightarrow \ref{en:i-dim-1}$) and ($\ref{en:i-dim-n} \Rightarrow \ref{en:i-dim-1}$) are clear.

For every $0 < n \in \Nat$, and $K \geq 2$, let $(K)_n$ be the instantiation at
$n$ of the $K$th statement.
For instance, $\eqref{en:i-imin-n}_1$ is equivalent to \enumref{en:i-imin-1}.
We will prove that $\eqref{en:i-imin-n}_n  \Rightarrow \eqref{en:i-B}_n
\Rightarrow \eqref{en:i-f}_n \Rightarrow \eqref{en:i-imin-n}_n$, 
that $\eqref{en:i-imin-n}_n \Rightarrow \eqref{en:i-imin-n}_{n+1}$,
that $\eqref{en:i-imin-n}_n \Leftrightarrow \eqref{en:i-bd}_n$,
and that $\eqref{en:i-imin-n}_n$ implies that every meager $X \subset \K^n$ is nowhere dense.

By induction on~$n$, the above would imply
($\ref{en:i-imin-1} \Leftrightarrow \ref{en:i-imin-n} 
\Leftrightarrow \ref{en:i-bd}
\Leftrightarrow \ref{en:i-dim-1} \Leftrightarrow \ref{en:i-f} \Leftrightarrow
\ref{en:i-B} \Rightarrow \mathrm{\ref{en:i-meager}}$).

($(\ref{en:i-imin-n})_n \Leftrightarrow (\ref{en:i-bd})_n$) is clear.

($(\ref{en:i-f}_n) \Rightarrow (\ref{en:i-bd}_n)$).
Let $X \subseteq \K^n$ be definable.
Let $f:= 1_X$ be the characteristic function of~$X$;
then, $\Dis(f) = \bd(X)$; thus, $\bd(X)$ is nowhere dense.

($(\ref{en:i-B}_n) \Rightarrow (\ref{en:i-f}_n)$).
Let $f: U \to \K$ be definable, with $U \subseteq \K^n$ open.
We want to prove that $\Dis(f)$ is nowhere dense;
\wloG, $f$ is bounded and $U = \K^n$.
Let $A := \Graph(f)$; then, $\Dis(f) = \B_n(A)$.
Thus, $\Dis(f)$ is nowhere dense.

Assume now that we have $(\ref{en:i-imin-n}_n)$.

Note that if $Y \subseteq \K^n$ is meager and definable, then, since $\K$ is
definably Baire, $Y$~has empty interior, and thus
$Y$ is nowhere dense; therefore, we have proved $\eqref{en:i-meager}_n$.
%the ``moreover''
%statement. 

We prove now $\eqref{en:i-B}_n$.
Let $A \subseteq \K^{n + m}$ be definable, and $B := \B_n(A)$.

If we prove that $B$ is meager, then, since $\K$ is definably Baire, 
$B$~has empty
interior, and thus, by inductive hypothesis, $B$~is nowhere dense.
\Wlog, $A$~is bounded (because, after using a definable homeomorphism from
$\K$ to $(0, 1)$, $B$~can only become larger).

Let $\pi := \Pi^{n + m}_n$, and $U := \pi(A)$.
If $U$ has empty interior, then, by $\eqref{en:i-imin-n}_n$,
$U$~is nowhere dense; thus, $B$ is also nowhere dense, because $B \subseteq \cl{\pi(A)}$.
Therefore, we can assume that $U$ has nonempty interior, and hence, \wloG,
that $U$ is open.

Thus, for every $r > 0$, let
\[
C(r) := \set{(x,y) \in \cl A: d(y, A_x) \geq r},
\]
and $B(r) := \pi\Pa{C(r)}$.
Since $B = \bigcup_r B(r)$, and by $\eqref{en:i-imin-n}_n$,
it suffices to prove that each $B(r)$ has empty interior.

Fix $r > 0$, and assume, for a contradiction, that $U' \subseteq B(r)$ is a
nonempty open definable set.
\Wlog, $U' = U$.
Let $C' := \cl{C(r)}$: notice that $C'$ is \dcompact.
Define
$g: U \to \K^m$, $x \mapsto \lexmin (C'_x)$.
%Since $g$ is lower semi-continuous, 
By \cite{DMS}*{2.8(1)}, $\Dis(g)$ is meager, and therefore nowhere dense;
thus, there exists $U' \subseteq U$ open and nonempty, 
such that $g$ is continuous on~$U'$; \wloG, $U = U'$.
Define also $f : U \to \K^m$; $x \mapsto \lexinf\Pa{C(r)_x}$.
Note that $f(x) \geq g(x)$ for every $x \in U$.

\begin{claim}\label{cl:disc-nd}
The set $D := \set{x \in U: f(x) > g(x)}$ is nowhere dense.
\end{claim}
We will do only the case $m = 1$.
It suffices to prove that $D$ is meager.
For every $s > 0$, let $D(s) := \set{x \in U: f(x) \geq g(x) + s}$.
If we prove that each $D(s)$ is nowhere dense, we have the claim.
By $\eqref{en:i-imin-n}_n$, it suffices to prove that $D(s)$ has empty interior.
Assume, for a contradiction, that $V$ is a nonempty subset of~$D(s)$, 
and let $x \in V$.
Since $g$ is continuous, we can assume that $d\Pa{g(x'), g(x)} < s/2$ for
every $x' \in V$.
By definition of $f$ and~$g$,
there exists $x' \in V$ such that $d\Pa{f(x'), g(x)} < s/2$. 
Hence, $d\Pa{f(x'), g(x')} < s$, absurd.

Thus, after shrinking~$U$, we can assume that $f = g$.
Fix $x \in U$, and let $y := f(x)$.
Since $g$ is continuous on $U$, after shrinking $U$ we can assume that
$d\Pa{f(x'), y} < r/3$ for every $x' \in U$. 
Since $\Graph(f) \subseteq \cl A$, there exists $(x', y') \in A$, such that
$(x', y')$ is near $(x, y)$; that is, $x' \in U$ and $d(y', y) < r/3$.
Moreover, by definition of~$g$, there exists $y'' \in C(r)_{x'}$, such that
$d \Pa{f(x'), y''} < r/3$.
However, this implies that $d(y', y'') < r$, contradicting the definition
of $C(r)$.

% We prove now ($5_n$).
% Let $U \subseteq \K^n$ be open  and definable,
% and $f: U \to \K$ be definable; we have to prove that $\Dis(f)$ is nowhere
% dense.
% \Wlog, $f$ is bounded; moreover, by ($\ref{en:i-imin-n}_n$),
% it suffices to prove that there exists $V \subseteq U$ with nonempty
% interior, such that $f$ is continuous on~$V$.
% Let $V := \set{x \in U: \cll\Pa{\Graph(f)}_x = \set{f(x)}}$.
% By ($6_n$), $U\setminus V$ is nowhere dense; thus, $V$ has nonempty
% interior, and, since $A$ is bounded, $f$~is continuous on~$V$.

We prove now 
($\enumref{en:i-imin-n}_n \Rightarrow \enumref{en:i-imin-n}_{n+1}$): 
let $A \subseteq \K^{n+1}$ be definable and with
empty interior; we want to show that $A$ is nowhere dense.
If not, let
\[
E' := \set{x \in \K^n: \cll(A)_x \text{ has nonempty interior}}
\]
By assumption, $E'$ has nonempty interior.
\Wlog, $A$ is bounded.
Let $E := \set{x \in \K^n: \cll(A_x) \text{ has nonempty interior}}$.
Since, by $(\ref{en:i-B})_n$, $E \sdiff E'$ is nowhere dense,
$E$~has nonempty interior.
By \enumref{en:i-imin-1}, 
$E = \set{x \in \K^n: A_x \text{ has nonempty interior}}$.
Since $\K$ is definably Baire, there exists $0 < r \in \K$ such that
\[
E(r) := \set{ x \in \K^n: A_x \text{ contains an interval of length } r}
\]
has nonempty interior; let $U \subseteq E(r)$ be a nonempty open set.

%There exists a definable function $h: U \to \K$, such that, for every $x \in
%U$, $h(x)$ is the centre of a 
% For every $x \in U$, let
% \[
% F_x := \bigcup \set{I : I \subseteq \K \text{ open interval}, \abs I \geq r, I
% \subseteq A_x},
% \]
% and $F := \bigsqcup_{x \in U} \Pa{\set x \times F_x} \subseteq A$.
After shrinking $A$ if necessary, we can assume that, for every $x \in U$,
$A$ is bounded and $A_x$ is an open interval of length~$r$.
For every $x \in U$, let $h(x)$ be the center of~$A_x$.

By $\eqref{en:i-f}_n$,
the set of points where $h$ is continuous has nonempty interior.
Thus, after shrinking~$U$, we can assume that $h$ is
continuous.
But then the set $\set{\pair{x, y} \in U \times \K: h(x) - r/2 < y < h(x) + r/2}$ is
open and contained in~$A$, absurd.

($\ref{en:i-B} \Leftrightarrow \ref{en:i-fr}$) is clear.

($\ref{en:i-imin-n} \Rightarrow \ref{en:i-meager-n} \Rightarrow
\ref{en:i-meager-1})$ are also clear.

($\ref{en:i-meager-1} \Rightarrow \ref{en:i-imin-1}$).
Let $X \subseteq \K$ be definable and have empty interior.
By hypothesis, $X$~is meager.
Moreover, $\fr X = \cl X \setminus X$ has also empty interior, and thus it is
meager.
Therefore, $\cl X$ is meager.
Since $\K$ is definably Baire, $\cl X$ has empty interior.

$(\ref{en:i-B} \Rightarrow \ref{en:i-dim-n}$).
Let $X \subseteq \K^n$ be definable, and let $d := \dim \cl X$.
We want to prove that $\dim X = d$.
\Wlog, $\pi(\cl X)$ contains an open subset of~$\K^d$, where $\pi := \Pi^n_d$.
If, for a contradiction, $\dim X < d$, then, by ($\ref{en:i-imin-n}$),
$\pi(X)$~is nowhere dense.
Notice that $\pi(\cl X) \setminus \pi(X) \subseteq \B_d(X)$.
Since $\B_d(X)$ is nowhere dense, we get a contradiction.

($\ref{en:i-imin-n} \Rightarrow \ref{en:i-dim-fiber}$).
Let $X := \set{x \in \K^d: \dim A_x > 0}$, where $d := \dim A$, and
assume, for a contradiction, that $X$ is somewhere dense (and thus, by
\enumref{en:i-imin-n}, $X$~has nonempty interior).
If $d = 0$, then $X = \emptyset \subset \K^0 = \mset{0}$,
and we have a contradiction.
%For every $i = 1, \dotsc, n$, let 
%$X_i := \set{x \in \K^d}$
Thus, \wloG, $A$~is closed (because $\dim \cl A = \dim A)$, 
and $Y := \Pi^n_{d + 1}(A)$ satisfies
\[
\forall x \in X\ \dim(Y_x) > 0.
\]
By Kuratowski-Ulam theorem, this implies that $Y$ is not meager,
and thus has nonempty interior, contradicting $\dim A = d$.

($\ref{en:i-imin-n} \Rightarrow \ref{en:i-dim-union-n}$).
Let $A_1, A_2 \subseteq \K$ be definable, such that $\dim A_i < d$, $i = 1, 2$.
We have to prove that $\dim(A_1 \cup A_2) < d$.
Assume, for a contradiction, that $B := \Pi^n_d(A_1 \cup A_2)$ has nonempty
interior. 
Let $B_i := \Pi^n_d(A_i)$;
notice that $\dim B_i < d$, and $B_1 \cup B_2 = B$.
By (\ref{en:i-imin-n}), the $B_i$ are nowhere dense in~$\K^d$;
thus, $B$~is nowhere dense, absurd.

($\ref{en:i-dim-union-n} \Rightarrow \ref{en:i-dim-union-1}$) is obvious.

($\ref{en:i-dim-union-1} \Rightarrow \ref{en:i-dim-1}$).
Let $A \subseteq \K$ be definable.
We have to prove that $\dim (\cl A) = \dim(A)$.
However, $\cl A = A \cup \fr A$.
Since $\fr A$ has empty interior, $\dim \fr A = 0$, and we are done.

($\ref{en:i-dim-fiber} \Rightarrow \ref{en:i-dim-union-1}$).
Let $A, B$ be definable subsets of $\K$ with empty interior.
We have to prove that $A \cup B$ has also empty interior.
Define $X := \Pa{A \times (0, 1) } \cup \Pa{B \times (2, 3)} \subset \K^2$.
Notice that $\dim X = 1$; thus, by~$\enumref{en:i-dim-fiber}$, the set
$Y := \set{y \in \K: \dim(X_y) > 0}$ is nowhere dense.
However, $Y = A \cup B$; thus, $\dim(A \cup B) = 0$.

($\ref{en:i-imin-n} \Rightarrow \mathrm{\ref{en:i-increasing}}$).
Let $\Pa{A_x}_{x \in \K}$ be an increasing definable family of subsets of~$\K^n$, each of them of dimension less or equal to~$d$.
Let $A := \bigcup_x A_x$.
Assume, for a contradiction, that $\dim A > d$; \wloG, $U := \Pi^n_{d+1}(A)$ has
nonempty interior.
However, $U = \bigcup_x \Pi^n_{d + 1}(A_x)$.
Since $\dim A_x \leq d$, each $\Pi^n_{d + 1}(A_x)$ is nowhere dense, and thus
$U$ is meager, contradicting the fact that $\K$ is definably Baire.

($\ref{en:i-dim-fiber-2} \Rightarrow \ref{en:i-dim-fiber}$) is obvious.

($\ref{en:i-dim-fiber} \Rightarrow \mathrm{\ref{en:i-dim-fiber-2}}$).
Assume, for a contradiction, that $\dim C > k$; \wloG, $U := \Pi^n_{k + 1}(C)$ has
nonempty interior.
Moreover, since, by \enumref{en:i-dim-n},
$\dim(A) = \dim(\cl A)$, %, by property (\ref{en:i-increasing}),
\wloG $A$ is \dcompact.
By~\enumref{en:i-dim-union-n}, \wloG the set
\[
C' := \set{x \in \K^n: \dim \Pa{\Pi^m_d(A_x)} \geq d }
\]
has dimension greater than~$k$, %where $\pi := \Pi^n_d$,
and $D' := \Pi^n_{k + 1}(C')$ has nonempty interior.
%\pi(x_1, \dotsc, x_n) := (x_{n - d + 1}, \dotsc, x_{n} is the projection onto the last 
Let $B := \Pi^{n + m}_{d + k + 1}(A)$; by assumption, $B$~is nowhere dense.
Hence, by Kuratowski-Ulam theorem, the set
\[
D := \set{u \in \K^{k + 1}: \dim(B_u) \geq d}
\]
has empty interior.
However, for every $u \in \K^{k + 1}$, $B_u = \Pi_d(A_u)$, and thus $D'
\subseteq D$, absurd.

($\ref{en:i-imin-1} \Rightarrow \mathrm{\ref{en:i-monotonicity}}$)
%and ($\ref{en:i-imin-1} \Rightarrow \mathrm{\ref{en:i-continuous}}$)
has  the same proof as \cite{miller05}*{Thm.~3.3}.

($\ref{en:i-countable-n} \Rightarrow \mathrm{\ref{en:i-countable-1}}$)
is clear (take $d = n = 1$ in \enumref{en:i-countable-n}).

($\ref{en:i-countable-1} \Rightarrow \mathrm{\ref{en:i-dim-union-1}}$)
is also clear: if $A$ and $B$ are subsets of $\K$ of dimension~$0$,
then the family $\pair{A, B}$ is an at most pseudo-enumerable family of unary
sets with empty interior, and thus their union $A \cup B$ has empty interior.

($\ref{en:i-imin-n} \Rightarrow \mathrm{\ref{en:i-countable-n}}$)
Let $\Pa{X_t: t \in N}$ be a definable family, such that $N$ is
at most pseudo-enumerable, and each $X_t$ is a subset of $\K^n$ of dimension at
most~$d$.
We have to prove that $Y \coloneqq \bigcup_{t \in N} X_t$ has dimension at
most~$d$.
By projecting onto some $\K^{d + 1}$, \wloG we can assume that $n = d + 1$;
by \enumref{en:i-imin-n}, each $X_t$ is nowhere dense, and we have to prove
that $Y$ is nowhere dense.
By definition of pseudo-enumerable, \wloG we can assume that $N$ is a closed,
discrete subset of $\K_{\geq 1}$.
Let $t_0 \in \K$.
Notice that the set $\set{t \in N: t < t_0}$ is pseudo-finite.
Thus, by \cite{fornasiero-lomin}*{Lemma~5.23}, $\bigcup_{t < t_0} X_t$ is
nowhere dense.
Thus, by \enumref{en:i-increasing}, $Y$~is also nowhere dense.

($\ref{en:i-meager-n} \Rightarrow \mathrm{\ref{en:i-ao}}$).
Let $X \subseteq \K^n$ be definable.
Notice that $Y \coloneqq X \setminus \inter X$ has empty interior; thus, it is meager.
Thus, $X = \inter X \cup Y$ is \ao.
\mbox{}\hfill \qedsymbol

\subsection{Constructible structures}\label{sec:ipmin}

\begin{definizione}[see \cite{miller05}*{\S3.2}]
$\K$ is a \intro{\ipminimal{}} structure if every definable subset of $\K^n$
is constructible, for every $n \in \N$.
% , for every $\K' \elem \K$, 
% every $0$-dimensional set definable in $\K'$ is constructible.
A theory $T$~is a \ipminimal  if every model of $T$ is \ipminimal.
\end{definizione}
See also \cite{robinson74} and \cite{pillay87}*{\S2} for related notions.

% \begin{fact}[\cite{pillay87}*{Lemma~2.3}]\label{fact:constructible-nd}
% Let $\K$ be \ipminimal.
% Let $Y \subseteq X \subseteq \K^n$ be a definable set.
% $Y$ is nowhere dense in $X$ iff $Y$ has no interior in~$X$.
% \end{fact}

\begin{thm}\label{thm:constructible}
The following are equivalent:
\begin{enumerate}
%\item\label{en:c-const-1} $\K$ is \ipminimal;
\item\label{en:c-const-n} $\K$ is \ipminimal;
%every definable set is constructible;
\item for every $\emptyset$-definable $A \subseteq \K^{m + n}$ there exists $N
\in \Nat$ such that, for all $x \in \K^m$, if $\dim A_x = 0$, 
then $\inlc{(A_x)}{N} = \emptyset$ \rom(see Definition~\ref{def:lc}\rom);
\item 
every $\emptyset$-definable set is a finite union of
$\emptyset$-definable locally closed sets;
\item\label{en:c-const-1} 
for every $\K' \elem \K$, 
every $0$-dimensional set definable in $\K'$ is constructible.
\item
every definable subset of $\K^n$ is a finite union of sets of the form
\[
\set{x \in \K^m: f(b,x) = 0 \et g(b,x) > 0},
\]
where $f$ and $g$ are $\emptyset$-definable and continuous, and $b \in \K^m$.
\end{enumerate}
Moreover, if $\K$ is \ipminimal, then it is \iminimal.
\end{thm}
\begin{proof}
The equivalence of the first 4 points is proved in the same way
as \cite{miller05}*{Thm.~3.2}, using Lemma~\ref{lem:i-lc}.
$(5 \Rightarrow 4)$ is obvious.
$(4 \Rightarrow 5)$ is proved in the same way as \cite{vdd-dense}*{Lemma~2.10}.
%The ``moreover'' clause is a particular case of Fact~\ref{fact:constructible-nd}.
The ``moreover'' clause follows from the fact that a constructible set 
with empty interior is nowhere dense.
\end{proof}

Notice that the equivalence ($\ref{en:c-const-1}\Leftrightarrow \ref
{en:c-const-n}$) in Theorem~\ref{thm:constructible} shows that, if $\K$ is
\ipminimal and $\K'$ is elementary equivalent
to~$\K$, then $\K'$ is also \ipminimal.

\begin{remark}
\Iminimality is not equivalent to \ipminimality.
In fact,
\begin{enumerate}\item 
it is easy to build an ultra-product of \ipminimal structures which
is not \ipminimal (while an ultra-product of \iminimal structures is
\iminimal);
\item \cite{FKMS}*{Theorem~A} produces an \iminimal structure that defines sets
on every level of the projective hierarchy.
\end{enumerate}
\end{remark}

\subsection{Pillay and Cantor-Bendixson rank in \ipminimal structures}
\label{subsec:rank}
\cite{pillay87}, extending the work in~\cite{robinson74}, studies topological
structures satisfying a weaker version of \ipminimality; that is, Pillay's
Condition~(A) asks that every definable unary set is constructible.

\begin{proviso}
In this subsection we assume that $\K$ is \ipminimal.%
\footnote{Some of the results in this subsection hold without this assumption.}
\end{proviso}

Pillay then defines a \intro{rank} for  definable sets (which he calls
the dimension rank), which we will denote by~$\rkP$,  in the following way:
\begin{enumerate}
\item If $X$ is nonempty, then $\rkP(X) \geq 0$.
\item $\rkP(X) \geq \lambda$ iff $\rkP(X) \geq \alpha$ for all 
$\alpha < \lambda$, where $\lambda$ is limit.
\item $\rkP(X) \geq \alpha + 1$ iff $X$ contains subset $Y$ which is
definable, closed,
nowhere dense (in $X$), and with $\rkP(Y) \geq \alpha$.
\end{enumerate}
Notice that $\rkP(X)$ might depend on the ambient structure~$\K$; that is,
if $\K' \elem \K$, then $\rkP(X^{\K'})$ might be different from
$\rkP(X^{\K})$.

%Let $M$ be a Hausdorff topological structure, such that every $M'$ elementarily
%equivalent to $M$ satisfies Condition (A) (\eg, $M$ is \ipminimal),
%Let $\K$ be an \ipminimal  structure and $X \subseteq \K^n$ be definable.
%and $X \subseteq M$ be definable and closed.

\begin{fact}[Pillay]\label{fact:rkP-closed}
Let $X \subseteq \K^n$ be definable.
%\footnote{Do we need $\K$ constructible here?}
\begin{enumerate}
\item $\rkP(X) = 0$ iff $X$ is discrete and nonempty;
\item $Y \subseteq X \Rightarrow \rkP(Y) \leq \rkP(X)$;
\item If $X = X_1 \cup \dots X_n$, where the $X_i$ are closed (in~$X$)
and definable,
then $\rkP(X) = \max_{1 \leq i \leq n} \rkP(X_i)$.
\end{enumerate}
\end{fact}

\begin{lemma}\label{lem:rkP-union}
Let $X \subseteq \K^n$ be definable.
Assume that $X = A  \cup B$, where $A$ and $B$ are definable, and
$A$ is open in $X$.
Then, $\rkP(X) \leq \rkP(B) + \rkP(A) + 1$ (where $+$ is the usual ordinal sum).
\end{lemma}
\begin{proof}
By induction on $\rkP(A)$.
If $A$ is empty, the conclusion is clear; thus, we can assume that $A$ is nonempty.
Assume, for a contradiction, that $Y \subseteq X$ is definable, closed, and
nowhere dense in $X$, 
but $\rkP(Y) \geq \rkP(B) + \rkP(A) + 1$.
Let $Y_A := Y \cap A$ and $Y_B := Y \cap B$.
Since $A$ is open, $Y_A$ is nowhere dense in~$A$; therefore (since $A$ is
nonempty) $\rkP(Y_A) < \rkP(A)$.
Moreover, $Y_A$ is open in~$Y$.
Thus, by inductive hypothesis, 
\[
\rkP(Y) \leq \rkP(Y_B) + \rkP(Y_A) + 1 \leq \rkP(B) + \rkP(Y_A) + 1
< \rkP(B) + \rkP(A) + 1,
\]
absurd.
\end{proof}

\begin{proposition}
%Assume that $\K$ is constructible.
Let $X \subseteq \K^n$ be definable.
Write $X = X_1 \cup \dots \cup X_r$, where each $X_i$ is locally closed
in~$X$.
Let $\gamma := \max_i \rkP(X_i)$.
Then,
\[
\rkP(X) \leq r \gamma + (r - 1).
\]
\end{proposition}\label{prop:rkP-lclosed}
\begin{proof}
By induction on~$r$.
If $r = 1$, the result is clear.
Thus, we can assume that $r \geq 2$ and that 
we have already proved the result for $r - 1$.
Fix $i \leq r$.
Let $A_i$ be the closure of $X_i$ inside~$X$.
\begin{claim}
$\rkP(A_i) \leq r \gamma + (r - 1)$.
\end{claim}
Since $X = A_1 \cup \dots \cup A_r$ and each $A_i$ is closed in $X$, the
conclusion then follows from Fact~\ref{fact:rkP-closed}.

Thus, it suffices to prove the claim; \wloG, we can assume $i = 1$.
Let $Y := (X_2 \cap A_1) \cup \dots (X_r \cup A_1)$.
Since $X_j \cap A_1 \subseteq X_j$, we have $\rkP(X_j) \leq \gamma$ for every
$j = 2, \dotsc, n$.
Therefore, by inductive hypothesis,
$\rkP(Y) \leq (r - 1)\gamma + r - 2$.
Moreover, $A_1 = X_1 \cap Y$, and $X_1$ is open in~$A_1$ (because $X_1$ is
locally closed).
Therefore, by Lemma~\ref{lem:rkP-union},
$\rkP(A_1) \leq \rkP(Y) + \rkP(X_1) + 1 \leq (r-1) \gamma + (r - 2) + \gamma +
1 \leq r \gamma + (r-1)$.
\end{proof}

\begin{fact}[Pillay]
%Assume that $\K$ is constructible.%
%\footnote{Check if the assumption is needed.}
\Tfae:
\begin{enumerate}
\item $\rkP(X) = \infty$;
\item 
there is a decreasing sequence $(X_i)_{i < \omega}$ of definable closed
subsets of $X$, such that $X_{i + 1}$ is definable, closed, and nowhere dense
in $X_i$, for all $i < \omega$;
\item
$X$ contains a definable closed nowhere dense subset~$Y$,
such that $\rkP(Y) = \infty$;
\item
$\rkP(X) > 2^{\card M}$.
\end{enumerate}
\end{fact}

\begin{lemma}\label{lem:constructible-isolated}
%Assume that $\K$ satisfies Condition~\rom(A\rom).
%Let $\K$ be \ipminimal, and 
Let $C \subseteq \K$ be nonempty, definable, closed,
and with empty interior.
Then, $C$~has at least one isolated point. 
In particular,
no nonempty closed perfect subset of $\K$ with empty interior is definable.
\end{lemma}
\begin{proof}
Assume, for a contradiction, that $C$ is perfect.
Let $A := \K \setminus C$; $A$ is an open set;
let $C^L$ be the set of left end-points of the connected components of~$A$
and $C^R$ be the set of right end-points.
Notice that $C^L$ and $C^R$ are definable subsets of~$C$.
Since $C$ is perfect and with empty interior, $C^L$ and $C^R$ have no isolated
points, and they are both dense in~$C$.
Notice that $C^L \cap C^R$ is the set of isolated points of~$C$:
since $C$ is perfect, $C^L$ and $C^R$ are disjoint.
Let $D \coloneqq \lc(C^L)$.
Since $\K$ is \ipminimal, 
%By Condition~(A), 
$D$ is dense in $C^L$; and therefore $D$ is dense in~$C$.
Let $a \in D$; since $D$ is locally closed, there exists $b', b'' \in \K$,
such that $b' < a < b''$ and $I \cap D$ is closed in~$I$, where $I \coloneqq
[b', b'']$.
Thus, $D$ is dense in $C$, $I \cap D = I \cap C$: thus,
$I \cap C \subseteq C^L$, which contradicts the fact that $C^R$ is dense 
in~$C$.
% there exists an open interval $I$ such that
% $I \cap C^L$ and $I \cap C^R$ are closed in~$I$ 
% and $I \cap C^L \neq \emptyset$.
% Let $a \in C^L \cap I$.
% Since $C$ is perfect and has empty interior, 
% $a$~is an accumulation point for~$C^R$; hence, $a \in C^L \cap C^R$,
% implying that $a$ is isolated in~$C$, absurd.
\end{proof}

\begin{lemma}\label{lem:constructible-derived}
%Let $\K$ be \ipminimal
Let $C \subset \K$ be definable with empty interior and
$D$ be the set of isolated points of~$C$.
Then, $D$~is discrete, definable, and dense in~$C$.
Moreover, $C' := C \setminus D$ is nowhere dense in~$C$.
\end{lemma}
\begin{proof}
That $D$ is discrete and definable is clear.
% Since $C$ is constructible, $\lc(C)$ is dense in~$C$.
% \begin{claim}
% It suffice to prove the conclusion for $\lc(C)$.
% \end{claim}
% In fact, the isolated points of $\lc(C)$ remain isolated in~$C$.
% Thus, \wloG $C$ is locally closed.
\begin{claim}
It suffices to prove the conclusion for $\cl C$.
\end{claim}
In fact, the isolated points of $\cl C$ and the isolated points of~$C$
coincide.

Thus, \wloG $C$ is closed.
Let $A := \K \setminus C$, $C^L$ be the set of left end-points of the
connected components of~$A$ and $C^R$ be the set of right end-points.
$C^L$ is dense in~$C$.
If, for a contradiction, $D$ is not dense in~$C$, let $I$ be a closed interval,
such that $C \cap D$ has no isolated points and is nonempty: but this
contradicts the Lemma~\ref{lem:constructible-isolated}.

The fact that $C'$ is nowhere dense in $C$ follows immediately from the first part.
\end{proof}

\begin{definizione}[Cantor-Bendixson Rank, see \cite{kechris}*{\S I.6.c}]
\label{def:CB}
Let $T$ be a Hausdorff topological space.
For every $X \subseteq T$, and every ordinal $\alpha$,
let $\CBd X 0 := X$,
$\CBd X {\alpha + 1}$ be the set of non-isolated points of~$\CBd X {\alpha}$,
and $\CBd X {\alpha} := \bigcap_{\beta < \alpha} \CBd X {\alpha}$ if $\alpha$
is a limit ordinal.
Let $\rkCB(X)$, the Cantor-Bendixson rank of~$X$,
be the smallest ordinal $\alpha$ such that
$\CBd X {\alpha} = \emptyset$ (or $\rkCB(X) = + \infty$ if such $\alpha$ does
not exist).
For every $a \in X$, let $\rkCB_X(a)$ be the supremum of ordinals $\alpha$ such
that $a \in \CBd X {\alpha}$.
\end{definizione}
Notice that each $\CBd X {\alpha} \setminus \CBd X {\alpha + 1}$ is discrete,
and that $X$ is a finite union of discrete sets iff $\rkCB(X) < \omega$.
Moreover,  $\CBd X {\rkCB(X)} = \emptyset$ (if $\rkCB (X) < + \infty)$.
Besides, $\rkCB(X) = 0$ iff $X$ is empty, $\rkCB(X) = 1$ iff $X$ is discrete
(and nonempty).

%Remember the definition of Cantor-Bendixson derivative and rank
\begin{remark}\label{rem:CB-union}%
\footnote{The present remark is folklore, but we could not find a reference for it. }
If $X$ and $X'$ are subsets of a Hausdorff topological space~$T$, 
then $\rkCB(X \cup X') \leq \rkCB(X) \oplus \rkCB(X')$, where $\oplus$ is the
Cantor sum of ordinals.
Hence, a set $X$ is a union of $n$ discrete sets iff $\rkCB(X) \leq n$.
\end{remark}
\begin{proof}
By induction on $\alpha := \rkCB(X)$ and $\beta :=\rkCB(Y)$.
The basic case when $X$ is a singleton is obvious.
% Choose $\alpha' < \alpha$ and $\beta' < \beta$.
% Let $X' := X^{[\alpha']}$, $Y' := Y^{[\beta']}$.
% By inductive hypothesis, $\rk(X' \cup Y') \leq \alpha' \oplus \beta' < \alpha
% + \beta$.
For a contradiction, let $a \in \CBd{(X \cup Y)} {\alpha + \beta}$.
\Wlog, $a \in X$; let $\gamma := \rkCB_X(a) < \alpha$.
%Similarly, if $a \in Y$, let $\beta' := \rkCB_Y(A)$.

Let $V$ be an open neighborhood of $a$ such that
$\CBd X{\gamma} \cap V = \mset a$, and $X' := X \cap V \setminus \mset a$.
Notice that $\alpha' := \rkCB(X') \leq \gamma < \alpha$.
Hence, by inductive hypothesis, $\rkCB(X' \cup Y) \leq \alpha' \oplus \beta <
\alpha \oplus \beta$.
Thus, $\rkCB(X' \cup Y \cup \mset a) \leq \alpha' \oplus \beta \oplus 1 \leq
\alpha \oplus \beta$.
Thus, $\rkCB_{X \cup Y}(a) < \alpha \oplus \beta$ for every $a \in X \cup Y$,
and we are done.
\end{proof}

\begin{exercise}
Given $N \in \Real$, find $A_1, \dotsc, A_N$ discrete subsets of $\Real$,
such that $\rkCB(A_1 \cup \dots \cup A_N) = N$.
%find a $A \subset \Real$ and $b \in \Real$,
%such that $\rkCB(A) = N$ 
%and $\rkCB(A \cup \mset b) = N + 1$.
\end{exercise}

From the above results, it is easy to deduce the following.
\begin{proposition}\label{prop:rkP-unary}
Assume that $\K$ is \ipminimal. 
\begin{enumerate}
\item 
Let $X \subseteq \K$ be definable and closed, such that
$X$ is a finite union of discrete sets. 
Then, $\rkCB(X) = \rkP(X)$.
\item $\K$ is locally o-minimal iff $\rkP(\K) = 1$.  
\item If $\K$ is not locally o-minimal, then $\rkP(\K) \geq \omega$.
\item If $\K$ is \dminimal but not locally o-minimal, 
then $\rkP(\K) = \omega$.
\item If $\K$ is $\omega$-saturated,  \ipminimal, but not \dminimal, 
then $\rkP(\K) = \infty$.
\end{enumerate}
\end{proposition}
\begin{proof}
The last point follows from the fact that, since $\K$ is not \dminimal, then,
by saturation, we can find $X \subset \K$ definable, closed, with empty
interior, and such that $\rkCB(X) \geq \omega$.
Hence, by Lemma~\ref{lem:constructible-derived}, 
$X \supset X^{(1)} \supset X^{(2)} \supset \dots$ is an infinite descending 
chain of definable sets, such that 
$X^{(i + 1)}$ is closed and nowhere dense in $X^{(i)}$.
\end{proof}
%We do not know what $\rkP(\K)$ is when $\K$ is \ipminimal.
%In particular, we do not know whether it is possible that $\rkP(\K) = \infty$.

We will compute later $\rkP(\K^n)$.

% ----------------------------------------------------

\section{Definable choice}\label{sec:DSF}
As usual, $\K$~is some definably complete expansion of an ordered field.

\begin{lemma}\label{lem:open-closed}
Let $P \subseteq \K$ be a set of parameters, and
$U \subseteq \K^n$ be $P$-definable and open.
Then, there exists $C \subseteq \K^{n+1}$ $P$-definable and closed,
such that $U = \Pi^{n+1}_n(C)$.
\end{lemma}
\begin{proof}
  If $U = \K^n$, take $C := \K^{n+1}$.
Otherwise, for every $r > 0$, let
\[\begin{aligned}
U(r) &:= \set{x \in U:  d(U, \K^n \setminus U) \geq r },\\
D &:= \bigcup_{r > 0} U(r) \times \set{1/r} \subseteq \K^n \times
\K_{>0},\\
C &:= \cl D.
\end{aligned}
\]
By definition, $C$ is closed, it is trivial that $C$ is $P$-definable,
and it is easy to see that $\Pi^{n+1}_n(C) = U$.
\end{proof}

%Some version of the following lemma has been proved by C.~Miller
\begin{lemma}[Definable Choice]\label{lem:skolem}
\begin{enumerate}
\item Let $P \subseteq \K$ be a set of parameters, and
$A \subset \K^n$ be $P$-definable, nonempty and constructible.
Then, there exists a $P$-definable point $a \in A$ .
\item Let $X \subseteq \K^{m + n}$ be definable and such that $X_b$ is
constructible for every $b \in \K^m$.
Then, $X$ has a definable $n$-choice function, that is a definable function 
 $f : \Pi^{m + n}_m(X) \to X$, such that $f(a) \in \mset a \times X_{a}$ for
 every $a \in \Pi^{m + n}_m(X)$.
 \item Let $X \subseteq \K^{n + m}$ be definable and $\Fs$.
Then, $X$ has a definable $n$-choice function.
\item Suppose that every unary definable set contains a locally closed point.
Then, $\K$~has definable Skolem functions \rom(\intro{\DSF{}}\rom).
\end{enumerate}
\end{lemma}
\begin{proof}
($1$).
Since $A$ is constructible, $A' := \lc(A)$ is nonempty; thus, since $A'$ is
also $P$-definable, it suffices to prove the conclusion for~$A'$;
therefore, \wloG $A$ is locally closed.

%First, we will prove the case when 
\Case 1 $A$ is closed in~$\K^n$.
For every $r > 0$, let $A(r) := \set{a \in A: \abs a \leq r}$, 
let $r_0 := \inf\set{r \in \K: A(r) \neq \emptyset}$, and let $A' := A(2r_0)$.
Notice that $A'$ is \dcompact, nonempty and $P$-definable.
Let $a := \lexmin(A')$: notice that $a$ is also $P$-definable, and in~$A$.

\smallskip

In the general case , since $A$ is locally closed, we can write
$A = U' \cap cl A$  for some open set~$U'$.

\begin{claim}
There exists $U \subseteq \K^n$ open and $P$-definable such that
$A = U \cap \cl A$.
\end{claim}
For every $a \in A$ let 
$r(a) := \sup \set{r > 0: A \cap B(a,r) = \cl A \cap B(a,r)}$.
Since $A$ is locally closed, $r(a) > 0$ for every $a \in A$.
Let $U := \bigcup_{a \in A} B\Pa{a, r(a)/2}$.

Thus, by Lemma~\ref{lem:open-closed}, there exists 
$D \subseteq \K^{n+1}$ closed and $P$-definable such that $U =
\Pi^{n+1}_n(D)$.
Let $E := D \cap (\Pi^{n+1}_n)^{-1}(A)$.

\begin{claim}
$E$ is closed (inside~$\K^{n+1)}$, nonempty, $P$-definable, and $\Pi^{n+1}_n(E) = A$.
\end{claim}
In fact, since $A$ is closed in~$U$ and $\Pi^{n+1}_n$ is continuous,
$E$ is closed in~$D$; since $D$ is closed in $\K^{n+1}$, $E$ is also
closed in $\K^{n+1}$. The remainder of the claim is clear.

Thus, by the case when $A$ is closed, we can find $e \in E$ which is
$P$-definable, and let $a := \Pi^{n+1}_n(e)$.

% \Case 2 $A$ is discrete.
% For every $r > 0$, let $A(r) := \set{a \in A: B(a,r) \cap A = \mset a}$, 
% let $r_0 := \sup\set{r \in \K: A(r) \neq \emptyset}$, 
% and let $A' := A(r_0/2)$.
% Notice that $A'$ is closed, nonempty, and $P$-definable, and apply Case~1.

% \Case 3 $A$ is open.
% For every $r > 0$, let 
% $A(r) \coloneqq \set{a \in A: d(a, \K^n\setminus A) \geq r}$,
% and $r_0 \coloneqq \sup\set{r \in \K: A(r) \neq \emptyset}$.
% Notice that $A(r_0/2)$ is $P$-definable, closed in $\K^n$, nonempty, and
% contained in~$A$; apply Case~1.

% We will prove the general case by induction on~$n$.
% If $n = 1$, for every $a \in A$ let 
% $r(a) := \sup \set{r > 0: A \cap B(a,r) = \cl A \cap B(a,r)}$.
% Since $A$ is locally closed, $r(a) > 0$ for every $a \in A$.
% Let $U := \bigcup_{a \in A} B\Pa{a, r(a)/2}$.
% Notice that $U$ is open, $P$-definable, and $A = \cl A \cap U$.
% By Case~3, there exists $u \in U$ $P$-definable.
% Let $I$ be the definably connected component of $U$ containing $u$:
% $I$ exists because $n = 1$, and it is a $P$-definable open interval;
% let $A' := A \cap I$.
% Therefore, $A'$ is closed in $I$, $P$-definable and nonempty;
% since $I$ is definably homeomorphic to $\K$, via a $P$-definable map, we can
% apply Case~1 to find a $P$-definable point in~$A'$.

% If $n > 1$, we proceed by induction, and we assume we have already proved the
% conclusion for $n - 1$.
% Let $A' := \Pi^{n}_{n - 1}(A)$.
% By the inductive hypothesis, there exists $a' \in A'$ $P$-definable.
% By the case $n = 1$, there exists $a'' \in A_{a'}$ also $P$-definable.
% Let $a := \pair{a', a''}$.

\smallskip

($2$).
The construction in ($1$) gives a definable way to choose $x_b \in X_b$ for
every $b \in \Pi^{n + m}_m(X)$.
(Equivalently: ($2$) follows from ($1$) by compactness).

\smallskip

($3$).
If $\K$ is locally o-minimal, the result is clear.
Otherwise, there exists $N \subset \K_{\geq 0}$ which is definable without
parameters,  closed, discrete, and unbounded.
Since $X$ is $\Fs$, there exists a definable family of \dcompact sets
$\Pa{X(i): i \in N}$ such that $X = \bigcup_{i \in N} X(i)$.
For every $i \in N$, let $Y \coloneqq \Pi^{m+n}_m(X(i))$,
and $Z(i) \coloneqq Y(i) \setminus \bigcup_{j < i} X(j)$.
Notice that $Y \coloneqq \Pi^{m + n}_m(X)$ is the disjoint union of all $Z(i)$.
For every $\z \in Z(i)$, let $f_i(\z) \coloneqq \pair{z, \lexmin(X(i)_\z)}$.
%\mset z \times\set{\lexmin(X(i)_\z)}$.
Let $f: Y \to X$ be the function that coincides with $f_i$ on each $Z(i)$.
Then, $f$~is a definable $n$-choice function for~$X$.

$(4)$.
Let $A \subseteq \K^{m + n}$ be definable.
We have to prove that there exists $f: B \to A$ definable, such that $f(b) \in
\mset b \times {A_b}$ for every $b \in B$, where $B := \Pi^{n + m}_m(A)$.
We proceed by induction on~$n$.

If $n = 0$, $f$ is the identity.
If $n = 1$, for every $b \in B$ let $C_b := \lc(A_b)$.
By hypothesis, $C_b \neq \emptyset$ for every $b \in B$, and it is
constructible. Thus, by ($2$), $C$ has a definable 1-choice function, and the
same function will work for~$A$.

Assume that $n > 1$ and we have already proved the conclusion for every
$n' < n$.
Let $C := \Pi^{m + n}_{m + 1}(A)$.
By inductive hypothesis, there exists a definable $(n - 1)$-choice function
$g: C \to A$.
By the case $n = 1$, there exists a definable 1-choice function
$h: \Pi^{m + 1}_m(C) \to C$.
Let $f := g \circ h : \Pi^{m + n}_m(A) \to A$: $f$ is an n-choice function
for~$A$.
\end{proof}

Therefore, locally o-minimal, d-minimal and \ipminimal structures
%(see \S\ref{sec:ipmin}) 
have \DSF.
On the other hand, 
structures with locally o-minimal open core might not have \DSF:
for instance, $\pair{\Real, \Ralg}$ does not have \DSF (see \cite{DMS}*{5.4}).
For the same reason, \ipminimal structures have elimination of imaginaries.
\cite{miller-definable-choice} already proved DSF and elimination
of imaginaries  for \dminimal structures.
%\footnote{C.~Miller ``Definable choice in d-minimal
%expansions of ordered groups''. Unpublished note, December 2006.}

%\emph{What about elimination of imaginaries?}

\begin{lemma}
Let $X \subseteq \K^n$ be definable.
Assume that $X$ is both an $\Fs$ and a $\Gd$, and that $\cl X$ is definably
Baire.
Then, $\lc(X) \neq \emptyset$.
\end{lemma}
\Cf \cite{kuratowski}*{\S34.VI}.
\begin{proof}
Let $Y := \cl X$.
We have to prove that the interior of $X$ inside $Y$ is nonempty.
Otherwise, $X$ is both dense and co-dense in~$Y$.
However, since $X$ is an~$\Fs$, this implies that $X$ is meager in~$Y$.
For the same reason, $Y \setminus X$ is meager in~$Y$,
contradicting the fact that $Y$ is definably Baire.
\end{proof}

\subsection{Sard's Lemma}

For every $\Cone$ function $f : \K^m \to \K^n$, define
$\Lambda_f(k) := \set{\x \in \K^m: \rkM \Pa{Df(\x))} \leq k}$, and
$\Sigma_f(k) := f\Pa{\Lambda_f(k)}$.
The set of singular values of $f$ is 
$\Sigma_f := \bigcup_{k = 0}^{n - 1} \Sigma_f(k)$.

\begin{lemma}[Sard's Lemma]
Let $f : \K^m \to \K^n$ be definable and~$\Cone$.
If $\K$ is \iminimal %and has \DSF (and, in particular, if $\K$ is \ipminimal),
then $\dim\Pa{\Sigma_f(d)} \leq d$.
\end{lemma}
\begin{proof}
If $d \geq n$, the conclusion is trivial.
Let $d < n$, and assume, for a contradiction, that
$\Pi^n_d\Pa{\Sigma_f(d)}$ contains a nonempty open box~$B$.
Since $\Lambda_f(d)$ is an $\Fs$ set,
there exists $g: B \to \Lambda_f(d)$, such that
$\Pi^n_d \circ f \circ g = \identity_B$.
Since $\K$ is \iminimal, we can apply Lemma~\ref{lem:imin-C1}, and,
by shrinking $B$ if necessary, we can assume that $g$ is~$\Cone$.
Hence, by differentiation, we have that, for every $x \in B$,
$\Pi^n_d(f(g(x))) \cdot Df(g(x)) \cdot Dg(x) = \idmatrix_d$.
Therefore, the matrix $Df(g(x))$ has rank at least $d$, a contradiction.
\end{proof}

\subsection{Dimension}
\begin{proviso}
In this subsection, $\K$ is \iminimal with \DSF.
\end{proviso}

\begin{lemma}\label{lem:imin-dimension}
The dimension $\dim$ is a dimension function in the sense of \cite{dries89}.
That is, $\dim$ satisfies the following axioms:
for every definable sets $A$ and $B \subseteq \K^n$ and $C \subseteq \K^{n+1}$,
\begin{enumerate}[\rm({Dim} 1)]
\item $\dim(A) = - \infty$ iff $A = \emptyset$,
$\dim(\mset{a}) = 0$ for each $a \in \K$,
$\dim(\K) = 1$.
\item $\dim(A \cup B) = \max(\dim(A), \dim(B))$;
\item $\dim(A^\sigma) = \dim(A)$ for each permutation $\sigma$ of
$\set{1, \dotsc, n}$.
\item Define $C(i) \coloneqq \set{\x \in \K^n: \dim(C_\x) = i}$, $i = 0, 1$.
Then, each $C(i)$ is definable, and
$\dim\Pa{C \cap (C(i) \times \K)} = \dim(C(i)) + i$, $i = 0,1$.
\end{enumerate}
\end{lemma}
\begin{proof}
The axioms (Dim 1,2,3) are either trivial, or follow immediately from
Theorem~\ref{thm:imin}.
It remains to prove Axiom (Dim 4).
We prove (Dim 4) by induction on $n$.
The fact that the $C(i)$ are definable is clear.
Let $D \coloneqq \Pi^{n+1}_n(C)$.
\Wlog, we can assume that either $D = C(0)$ or $D = C(1)$.
If $n = 0$  the result is clear.
Assume now that $n=1$.

1) If $D = C(0)$, we have to prove that $\dim(C) = \dim(D)$.
Assume not.
If $\dim(C) = 2$, then $C$ has nonempty interior; thus, there exists $d \in D$
such that $C_d$ has nonempty interior, contradicting $\dim(C_d) = 0$.
Since $\dim(C) \geq \dim(D)$, the only possibility left is that
$\dim(C) = 1$ and $\dim(D) = 0$.
Let $\rho: \K^2 \to \K$ be the projection onto the second coordinate;
our assumptions imply that $\rho(C)$ has nonempty interior; let $J$ be an open
interval contained in $\rho(C)$.

By definable choice, there exists a definable function
$f: J \to D$, such that, for every $j \in J$, $\pair{j, f(j)} \in C$.
After shrinking $J$ if necessary, \wloG $f$ is continuous and either constant,
or strictly monotone.
If $f$ is strictly monotone, then $D$ contains an open interval, contradicting
$\dim(D) = 0$.
If $f$ is constant, say $f = d$, then $J \subseteq C_d$, contradicting
$\dim(C_d) = 0$.

2) If $D = C(1)$, we have to prove that $\dim(C) = \dim(D) + 1$.
Assume not.
Remember that $\dim(C) \geq \dim(D)$.
If $\dim(D) = 0$, then $0 \leq \dim(C) \leq 1$, thus $\dim(C) = 0$,
contradicting the fact that $\dim(C_d) = 1$ for every $d \in D$.
If $\dim(D) = 1$, then the only possibility is that $\dim(C) = 1$.
However, Theorem~\ref{thm:imin}(13) implies that $\dim(C_d) = 0$ for at least
one $d \in D$, absurd.

Assume now that we have proved Axiom (Dim 4) for every $n' < n$; we want to
prove it for~$n$. Let $d \coloneqq \dim(D)$.

1) If $D = C(0)$, we have to prove that $\dim(C) = d$.
Assume not.
Since $C \subseteq D \times \K$ and $\dim(C) \geq d$,
the only possibility is that $\dim(C) = d + 1$.
Let $L \coloneqq \K^d \times {\mset 0}^{n-d} \times \K \subset \K^{n+1}$,
and $\Pi^{n+1}_L: \K^{n+1} \to L$ be the orthogonal projection onto~$L$.
Since $\dim(C) = d+1$, after a permutation of the first $n$ coordinates,
we can assume that $E$ has nonempty interior, 
where $E \coloneqq \Pi^{n+1}_L(C)$.
After shrinking $C$ if necessary, we can assume that $E$ is open.
Let $F \coloneqq \Pi^{n+1}_d(C) = \Pi^{L}_D(E)$.
Notice that $F$ has nonempty interior in $\K^d$.
Since $\dim(D) = d$, by Theorem~\ref{thm:imin}(13), for $\fv$ outside a
nowhere dense set, $\dim(D_\fv) \leq 0$.
Thus, there exists $\fv \in F$ such that $\dim(D_\fv) = 0$. 
We will to compute the dimension of $C' \coloneqq C_\fv$ in two different
ways, and reach a contradiction.
Since $E$ is open, $\dim(E_\fv) = 1$, and thus $\dim(C') \geq 1$.
For every $\x \in D_\fv$, $C'_\dv = C_{\fv, \x}$,
and thus $\dim(C'_\dv) = 0$.
Thus, by the case $n = 1$, $\dim(C') = 0$, absurd.

2) If $D = C(1)$, we have to prove that $\dim(C) = d + 1$.
Assume not.
The only possibility is that $\dim(C) = d$.
First, we will treat the case when $d = n$.
Thus, $D$ has nonempty interior.
On the other hand, $C$ is nowhere dense; thus, by Kuratowski-Ulam Theorem,
there exists $\dv \in D$ such that $C_\dv$ has empty interior, contradicting
the fact that $\dim(C_\dv) = 1$.

Assume now that $d < n$.
\Wlog, $D'$ has nonempty interior, where $D' \coloneqq \Pi^n_d(D)$.
Let $L \coloneqq \K^d \times {\set 0}^{n-d} \times \K$ and
$C' \coloneqq \Pi^{n+1}_L(C)$.
By inductive hypothesis, $\dim C \geq \dim(C') = \dim(D') + 1 = d + 1$, absurd.
\end{proof}

\begin{lemma}\label{lem:dim-function}
%Let $\K$ be \iminimal with \DSF, 
Let $X \subseteq \K^n$, $Y \subseteq \K^m$, and
$f: X \to Y$ be definable.
Then:
\begin{enumerate}
\item if $f$ is surjective, then $\dim Y \leq \dim X$;
\item if $f$ is injective, then $\dim Y \geq \dim X$;
\item if $f$ is bijective, then $\dim Y = \dim X$.
\end{enumerate}
\end{lemma}
\begin{proof}
1) Apply Lemma~\ref{lem:imin-dimension} to the graph of~$f$.

2) follows from 1) and definable choice;  3) follows from 1) and 2).
\end{proof}

% ------------

% ----------------------------------------------------

\section{D-minimal structures}\label{sec:dmin}
\subsection{Fundamental results}
\begin{definizione}\label{def:dmin-true}
$\K$ is d-minimal if (it is definably complete and)
 for every $\K' \elem \K$, every definable subset of~$\K'$
is the union of an open set and finitely many discrete sets.
\end{definizione}

\begin{remark}\label{rem:discrete}
Let $N \in \Nat$.
If $A$ is a Hausdorff topological space, such that $A$ is a union of $N$
discrete sets, then $\rkCB(A) \leq N$ \rom(see Definition~\ref{def:CB}\rom).
If $A \subseteq \K^n$ is definable and $\rkCB(A) \leq N$,
then $A$ is a union of $N$ disjoint definable and discrete sets.
%\footnote{Check!}
\end{remark}
\begin{proof}
The first part is immediate from Remark~\ref{rem:CB-union} and induction
on~$N$.
Let $A$ be as in the second part, and proceed by induction on~$N$.
If $N = 1$, $A$ itself is discrete, and we are done.
Assume that we have already proved the conclusion for $N-1$.
Let $B$ be the set of isolated points of~$A$.
Notice that $A$ is the disjoint union of $B$ and $\CBd A 1$, and that
$\rkCB(\CBd A 1) = N-1$.
Therefore, by inductive hypothesis, $\CBd A 1$ is the disjoint union of $N-1$
definable discrete sets, and we are done.
\end{proof}

\begin{remark}
Definition~\ref{def:dmin-true}
is equivalent to Definition~\ref{def:dmin-intro}.
\end{remark}
\begin{proof}
By compactness. 
More in details, assume that $\K$ does satisfy
Definition~\ref{def:dmin-true}.
%~\ref{def:dmin-intro}.
Let $X \subset \K^{n+1}$ be definable; 
we have to prove that there exists $N \in \Nat$ such that,
for each $\bv \in \K^n$, the set $Y_\bv \coloneqq X_\bv \setminus \interior
(X_\bv)$ is the union of at most $N$ discrete sets.
Assume not. 
Let $\K' \succeq \K$ be $\omega$-saturated.
By saturation, there exists $\bv \in {\K'}^n$ such that $Y_\bv$ is not the
union of any finite number of definable discrete sets, and thus, by the 
Remark~\ref{rem:discrete}, $Y_\bv$ is not the union of any finite number of
discrete sets, absurd.
%Thus, $\K$~does not satisfy Definition~\ref{def:dmin-true}.

Conversely, assume that $\K$ satisfies Definition~\ref{def:dmin-intro}.
Let $X \subseteq \K^{n+1}$ be $\emptyset$-definable and let $\K' \elem \K$.
We want to show that, for every $\bv \in {\K'}^n$, the set 
$Y_\bv \coloneqq X_\bv \setminus \interior(X_\bv)$ is the union of finitely
many discrete sets.
By our assumption and Remark~\ref{rem:discrete}, there exists $N \in \Nat$,
such that, for every $\bv \in \K^n$, $\rkCB(Y_\bv) \leq N$.
Since the above can be expressed by a first-order formula,
we have that, for every $\bv \in {\K'}^n$, $\rkCB(Y_\bv) \leq N$,
which (again, by Remark~\ref{rem:discrete}) is what we wanted.
\end{proof}

\begin{remark}\label{rem:dmin-constructible}
Let $\K$ be d-minimal.
Then, $\K$ is \ipminimal, and therefore 
every $\K$\hyph definable set is constructible.
\end{remark}

\begin{question}
Assume that $\K$ is constructible.
Is it \dminimal?
What if moreover $\K$ expands~$\R$?
\end{question}

\begin{definizione}[\cite{miller05}*{Def.~8.4}] %and\S3.4}
\label{def:regular}
Let $d \leq n \in \Nat$, $\Pi(n,d)$ be the set of projections form $\K^n$ onto
$d$-dimensional coordinate spaces, and $\mu \in \Pi(n,d)$.  
Given $A \subseteq \K^n$ definable and  $p \in \Nat$, 
let $\reg^p_\mu(A)$ denote the set of all $a \in A$ such that, for definable
some open neighbourhood~$U$ of~$a$, $A' := A \cap U$ is a $\Cp$ embedded
manifold of dimension~$d$, 
and $\mu \rest {A'}$ maps $A'$ 
$\Cp$-diffeomorphically onto an open subset of~$\K^d$.
$A$ is \intro{$\mu^p$-regular} if it is equal to $\reg^p_\mu(A)$,
and it is $d^p$-regular if it is $\mu^p$-regular for some $\mu \in \Pi(n,d)$
(and will drop the superscript $p$ if it is clear from the context).
We define $\reg^p(A) := \bigcup_{i=0}^n \reg^p_i(A)$.
%If we don't specify the $d$, we take $d = \dim(A)$.

Notice that $\reg^p_0(A)$ is the set of isolated points of~$A$: we will use
the notation $\isol(A) := \reg^0_0(A)$.
%For every $A \subseteq \K^n$, let $\isol(A)$ be the set of isolated points
%of~$A$. Notice that $\isol(A)$ is discrete.
\end{definizione}
% Let $\reg^p(A) \coloneqq \bigcup_{\mu \in \Pi(n,d)} \reg^p_\mu(A)$. 
%be defined as in~\cite[Def.~8.4 and\S3.4]{miller05}.  
As in the case when $\K$ is an expansion of~$\Rbar$,
$\reg^p_\mu(A)$ is definable, open in~$A$, and a $\Cp$-submanifold of $\K^n$
of dimension~$d$.
Hence, if $A$ is $d^p$-regular, then it is an embedded manifold of
dimension~$d$.
For instance, $\reg^n(A)$ is the interior of~$A$ (inside~$\K^n)$.

\begin{remark}\label{rem:regular-local}
Let $\pi \coloneqq \Pi^n_d$, and $A \subset \K^n$ be definable.
$A$ is $\pi^p$-regular iff, for every $y \in \pi(A)$ and $x \in A_y$, 
there exist $U \subseteq \K^d$ open box around $y$ and $W \subseteq \K^{n-d}$
open box around~$x$, such that $A \cap(U \times W) = \Graph(f)$ 
%(the graph of~$f$) 
for some (definable) $\Cp$-map $f: U \to W$.
\end{remark}

\begin{lemma}\label{lem:regular}
Suppose that $\K$ is \iminimal, and let $A \subseteq \K^{m + n}$ be definable,
such that $B := \set{x \in \K^m: \isol(A_x) \neq \emptyset}$ has interior.
Then, for every $p \in \Nat$, $\reg^p_\pi(A) \neq \emptyset$, 
where $\pi := \Pi^{n + m}_m$.
\end{lemma}
\begin{proof}
Fix $p \in \Nat$; let $V \subseteq B$ be a nonempty open box, and
$C := \bigsqcup_{v \in V} \Pa{\mset v \times \isol(A_v)}$.
Notice that $V \subseteq \pi(C)$.
By Definable Choice, there exists a definable function $f : V \to \K^n$ such
that $\pair{x, f(x)} \in C$ for every $x \in V$.
For every $x \in V$, define
\[\begin{aligned}
f^+(x) &:= \min \Pa{f(x) + 1, \inf\set{y \in A_x: y > f(x)}},\\
f^-(x) &:= \max \Pa{f(x) - 1, \sup\set{y \in A_x: y < f(x)}}.
\end{aligned}\]
Notice that $f^- < f < f^+$ on all~$V$.
By \iminimality, after shrinking~$V$, we can assume that $f$, $f^+$ and $f^-$
are $\Cp$ on~$V$.
It is easy to see that $\Graph(f) \subseteq \reg^p_\pi(A)$.
\end{proof}

\begin{lemma}\label{lem:i-regular-2}
Suppose that every 0-dimensional definable subset of $\K$ has an isolated
point, and let $p \in \Nat$. Then:
\begin{enumerate}
\item $\K$ is \iminimal.
\item Let $\Afam$ be a finite collection of definable subsets of~$\K^n$.
Then, there is a $\Pi$-good partition  $\Part$ of~$\K^n$, compatible
with~$\Afam$, such that $P \setminus \reg^0_\mu(P)$ is
nowhere dense in $P$ for every projection $\mu$ and every $P \in \Part$ such
that $P$ is $\mu$-good.
\item $A \setminus \reg^p(A)$ is nowhere dense in~$A$, for every definable
set~$A$. %of dimension~$d$.
\end{enumerate}
\end{lemma}
\begin{proof}
(1): let $X \subseteq \K$ be definable and with empty interior.
Suppose, for a contradiction, that $\cl X$ contains a nonempty open
interval~$I$, and let $Y := X \cap I$.
Notice that $\dim Y = 0$ and $Y$ is dense in $I$, and therefore it has no
isolated points.

(2) and (3) have the same proof as \cite{miller05}*{Prop.~8.4}
(using Thm.~\ref{thm:imin} and Lemmas~\ref{lem:partition}
and~\ref{lem:regular}).
%\footnote{Check: we added the extra condition on~$d$.}
\end{proof}

% We don't know if the following converse of the above Lemma holds (at least,
% for expansions of $\R$).
% \begin{question}
% Assume that $\K$ is \iminimal.
% Is it true that every $0$-dimensional definable subset of $\K$ has an isolated
% point?%
% \footnote{Check if \cite{FKMS} gives a counter-example.}
% \end{question}

\begin{lemma}
The following are equivalent:
\begin{enumerate}
\item $\K$ is d-minimal;
\item for every $\K' \elem \K$, every subset of~$\K$
is the union of a definable open set and finitely many definable discrete sets;
\item for every $m \in \Nat$ and every definable $A \subseteq \K^{m+1}$ there
exists $N \in \Nat$ such that, for all $x \in \K^m$, either $A_x$ has
interior or is a union of $N$ definable discrete sets;
\item for every $m, n \in \Nat$ and definable $A \subseteq \K^{n + m}$ there
exists $N \in \Nat$ such that for every $x \in \K^m$, either $\dim A_x > 0$,
or $A_x$ is a union of $N$ definable discrete sets.
\end{enumerate}
\end{lemma}
\begin{proof}
($1 \Leftrightarrow 2$) follows from Remark~\ref{rem:discrete}.

% ($2 \Rightarrow 1$).
% Let $A \subseteq \K'$ be definable, with $\K' \elem \K$.
% Define $A_0 := \int A$.
% If $A$ is a finite union of discrete sets, then there exists $N \in \Nat$ such
% that $A^({i}) = 0$ for every $i > N$;%
% \footnote{Check!}
%  moreover, each $A^{(i)}$ is definable
% and discrete (for $i \geq 1)$.
% Thus, $A = A_0 \cup A^{(1)}= \cup \dots \cup A^{(n)}$.

$(2 \Leftrightarrow 3)$ is a routine compactness argument.

(3) is the case $n = 1$ of (4).

$(3 \Rightarrow 4)$.
Induction on~$n$.
The case $n = 1$ is the hypothesis.
Let $n > 1$, and assume we have already proved $(3)$ for each $n'< n$.
Let $C := \set{x \in \K^m: \dim(A_x) = 0}$.
Let $x \in C$.
For each $y \in K$, we have $\dim(A_{\pair{x,y}}) = 0$, and therefore $A_{\pair{x,y}}$
is a union of $N$ discrete definable sets (for some $N$ independent from $x$
and~$y$). 
Moreover, for each $x \in C$, the set
$D(x) := \set{y \in \K: A_{\pair{x,y}} \neq  \emptyset}$ has empty interior,
because $\K$ is \iminimal. 
Thus, for every $x \in C$, $D(x)$~is a union of $M$ definable discrete sets
(for some $M$ independent from~$x$).
Hence, for every $x \in C$, $A_x$~is a union of $NM$ discrete sets, which can
be taken definable.
\end{proof}

% If $X$ is compact and $\rkCB(X) < + \infty$ 
% then it must be a successor ordinal:
% denote by $\rkCBs(X)$ the predecessor of $\rkCB(X)$.

% $\K$~has locally o-minimal open core
% iff every definable subset of $\K$ of dimension $0$
% has rank at most 1 \textbf{to be proved}.
Recall the definition of the Cantor-Bendixson rank (Definition~\ref{def:CB}).

\begin{corollary}
If $\K$ d-minimal iff for every $n \in \Nat$ and every definable
$X \subseteq \K^{n+1}$, there exists $N \in \Nat$ such that,
for all $a \in \K^n$, $\rkCB \Pa{X_a \setminus \interior(X_a)} \leq N$.
\end{corollary}

% \begin{definizione}
% Let $\mu \in \Pi(n,d)$ (see Definition~\ref{def:regular}).
% A definable $d$-dimensional $\Cp$-submanifold $M \subseteq \K^n$ is 
% \intro{weakly $\mu$-special} if, for each $y \in \mu(M)$ and $x \in M_y$, 
% there exist $U \subseteq \K^d$ open box
% around $y$ and $W \subseteq \K^{n-d}$ open box around~$x$, such that
% $A \cap(U \times W) = \Graph(f)$ for some (definable) $\Cp$-map $f: U \to W$.
% We say that $M$ is $\mu$-special if the box $U$ in the above definition does
% \emph{not} depend on~$x$ \rom(but only on~$y$\rom).
% $M$~is (weakly) special if it is (weakly) $\mu$-special for some 
% $\mu \in \Pi(n,d)$.
% \end{definizione}

\begin{proposition}\label{prop:d-special}
Assume $\K$ is d-minimal, let $p \in\Nat$, and $\Afam$ be a finite collection
of definable subsets of $\K^n$.
Then, there exists a finite partition of $\K^n$ into %\emph{weakly} special
regular definable $\Cp$-submanifolds compatible with~$\Afam$.
\end{proposition}

\begin{proof}
As in the proof of \cite{miller05}*{Thm.~3.4.1}, we are reduced to show that if
$A \in \Afam$, with $0 < d := \dim A < n$; then $A$ can be partitioned into
regular $\Cp$-manifolds.
Moreover, we can further assume that $A$ is a $\pi$-good, where $\pi :=
\Pi^n_d$, that $A \setminus \reg^p_\pi(A)$ is nowhere dense in~$A$,
and each $A_x$ is discrete, for every $x \in \K^d$.

Let $M := \reg^p_\pi(A)$.
By definition, $M$~is a $\pi$-regular $\Cp$-manifold.
It suffices to prove that $\pi(A \setminus M)$ is nowhere dense (in~$\K^d$)
to conclude the proof (since then $\fdim(A \setminus M) < \fdim(A)$ and we can
proceed by induction on $\fdim(A)$).
Assume, for a contradiction, that $B \subset \pi(A \setminus M)$ is a nonempty
open box, and let $N := \reg^p_\pi(A \setminus M)$.
By shrinking~$B$, we might assume that $\cll(M_x) = (\cll M)_x$ for every
$x \in B$. 
By Lemma~\ref{lem:i-regular-2}, $M$~is dense in~$A$, and therefore, for every
$x \in B$,
$\cll(A_x) = (\cll A)_x = (\cll M)_x = \cll(M_x)$, that is $M_x$ is dense
in~$A_x$.
However, $A_x$~is discrete, and thus $A_x = M_x$.
%
% Let $S \subseteq A$ be the (definable) set of $\pair{x,y} \in A$ such that,
% for some open bounded boxes $U\subset \K^d$ around~$x$ and $W \subset
% \K^{n-d}$ around~$y$, there exists a (definable) function $f: U \to W$ such
% that $A \cap (U \times W) = \Graph(f)$.
%
% As in~\cite[\S8.5]{miller05}, $A \cap \pi^{-1}(U) \subseteq \reg^0_\pi(A)$.
% Let $M := A \cap \pi^{-1}(\good)$.
% As in~\cite[\S8.5]{miller05}, $M$~is a $\pi$-special $C^0$-submanifold.
% It suffices to prove that $S$ is dense in $\mathrm{int}(\pi A)$ to conclude
% the case $p = 0$.
%
% Assume, for a contradiction, that $B \subseteq A \setminus \pi(S)$ is a
% nonempty open box.
\end{proof}

We say that $A$ is $\pi$-special if it is $\pi$-regular and the box $U$ in
Remark~\ref{rem:regular-local} does \emph{not} depend on~$x$ \rom(but only
on~$y$\rom).

\begin{conjecture}
Assume $\K$ is d-minimal, let $p \in\Nat$, and $\Afam$ be a finite collection of definable
subsets of $\K^n$.
Then, there exists a finite partition of $\K^n$ into \emph{special} definable
$\Cp$-submanifolds compatible with~$\Afam$.
\end{conjecture}
The proof in \cite{miller05}*{Theorem~3.4.1} unfortunately has a mistake, and
we were not able to adjust it to prove the above conjecture.

\subsection{Subsequence selection}\label{sec:sequence-selection}
Again, $\K$ is a \dminimal structure.
If $\K$ were \ominimal (or, more generally, locally \ominimal), then Curve
Selection would hold:
%\footnote If $\K$ is locally o-minimal, then, 
given a definable set $A \subseteq \K^n$ and $b \in \fr
A$, there exists a definable continuous function $\gamma: (0,1) \to A$ such
that  $\lim_{t \to 0^+} \gamma(t) = b$.
If $\K$ is not locally o-minimal, then Curve Selection does not hold: for
instance, let $A \subset (0,1)$ be discrete and definable, such that $0
\in \fr A$: then no such curve $\gamma$ exists.
However, in that case we can use definable sequences instead of definable
curves (see Definition~\ref{def:sequence}).

\begin{definizione}
A \pN set is a definable, discrete, closed, and unbounded subset of~$\K_{\geq 1}$.
\end{definizione}
\begin{fact}
A definably complete structure $\F$ has locally \ominimal open core iff a \pN
subset of $\F$ does not exist.
In particular, $\K$ is not locally o-minimal iff a \pN subset of $\K$ exists.
\end{fact}

\begin{proviso}
For the remainder of this subsection we will assume that $\K$ is \dminimal but
not locally o-minimal (and hence a \pN subset of $\K$ exists).
\end{proviso}

\begin{definizione}\label{def:sequence}
A definable sequence is a definable function $f: D \to \K^n$, such that $D$ is
a \pN set.

Given a definable sequence $f: D \to \K^n$, 
a definable subsequence is the restriction of $f$ to $E$, such that $E
\subseteq D$ and $E$ is unbounded (and hence a \pN set).
\end{definizione}

\begin{lemma}\label{lem:sequence-selection}
Let $A \subset \K^n$ be definable, $b \in \fr A$, and $D$ be a
\pN set.
Then, there exists a definable sequence $f: D \to A$, such that
$\lim_{t \to +\infty} f(t) = b$.

Moreover, $f$ can be chosen uniformly: given  $\Pa{A_x: x \in \K^m}$ a
definable family of subsets of $\K^n$, and $b:  \K^m \to \K^n$ a definable
function, such that for every $x \in \K^m$ $b(x) \in \fr(A_x)$, there exists a
definable function
$f: D \times \K^m \to \K^n$, such that %for every
\[
\begin{aligned}
\forall x \in \K^n\ \forall t \in D \ \Pa{f(x,t) \in A_x}\\
\forall x \in \K^n\ \Pa{\lim_{t \to + \infty} f(t,x) = b(t)}.
\end{aligned}
\]
\end{lemma}
\begin{proof}
By Definable Choice, there is a definable function $f: D \to A$ such that, for
every $d \in D$, $\abs{f(d) - b} < 1/d$.
\end{proof}

In particular, from a definable sequence we can  extract a converging (in $\K
\cup \set{\pm \infty}$) definable subsequence.

\begin{definizione}
Given a definable sequence $f: D \to \K^n$ and $b \in \K^n$, we say that $b$
is an accumulation point for $f$ if 
\[
\forall \eps >0\ \forall N \in D\ \exists d \in D\
\Pa{d > N \et \abs{f(d) - b} < \eps}
\]
(or equivalently
\[
\forall N \in D\ \exists d \in D\
\Pa{d > N \et \abs{f(d) - b} < 1/N}).
\]
\end{definizione}

\begin{lemma}\label{lem:subsequence-selection}
%Let $A \subset \K^n$ be definable and $b \in \fr A$.
%Then, there exists a definable sequence $f: D \to A$, such that
%$\lim_{t \to +\infty} f(t) = b$.
Let $f: D \to \K^n$ be a definable sequence and
$b$ be an accumulation point for~$f$.
Then, there exists a definable subsequence $g$ converging to~$b$
(that is, there exists $E \subseteq D$ definable and unbounded, such that
$\lim_{t \to + \infty, t \in E} f(t) = b$).
\end{lemma}
\begin{proof}
Trivial.
\end{proof}

The following proposition is a uniform version of  Lemma~\ref{lem:subsequence-selection}
(uniform both in the sense that the given convergence is uniform in the parameter~$x$, and
that the domain $E$ of the subsequence does not depend on~$x$).
In \S\ref{subsec:Whitney} we will see an application of it. 

\begin{proposition}\label{prop:uniform-subsequence}
Let $D$ be a \pN set.
Let $f: D \times \K^m \to [0,1]^n$ and $g: \K^m \to [0,1]^n$ be definable
functions.
For every $x \in \K^n$ consider the definable sequence $f_x(t) := f(t,x)$.
Assume that, for every $x \in \K^n$, $g(x)$ is an accumulation point
for~$f_x$.
Then, there exists 
$U \subseteq \K^m$ definable, open, and dense, such that $f$ is continuous on
$D \times U$ and for every
$B \subseteq U$ definably connected $g$ is uniformly continuous on $B$ and
there exists $E \subseteq D$ definable and
unbounded, such that
\[
\forall x \in B\ \lim_{t \to +\infty, t \in E} f(t,x) = g(x)
\]
uniformly on~$B$, and therefore
\[
\forall x \in B\ \lim_{\substack{t \to +\infty, t \in E\\ y \to x, y \in B}} f(t,y) = g(x).
\]
\end{proposition}

%We think that 
In Proposition~\ref{prop:uniform-subsequence} we
\emph{cannot} take $E$ to be independent from the choice of~$B$.%
\footnote{For a counter-example, let $\K$ be some \dminimal non \ominimal
expansion of~$\R$, $A$~be a \pN set,
$X$ be a definable $0$-dimensional subset of $\R$ such that $A$ is the set of
accumulation points of~$X$, and
$h: D \to X \setminus A$ be a definable bijection with some other \pN set~$D$.
Define $g: \R \to A$ as $g(x) := \min\set{a \in A: a \geq x}$
%Extend $h$ to $\R$ by setting $h(x) := h(\min\set{d \in D: d \geq x})$.
and
$f: D \times \R \to X \setminus A$ as $f(t, x) := h(t)$.
Let $E$ be any \pN subset of $D$ and $U$ be any definable open dense subset
of~$\R$.
We claim that there exists some $x \in \R$ such that 
$g(x) \neq \lim_{t \to +\infty, t\in E} f(t,x)$.
In fact, \wloG we can assume that $A \subset U$.
If the claim were false, then for every $a \in A$ we would have
\[
a = g(a) = \lim_{t \to +\infty, t \in E} f(t,a) = \lim_{t \to \infty, t \in E} h(t),
\]
which is absurd, since the left side depends on $a$ while the right side is constant.
}

\begin{example}\label{ex:needle}
Let $D := \set{2^n: n \in \N} \subset \R$ and $\K := \tuple{\R, +, \cdot, <,
  D}$.
Then, $\K$ is \dminimal and $D$ is a \pN subset of $\K$.
Define $f: D \times [0,1] \to [0,1]$ by
\[
f(t,x) :=
\begin{cases}
0 & \text{ if } 1/t \leq x \leq 1\\
1 - t x & \text{ if } 0 \leq x \leq 1/t.
\end{cases}
\]
Then, for every $x \in (0,1]$, $\lim_{t \to + \infty} f_x(t) = 0$, but not uniformly
on $(0,1)$.
Let $F := 1/D$ and
$U := (0,1) \setminus F$.
Then, for each interval $B \subset U$,
$\lim_{t \to + \infty} f_x(t) = 0$ uniformly on~$B$.
\end{example}

Putting together Proposition~\ref{prop:uniform-subsequence}
and Lemma~\ref{lem:sequence-selection}, we have the following result.
\begin{corollary}\label{cor:uniform-selection}
Let $D$ be a \pN set.
Let $\Pa{A_x: x \in \K^m}$ be a
definable family of subsets of $[0,1]^n$, and $g:  \K^m \to \K^n$ a definable
function, such that $\forall x \in \K^m$ $\Pa{g(x) \in \fr(A_x)}$.
Then, there exist
\begin{enumerate}[\normalfont(a)]
\item a definable function $f: D \times \K^m \to [0,1]^n$
\item
$U \subseteq \K^m$ definable, open, and dense
\end{enumerate}
such that
\begin{enumerate}
\item $
\forall x \in \K^n\ \forall t \in D \ \Pa{f(x,t) \in A_x}$,
\item $g$ is uniformly continuous on each definably connected subset of~$U$,
\item $f$ is continuous on $D \times U$,
\item  and for every $B \subseteq U$ definably connected
there exists $E \subseteq D$ definable and unbounded such that
\[
\forall x \in B\ \lim_{t \to + \infty, t \in E} f(t, x) = g(x)
\]
uniformly on~$B$, and in particular
\[
\forall x \in B\ \lim_{\substack{t \to +\infty, t \in E\\ y \to x, y \in B}} f(t,y) = g(x).
\]
\end{enumerate}
\end{corollary}

Before proving Proposition~\ref{prop:uniform-subsequence}, we need an easy
lemma: notice that  Lemma~\ref{lem:uniform-continuous} is false if $\K$ is locally o-minimal!
\begin{lemma}\label{lem:uniform-continuous}
Let $g: (0,1)^m \to \K^n$ be definable.
Then, there exists $U \subseteq (0,1)^m$ open and dense such that,
for every $B \subseteq U$ definably connected,
$g \rest B$ is uniformly continuous.
\end{lemma}
\begin{proof}
By Theorem~\ref{thm:imin}\eqref{en:i-f}, there exists $U' \subseteq (0,1)^m$
definable, open, and dense, such that $g \rest{U'}$ is continuous.
%\wloG we can assume that $g$ is continuous.
For uniform continuity, the only problematic points are the ones in 
$C := \fr{U'}$ (where $\fr$ is taken inside $[0,1]^m$).
Let $D$ be a \pN set.
Let $U := \set{x \in U: 1/d(x, C) \notin D}$.
Then, $U$ is open and dense in $(0,1)^m$, and, for every $B \subseteq U$
definably connected, $\fr B \cap C = \emptyset$: therefore, $g \rest B$ is
uniformly continuous.
\end{proof}

\begin{example}
Let $\K$ and $U$ be as in Example~\ref{ex:needle}.
Let $g: (0,1) \to \K$, $g(x) := 1/x$.
Then, $g$ is continuous but not uniformly continuous; however,
it is easy to check that $U$ satisfies the conclusion of Lemma~\ref{lem:uniform-continuous}.
\end{example}

\begin{proof}[Proof of Proposition~\ref{prop:uniform-subsequence}]
% By replacing $\K^n$ by a dense open definable subset $U$, we can assume that
% $g$ is continuous; hence, if $B \subseteq U$ is definably connected, we can
% assume that $g \equiv 0$ on~$B$.%
% \footnote{To be precise, we need that $g$ is uniformly continuous on~$B$.
% If $g$ is not uniformly continuous, we can further shrink $U$ and get that $g$
% is uniformly continuous on each definably connected subset of~$U$.
% For instance, if $U = (0,1)$ and $g$ is continuous on $U$,
% the problematic points (for uniform continuity of $g$) are only $0$
% and~$1$: 
% we can replace $U$ with $U' := U \setminus 1/D \cup 1 - 1/D$: then if $B$ is
% a subinterval of $U'$, then $0$ and $1$ are not boundary points of~$B$, and
% hence $g$ is uniformly continuous on~$B$.}
Let $V \subseteq D \times \K^m$ be open and dense in $D \times \K^m$, such that
$f$ is continuous on~$V$.
By applying Theorem~\ref{thm:imin}\eqref{en:i-B} to the set $V$ (and the fact
that $D$ is a \pN set), there exists $U \subseteq \K^m$ definable, open, and
dense, such that $D \times U \subseteq V$.
%we can assume that $V = D \times B$, and hence that $\tau$ is continuous on all $D
%\times B$.

By Lemma~\ref{lem:uniform-continuous}, after shrinking $U$ if necessary,
%there exists $U \subseteq \K^m$ open, definable, and dense, such that, 
we can assume that,
given $B \subseteq U$ definably connected,
$g \rest B$ is uniformly continuous.
%Hence, \wloG we can assume that $g \equiv 0$ on~$B$.
%By definable choice, there exists a definable function
Define $\tau: D \times U \to D$ as
\[
\tau(t,x):= \min\set{s \in D:  \abs{f(s,x) - g(x)} < 1/t \et s > t}. 
\]
%(where $B$ is some definably connected subset of~$U$), 
Thus,  
%for every $d \in D$ and $x \in U$,
\[
\forall d \in D\ \forall x \in U \ \Pa{\abs{f(\tau(d, x),x) - g(x)}  < 1/d \et \tau(d,x) > d}.
\]
By applying the same trick as at the start of the proof, 
after further shrinking~$U$,
we can further assume that $\tau$ is continuous
on $D \times U$.
%Let $V \subseteq D \times B$ be open and dense in $D \times B$, such that $\tau$ is
%continuous on~$V$.
%By applying Theorem~\ref{thm:imin}\eqref{en:i-B} to the set $V$ (and the fact that $D$ is a \pN set), after shrinking $B$ if necessary,
%we can assume that $V = D \times B$, and hence that $\tau$ is continuous on
%all $D \times B$.
Let $B \subseteq U$ be definably connected.
Since $D$ is discrete, for every $d \in D$, the function $\tau(d, \cdot)$
%does not depend on $x \in B$: 
is constant on~$B$:
we can therefore denote $\lambda(d) := \tau(x,d)$.
Let $E \subseteq D$ be the image of~$\lambda$.
Thus, for every $d \in D$, we have
\begin{equation}
\label{eq:uniform-lambda}
\begin{aligned}
\abs{f(\lambda(d), x) - g(x)} & < 1/d,\\
\lambda(d) &\in E;\\
\lambda(d) &> d. 
\end{aligned}
\end{equation}
Fix $\eps > 0$.
By \eqref{eq:uniform-lambda} we have that there exists $M > 0$ \st
\begin{equation}
\label{eq:uniform-lambda-2}
\forall x \in B\ \forall d \in D\
\Pa{d > M \rightarrow \abs{f(\lambda(d), x) - g(x)} < \eps}.
\end{equation}
Since $\lambda \to + \infty$, the following function $\mu: E \to D$ is
well-defined
\[
\mu(e) := \max\set{d \in D: \lambda(d) = e}.
\]
Notice that $\forall e \in E$ $\Pa{e = \lambda(\mu(e))}$;
therefore, 
\begin{equation}
\label{eq:uniform-mu}
\forall x \in B\ \forall e \in E\ \Pa{\abs{f(e, x) - g(x)} = \abs{f(\lambda(\mu(e)), x) - g(x)} < 1/\mu(e)}.
\end{equation}

\begin{claim}\label{cl:uniform-mu-infty}
$\lim_{e \to + \infty, e \in E} \mu(e) = + \infty$.
\end{claim}
Otherwise, there exists $N > 0$ \st
\begin{equation}
\label{eq:uniform-mu-1}
\forall L >0 \ \exists e \in E\ \Pa{e > L \et \mu(e) \in D_{\geq N}},
\end{equation}
%where $F := \set{d \in D: d < N}$.
Notice that $D_{\geq N}$ is pseudo-finite: hence $F := \lambda(D_{\geq N})$ is
also pseudo-finite, and in particular $F$ is bounded by some $L > 0$.
Take $e \in E$ as in \eqref{eq:uniform-mu}: then, $e = \lambda(\mu(e)) \in F$
and $e > L$, absurd.

By Claim~\ref{cl:uniform-mu-infty}, there exists $N > 0$ such that
$\forall e \in E$ $\Pa{e > N \rightarrow \mu(e) > M}$,
and therefore by \eqref{eq:uniform-lambda-2} we have
\[
\forall x \in B\ \forall e \in E\
\Pa{e > N \rightarrow \abs{f(e, x) - g(x)} < \eps}. \qedhere
\]
%there exists $N >0$ such that
%$\forall d \in D$ $\Pa{\lambda(d) > N \rightarrow d > M}$.
%We have to show:
%\[
%\forall \eps > 0\ \exists N \in E\ \forall x \in B\ \forall s \in E\
%\Pa{s > N \rightarrow \abs{f(s,x) - g(x)} < \eps}}.  
%\]
\end{proof}

\begin{corollary}\label{cor:zero-derivative}
Let $N$ be a \pN set.
Let $f: N \times \K^m \to \K^n$ be a definable function; write $f_t(x) :=
f(t,x)$.
Assume that, for every $x \in \K^m$, $\lim_{t \to + \infty} f_t(x) = 0$.
Then, there exists $U \subseteq \K^m$ open, definable, and dense, such that
$\forall x \in U$ $\lim_{t \to +\infty} D f_t(x) = 0$ (where $D$ denotes the
derivation with respect to the $x$ variables).
\end{corollary}
\begin{proof}
Assume not.
By Proposition~\ref{prop:uniform-subsequence}, after shrinking $N$ and taking
an open subset of $\K^m$ if necessary, there exists $C > 0$ such that
$\forall x \in \K^m$ $\forall t \in N$ large enough
$D f_t(x) > C$.
Fix $x \neq x' \in \K^n$; let $b := \abs{x-x'}$.
Then,
\[
\abs{f_t(x) - f_t(x')} \geq C b
\]
while the left side goes to $0$ as $t \to + \infty$, absurd.
\end{proof}

Notice that Proposition~\ref{prop:uniform-subsequence}   remains true with the weaker hypothesis
that $\K$ is \iminimal but not locally \ominimal (with the same proof), while
Lemma~\ref{lem:sequence-selection} requires also DSF.

\subsection{Pillay rank in d-minimal strustures}
\label{subsec:rank-2}
Remember that $\K$ is a \dminimal structure.
Let $X \subseteq \K^n$ be definable and closed.
In \S \ref{subsec:rank} we defined the Pillay rank of $X$, and computed the
rank of $\K$.
Here we want to give an upper bound on $\rkP(X)$, and in particular prove that
it's not $\infty$ (but it will be, in general, an ordinal).

First, let's consider the locally o-minimal case, when the rank is finite.

\begin{lemma}\label{lem:rkP-lomin}
Let $\K$ be locally o-minimal.
\begin{enumerate}
\item $\rkP(\K^n) = n$.  
\item Let $X \subseteq \K^n$ be definable.  Then, $\rkP(X) = \dim(X)$.
\end{enumerate}
\end{lemma}
\begin{proof}
It's an easy exercise to show that $\rkP(\K^n) \geq n$.
Thus, it suffices to prove 2).
Proceed by induction on $d := \dim(X)$.
If $d = 0$, then $X$ is discrete, and therefore $\rkP(X) = 0$.
For the inductive step, assume $d > 0$.
%First, we do the case when $X$ is an embedded manifold.
Let $Y \subset X$ be definable and nowhere dense.
Since $\K$ is locally o-minimal, $\dim(Y) < d$.
%Since $X$ is a manifold, $\dim(Y) < d$.
Hence, by inductive hypothesis, $\rkP(Y) < d$, and therefore $\rkP(X) \leq d$.
For the opposite inequality, it's easy to find $Y \subset X$ definable such
that $\dim(Y) = d-1$ and $Y$ is a manifold.
Therefore, by inductive hypothesis,
$\rkP(Y) = d-1$.
Since $Y$ is a manifold and $\dim(Y) < \dim(X)$,
$Y$ is nowhere dense in~$X$:
therefore, $\rkP(X) \geq d$.
%
% Now consider the general case.
% Let $Y := \reg^0(X)$.
% Notice that $Y$ is an embedded manifold of dimension~$d$, it is open in~$X$,
% and that, by Lemma~\ref{lem:i-regular-2}, $Y$ is dense in~$X$.
% Therefore, since $\K$ is locally o-minimal, $\dim(X -\setminus Y) < d$%
% \footnote{Give reference.}
% Thus, by the previous case, $\rkP(Y) \leq d$.
%
\end{proof}

Let us consider now the general \dminimal case.
\begin{thm}\label{thm:rkP}
Let $\K$ be \dminimal.
\begin{enumerate}
\item If $\K$ is locally o-minimal, then $\rkP(\K^n) = n$.
\item Assume that is not locally o-minimal, then
$\rkP(\K^n) = \omega^n$.
\item If $Y \subseteq \K^n$ is an embedded manifold of dimension~$d$, then
$\rkP(Y) \leq \omega^d$, with equality holding when $\K$ is not locally o-minimal.
\item If $X \subseteq \K^n$ is definable of dimension~$d$, then there exists a
natural number $m$ such that $\rkP(X) \leq m \omega^d$.
\end{enumerate}
\end{thm}
\begin{proof}
1) is Lemma~\ref{lem:rkP-lomin}.

For 2), like in Proposition~\ref{prop:rkP-unary}, it is easy to see that
$\rkP(\K^n) \geq \omega^n$.
Let $Y \subseteq \K^n$ be definable and nowhere dense.
Thus, by 4), $\rkP(Y) < \omega^{\dim(Y) + 1} \leq \omega^n$, and therefore
$\rkP(\K^n) \leq \omega^n$.

3) and 4) are proved together by induction on~$d$.

If $d = 0$, then $Y$ is discrete, and therefore $\rkP(Y) = 0$.
$X$ is a union of finitely many discrete sets, and therefore $\rkP(X) \leq m$
for some natural number~$m$.

Assume now that $d > 0$ and that
we have proven 3) and 4) for every $d' < d$: we want to prove
them for~$d$.
Let $Z \subseteq Y$ be definable and nowhere dense.
Since $Y$ is a manifold, $\dim(Z) < d$; thus, by inductive hypothesis,
$\rkP(Z) < \omega^d$, and hence $\rkP(Y) \leq \omega^d$.
Thus, it remains to prove 4).

By Proposition~\ref{prop:d-special}, there exists a partition
$X = X_1 \cup \dots \cup X_r$ of $X$ into embedded manifolds $X_i$.
Let $\gamma := \max_{i=1,\dotsc r}\rkP(X_i)$.
By 3), $\gamma \leq \omega^d$.
Notice that each $X_i$ is locally closed in $X$ (since it is locally closed in
$\K^n$).
Thus, by Proposition~\ref{prop:rkP-lclosed}, 
$\rkP(X) \leq r \gamma + (r-1) < (r+1) \omega^d$, and we
are done.
\end{proof}

Let $Y \subset X$ be definable sets with $Y$ nowhere dense in~$X$.
In o-minimal structures we have the inequality $\dim(Y) < \dim (X)$,
which is quite useful in proving theorem by induction on the dimension.
In \dminimal structure such inequality fails.
The usefulness of the Pillay rank is that we have $\rkP(Y) < \rkP(X)$, and
therefore we can proceed by induction on it.
\S\ref{subsec:stratification} will show an application of the above idea.

Moreover, $\rkP$ is a refinement of both the dimension and the
Cantor-Bendixson rank: that is, if $\dim(X) < \dim(Y)$, then $\rkP(X) <
\rkP(Y)$, and similarly for~$\rkCB$.

\subsection{Stratification}
\label{subsec:stratification}
\begin{definizione}
Let $X \subseteq \K^n$ be a definable set.  A \intro{weak stratification} of
$X$ is a finite partition $X= C_1 \sqcup \dots \sqcup C_k$, such that each
$C_i$ is definable, and for each $i$, $\partial C_i$ is a union of some of the
$C_j$ (where $\partial C_i$ is the set of points in the closure of $C_i$
inside $X$ but not in~$C_i$).  
The $C_i$ are the strata of the given stratification.
\end{definizione}

The difference with the usual definition of stratification (for subsets of
$\R^n$) is that in the latter case the strata are required to be connected.
Since we will only consider weak stratifications, we will drop the adjective
weak: a stratification will be a weak stratification.

\begin{definizione}
Let $X \subseteq Y$ be sets, $\Afam$ a family of subsets of~$Y$.
Let $\Bfam$ be partition (or a stratification) of~$Y$; we say that $\Bfam$ is
compatible with $\Afam$ if, for every $A \in \Afam$,
$A \cap X$ is a union of sets in~$\Bfam$.

We also denote $\partial \Afam := \set{\partial A: A \in \Afam}$.
\end{definizione}
Thus, a stratification is a partition $\Bfam$ which is compatible with
$\partial \Bfam$.

\begin{thm}\label{thm:stratification}
Let $\Afam$ be a finite collection of definable subsets of~$\K^n$ and $p \in \Nat$.
Then, there exists a stratification of $\K^n$ into embedded $\Cp$ manifolds
compatible with~$\Afam$.
\end{thm}
\begin{proof}
We prove (by induction) the following statement.
\begin{sentence}[(*)]
Let $X \subseteq \K^n$ be a definable set.
Let $\Afam$ be a finite collection of definable subsets of~$\K^n$ and $p \in
\Nat$.  Then, there is a stratification of $X$ into embedded $\Cp$
manifolds compatible with~$\Afam$.
\end{sentence}
By Proposition~\ref{prop:d-special}, there exists a finite partition of $\K^n$
in $\Cp$ manifolds compatible with~$\Afam \cup \mset X$.
\Wlog, we can assume that $\Afam$ itself is such partition: we want to find a
stratification refining~$\Afam$.

Let $\Afam' := \set{A \in \Afam: A \subseteq X}$.
We want to find a stratification of $X$ refining $\Afam'$.
Notice that if $A$ were open in $X$ for every $A \in \Afam'$, then $\Afam'$
would already be a stratification of~$X$.

% We proceed by induction first on $d := \dim(X)$.
% If $d = n$, let $\Bfam := \set{A \in \Afam': \dim(A) = n}$, and
% $Y := \bigcup \Bfam$.
% Notice that $Y$ is open in $X$ and $\Bfam$ is a stratification of $Y$ compatible
% with~$\Afam$ (because each $A \in \Bfam$ is open in $\K^n$ and a fortiori in $Y$).
% Moreover, $\dim(X \setminus Y)  < n$.
% Thus, by inductive hypothesis, there exists a stratification
% $\Bfam'$ of $Y$ which is compatible with $\Afam$ and with 
% $\set{\partial B: B \in \Bfam}$.
% Then, $\Cfam := \Bfam \cup \Bfam'$ is a stratification of $X$ compatible
% with~$\Afam$ (since for every $B \in \Bfam$ and $B' \in \Bfam'$, 
% we have $B$ disjoint for  $\cl{\Bfam'}$, because $B$ is open).

% Thus, we can assume $d < n$.
% If $d = 0$, then $m := \rkCB(X)$ is finite: we proceed by induction on~$m$.
% If $m = 1$, then $X$ is discrete, and $\Afam'  := \set{A \in \Afam: A
%   \subseteq X}$ is  a stratification of $X$ compatible with $\Afam$.
% If $m > 1$, let $Y$ be the set of isolated point of $X$ and $Z := X \setminus
% Y$.
% Again, $\Afam \rest Y := \set{A \cap Y: A \in \Afam}$ is a stratification of $Y$ compatible with $\Afam$.
% Moreover, $\rkCB(Z) < m$, and therefore, by inductive hypothesis,
% there exists a stratification $\Bfam$ of $Z$ compatible with $\Afam$ and with
% $\set{\partial(A \cap Y): A \in \Afam}$.
% Then, $\Afam\rest Y \cup \Bfam$ is a stratification of $X$ compatible with
% $\Afam$ (as before, because  $A \cap Y$ is open in $X$ for every $A \in \Afam$).

% Thus, we can assume $0 < d < n$.
We proceed by induction on $\alpha := \rkP(X)$ (by Theorem~\ref{thm:rkP},
$\alpha$ is an ordinal number).
Let $d := \dim(X)$; for each $i = 0, \dotsc, d$, let $Y_i := \reg^p_d(X)$ and
$Y := \bigcup_{i=0}^d Y_i = \reg^p(X)$.
%Let $Y := \reg^p(X)$.
Notice that each $Y_i$ is a $\Cp$-manifold of dimension $i$ and is open in~$X$.
For $i = 0 \dotsc d$, let $\Bfam_i := \set{A \cap Y_i: A \in \Afam'  \et
  \dim(Y_i \cap A) = i}$ and $Y_i' := \bigcup \Bfam_i$.
Notice that, for every $B \in \Bfam_i$, $B$~is an open subset both of the
corresponding~$A \in \Afam$ (because $Y_i$ is open in $X$ and $A$ is a subset of~$X$)
and of~$Y_i$ (because $Y_i$ is a manifold and $B$ is a submanifold of the same
dimension as~$Y_i$) and hence of $X$ (because $Y_i$ is open in~$X$).
Let $\Bfam := \bigcup_{i=0}^d \Bfam_i$ and $Y' := \bigcup_{i=0}^d Y'_i =
\bigcup \Bfam$.
Thus, $\Bfam$ is a stratification of $Y'$ compatible with $\Afam$.
Moreover, for each $i = 0, \dotsc, d$,
$\dim(Y_i \setminus Y_i') < i$,  thus, since $Y_i$ is a manifold,
$Y_i'$~is dense in $Y_i$: therefore, $Y'$ is dense in~$Y$.
By Lemma~\ref{lem:i-regular-2}, $Y$~is dense in~$X$,
 and therefore $Y'$ is dense in~$X$.
Let $Z := X \setminus Y'$.
Since $Z$ is nowhere dense in~$X$, we have $\rkP(Z) < \alpha$.
Thus, by inductive hypothesis, there exists $\Cfam$ stratification of $Z$
compatible with $\Afam \cup \partial \Bfam$.
Then, $\Bfam \cup \Cfam$ is a stratification of $X$ compatible with~$\Afam$
(because, for every $B \in \Bfam$ and $C \in \Cfam$, $B$ and $C$ are disjoint
and $B$ is open in~$X$, and hence $B$ and $\partial C$ are also disjoint).
\end{proof}

Among all possible stratifications of $X$ there is a canonical one.
% a coarsest one.
\begin{definizione}\label{def:canonical-stratification}
Let $X \subseteq \K^n$ be a definable set, $p \in \N$, and $\Afam$ be a finite
collection of definable subsets of~$\K^n$.  We now define $\Sp(X,\Afam)$ the
canonical $\Cp$ stratification of $X$ compatible with $\Afam$ by induction on
$\rkP(X)$.  By replacing $\Afam$ with the atoms of the Boolean algebra
generated by it, \wloG we can assume that $\Afam$ is a partition of~$\K^n$.
Let $d := \dim(X)$; define
\[
\begin{aligned}
Z &:= \reg^p_d(X);\\
\Bfam &:= \Bfam^p(Z, A) := \set{\reg^p_d(A \cap Z): A \in
  \Afam \et \dim(A \cap Z) = d },\\
W &:= \bigcup \Bfam,\\
Y & := X \setminus W,\\
\Afam' &:= \Afam \cup \Bfam.
\end{aligned}
\]  
%Let $Y := X \setminus \bigcup \Bfam$ and $\Afam' := \Afam \cup \partial \Bfam$.
%Notice that $Y$ is nowhere dense in~$X$: 
%Notice that 
Notice that $W$ is a $\Cp$ manifold, it is open in~$X$, 
and that $\Bfam$ is a partition of~$W$ by open sets,
and therefore $\Bfam$ is a stratification of~$W$ compatible with~$\Afam$.
By Lemma~\ref{lem:rkP-regular}, $\rkP(Y) < \rkP(X)$: 
therefore, by induction on $\rkP$, we can assume that we have
already defined $\Sp(Y, \Afam')$, the canonical $\Cp$ stratification of~$Y$
compatible with~$\Afam'$.  
Then, $\Sp(X, \Afam) := \Bfam \cup \Sp(Y, \Afam')$
is the canonical $\Cp$ stratification of $X$ compatible with~$\Afam$ (notice
that $\Sp(X, \Afam)$ is indeed a stratification of $X$ compatible
with~$\Afam$).

If we don't specify the family~$\Afam$, the canonical $\Cp$ stratification of
$X$ is $\Sp(X) := \Sp(X, \emptyset)$.
\end{definizione}

\begin{lemma}\label{lem:rkP-regular}
In the setting of Definition~\ref{def:canonical-stratification}, we have
$\rkP(Y) < \rkP(X)$.
\end{lemma}
\begin{proof}
Notice that it's not true in general that $Y$ is nowhere dense in~$X$: thus,
we do need to give an argument.

Assume, for a contradiction, that $\rkP(Y) = \rkP(X) =:\alpha$.

Let $T := \cl Z$ (where the closure is taken inside~$X$).
\begin{claim}\label{cl:dim-T}
$\dim(X \setminus T) < d$.
\end{claim}
In fact, by Lemma~\ref{lem:i-regular-2}, we have
$X \setminus T \subseteq \bigcup_{i=0}^{d-1} \cl{\reg_i^p(X)}$.
Since, for every $i = 0, \dotsc, i-1$,
 we have that either $\reg_i^p(X)$ is empty, or 
$\dim(\cl{\reg_i^p(X)}) = \dim(\reg_i^p(X)) = i < d$, the claim follows.

Let $Y_1 := Y \cap T$ and $Y_2 := \cll_Y(Y \setminus T)$.
\begin{claim}\label{cl:rkP-Y}
$\rkP(Y_1) = \alpha$.
\end{claim}
In fact, we have that $Y$ is the union of it two closed (in~$Y$) subsets $Y_1$
and~$Y_2$, and therefore, 
by Fact~\ref{fact:rkP-closed}, either $\rkP(Y_1) = \alpha$, or $\rkP(Y_2) =
\alpha$.
However, $\dim(Y_2) \leq \dim(X \setminus T) < d$, and therefore, by
Theorem~\ref{thm:rkP}, $\rkP(Y_2) < \omega^d$; on the other hand $\dim(X) =
d$, and therefore $\alpha \geq \omega^d$.

% \begin{claim}\label{cl:rkP-Z}
% $\rkP(T) = \alpha$.
% \end{claim}
% By  Claim~\ref{cl:dim-T}, $\dim(X \setminus T) < d$; we conclude by reasoning
% as in Claim~\ref{cl:rkP-Y}.

Thus, it suffices to prove the following claim to get a contradiction.
\begin{claim}\label{cl:Y1-nd-T}
$Y_1$ is nowhere dense in~$T$.
\end{claim}

We need some further results before proving Claim~\ref{cl:Y1-nd-T}.

\begin{claim}
$W$ is dense in~$T$.
\end{claim}
Since $Z$ is dense in $T$, it suffices to show that $W$ is dense in~$Z$.
Assume, for a contradiction, that $U$ is nonempty definable subset of 
$Z \setminus W$ which is open in~$Z$.
Since $W$ is a manifold of dimension~$d$, we have that $\dim(U) = d$.
Thus, there exists $A \in \Afam$ such that $\dim(A \cap U) = d$.
By Lemma~\ref{lem:i-regular-2}, $\reg^p_d(A \cap U)$ is nonempty
(or, equivalently, $\interior_U(A \cap U)$ is nonempty).
However, $\reg^p_d(A \cap U) \subseteq \reg^p_d(A \cap Z)$ 
(again, because $Z$ is a manifold of dimension~$d$, and hence
$\reg^p_d(A \cap Z) = \interior_{Z}(A \cap Z)$, and similarly for~$U$).
Thus, $\reg^p_d(A \cap Z) \in \Bfam$ and therefore
$\reg^p_d(A \cap Z) \subseteq W$; hence, $\reg^p_d(A \cap U)$ is a nonempty
subset of $U \cap W$, absurd.

Thus, since $W$ is also open in~$T$, Claim~\ref{cl:Y1-nd-T} follows from the
fact that $Y \cap W = \emptyset$.
%it suffices to prove the following.
%\begin{claim}
%$Y \cap W$ is nowhere dense in~$W$.
%\end{claim}
\end{proof}

\subsection{Verdier and Whitney stratifications}\label{subsec:Whitney}

Whitney stratifications are particularly important in semialgebraic and
subanalytic geometry.

\cite{Loi} proved that in o-minimal structures every definable set admits a
Whitney stratification (he stated his result for expansions of~$\R$, but the
same proof works in general o-minimal structures).
To do it, he introduced Verdier stratifications
(see Def.~\ref{def:Verdier}), 
and proved that, for an o-minimal structure,
\begin{enumerate}
\item Every Verdier stratification (with finitely many definable sets) is a Whitney
stratification
\item Every definable set admits a Verdier stratification 
(into finitely many definable sets).
\end{enumerate}

Unfortunately, Loi's result does not extends to \dminimal structures
(see Example~\ref{ex:no-Verdier}).

% \begin{conjecture}
% There exists a Whitney stratification of~$X$.
% \end{conjecture}

% We try to prove the above conjecture; we follow the ideas in \cite{Loi}, where
% the conjecture is proven for the case when $\K$ is o-minimal.

However, some partial result is still true.

\begin{definizione}\label{def:pair-property}
Let $P(X,Y;y)$ be a property, where $X$ and $Y$ are definable embedded
submanifolds of $\K^n$, $Y
\subseteq \partial X$, and $y \in Y$.

We say that $P$ is \intro{local} if, for every $y \in Y$ and
every definable open neighbourhood $U$ of~$y$, 
$P(X,Y;y)$ holds iff $P(X \cap U, Y \cap U; y)$ holds.

We denote $P(X,Y) := \set{y \in Y: P(X,Y;y) \text{ holds}}$.
\end{definizione}

\begin{definizione}
In the setting of Definition~\ref{def:pair-property}, we say that: 
\begin{enumerate}
\item Whitney Property (a) holds if\\
$a(X,Y;y) :=$
for every $\eps > 0$ there exists a neighbourhood $U$ of $y$ in $\K^n$
such that
\[
\forall x \in X \cap U\ \forall z \in Y \cap U\
\delta(T_z Y, T_x X) \leq \eps,
\]
where
\[
\delta(A,B) := \sup_{a \in A, \norm a = 1} d(a,B)
\]
is the distance between the vector subspaces $A, B \subseteq \K^n$.
\item Whitney Property (b) holds if\\
% \footnote{Check that it is indeed Whitney (b).}  the following happens.  Let
% $k:= \dim X$, $T: X \to \Gras_k(\K^n)$ be the map $x \mapsto T_x X$, and
% $\ell: X \times Y \to \Gras_1(\K^n)$ be the map $\pair{x,z} \mapsto \K \cdot
% (x-z)$.
%
$b(X,Y;y) :=$ for every $\eps >0$ there exists a neighbourhood $U$ of $y$ in
$\K^n$ such that
%\exists U \text{ neighbourhood of } y\ 
\[
\forall x
\in X \cap U\ \forall z \in Y \cap U\ \delta(\K\cdot (x - z), T_x X) <
  \eps.
\]
\end{enumerate}
A $\Cone$ stratification $\Sfam$ is an a-stratification if for every $X,
Y \in \Sfam$ such that $Y \subseteq \partial X$, we have $a(X,y) = Y$.
If moreover $b(X,y) = Y$, then $\Sfam$ is a Whitney stratification.
\end{definizione}

\begin{remark}
If $X$ and $Y$ are $\Ctwo$ manifolds and $\K$ is an expansion of the real
field, then the definition of the Whitney properties (a) and (b) given above
coincides with the usual Whitney conditions (defined in \cite{Whitney}).
\end{remark}

\begin{definizione}[\cite{Loi}*{p.~348}]\label{def:Verdier}
In the setting of Definition~\ref{def:pair-property}, we say that $w(X,Y;y)$ holds
if there exists a constant $C > 0$ and a neighbourhood $U$ of $y$ in $\K^n$
such that
\[
\forall x \in X \cap U\quad \forall z \in Y \cap U\quad
\delta(T_z Y, T_x X) \leq C \norm{z - x} 
\]
A $\Cone$ stratification $\Sfam$ is a Verdier stratification if for every $X,
Y \in \Sfam$ such that $Y \subseteq \partial X$, we have $w(X,y) = Y$.
\end{definizione}

% \begin{definizione}
% In the setting of Definition~\ref{def:pair-property}, we say that $T(X,Y;y)$ holds
% if for every $\eps > 0$ there exists a neighbourhood $U$ of $y$ in $\K^n$
% such that
% \[
% \forall x \in X \cap U\ \forall z \in Y \cap U\
% \delta(T_z Y, T_x X) \leq \eps.
% \]
% A $\Cone$ stratification $\Sfam$ is an  a-stratification if for every $X,
% Y \in \Sfam$ such that $Y \subseteq \partial X$, we have $T(X,y) = Y$.
% \end{definizione}

Notice that $w(X,Y;y)$ implies $a(X,Y;y)$: thus, by Loi's Theorem,
every  \ominimal
structure admits an a-stratification.
We will show that every \dminimal structure admits an a-stratification.

We need first a general result about stratifications (\cf \cite[Prop.~2]{LSW}).

\begin{remark}
Properties (a), (b), and (w) are local properties.
\end{remark}

\begin{proposition}\label{prop:local-stratification}
Let $P(X,Y;y)$ be a local property.
Suppose that for every pair $\pair{X,Y}$ of definable $\Cp$ submanifolds of
$\K^n$ with $Y \subseteq \partial X$ and $Y \neq \emptyset$, we have
\begin{sentence}[(*)]
the set $P(X,Y)$ is definable and nonempty.
\end{sentence}
Then, with the same assumptions on $X$ and $Y$, we have
\begin{sentence}[(**)]
the set $P(X,Y)$ is definable and dense in~$Y$.
\end{sentence}
Moreover, let $F \subseteq \K^n$ be definable and
$\Afam$ be a family definable subsets of~$\K^n$.
Then, there exists a $\Cp$ stratification $\Dfam$ of $F$ compatible with
$\Afam$, and such that 
\begin{sentence}[(***)]
$P(X,Y) = Y$ for every $X,Y \in \Dfam$ with $Y \subseteq \partial X$.
\end{sentence}
\end{proposition}
\begin{proof}
We first show that (*) implies (**).
If not, then there exists $X$ and $Y$ as in the assumption such that $P(X,Y)$
is not dense in $Y$: therefore, 
% since $Y$ is a manifold and $\K$ is \dminimal,
there exists $Y' \subseteq Y$ definable, open, and nonempty, such that $P(X,Y)
\cap Y' = \emptyset$.
Let $U \subseteq \K^n$ be definable and open, such that $Y' = Y \cap U$, and
let $X' := X \cap U$.
Then, since $P$ is local, we have that $P(X',Y') = \emptyset$, contradicting (*). 

We now show that (**) implies (***).
As usual, we can assume that $\Afam$ is a partition of~$\K^n$ into $\Cp$ manifolds.
We proceed by induction on $\rkP(F)$.
Notice that, since $Y$ is a manifold, (**) is equivalent to
saying that  $Y \setminus P(X,Y)$ is definable and nowhere dense in~$Y$, or
that $\dim(Y \setminus P(X,Y)) < \dim(Y)$.

For the induction to work, we strengthen the conditions we require for the
stratification~$\Dfam$.
\begin{enumerate}[(i)]
\item $\Dfam$ is a $\Cp$ stratification of $F$ compatible with $\Afam \cup \partial
\Afam$;
\item For every $X \in \Afam\cup \Dfam$ and $Y \in \Dfam$, if $Y
\subseteq \partial X$ (equivalently, $Y$ meets $\partial X$), then
$P(X,Y) = Y$.
\end{enumerate}
So, assume we have already proven that, for every $F'$ with $\rkP(F') <
\rkP(F)$ and every finite family of definable sets~$\Afam'$ there exists a
$\Dfam'$ satisfying (i) and (ii) for $F'$ and $\Afam'$; 
we want to find $\Dfam$ for $F$ and~$\Afam$.

%As usual, we can assume that $\Afam$ is a partition of~$\K^n$.
Let $\Afam' := \Afam \cup \partial \Afam$.
We build $\Dfam$ inductively.
Let $\Bfam := \Bfam^p(F, \Afam')$ (as in
Definition~\ref{def:canonical-stratification}).
$\Bfam$ is not yet a good enough starting point, since given $Y \in \Bfam$ and $X \in
\Afam$ with $Y \subseteq \partial X$, it may happen $P(X,Y) \neq Y$.
Thus, given $B \in \Bfam$, let 
\[
e(B) := \bigcap \set{P(A,B): A \in \Afam \et B \subseteq \partial A};
\]
by (**), $e(B)$ is dense in~$B$; since $P$ is a local condition, 
$e(B)$ is also open in~$B$.
%Let
% \[
% \begin{aligned}
% \Efam_B &:= \Bfam^p(e(B), \Afam')\\
% E_B &:= \bigcup \Efam_B.
% \end{aligned}
% \]
% Notice that, for every $A \in \Afam$ and $E \in \Efam_B$, if $E
% \subseteq \partial A$, then $P(A,E) = E$.
% Let 
% \[
% \begin{aligned}
% \Efam &:= \bigcup_{B \in \Bfam} \Efam_B,\\ 
% F &:= \bigcup \Efam,\\
% D' &:= D \setminus F.
% \end{aligned}
% \]
% Notice that $\Efam$ is a stratification of $F$ satisfying (i) and (ii),
% and that $\rkP(D') < \rkP(D)$.%
Let
\[\begin{aligned}
\Efam &:= \set{e(B): b \in \Bfam}\\
F &:= \bigcup \Efam\\
D' &:= D \setminus F.
\end{aligned}\] 
Notice that $\Efam$ is a stratification of $F$ satisfying (i) and (ii).

\begin{claim}
$\rkP(D') < \rkP(D)$.
\end{claim}
%\footnote{Give proof.}
Let $W := \bigcup \Bfam$ and $T := \cl W$.
Reasoning as in the proof of Lemma~\ref{lem:rkP-regular}, it suffices to prove
that $\rkP(D' \cap W) < \rkP(D)$.
But $D' \cap W$ is nowhere dense in~$W$, and the claim follows.

Thus, by inductive hypothesis, we can find a stratification $\Dfam'$
satisfying (i) and (ii) for $D'$ and the family 
$\Afam \cup \Efam$.
Finally, $\Efam \cup \Dfam'$ satisfies (i) and (ii). 
\end{proof}

In general, not every definable set admits a Verdier or a Whitney stratification. 
\begin{examples}\label{ex:no-Verdier}
Let $\K$ be a \dminimal non \ominimal expansion of $\tuple{\R, < +, \cdot,
  \exp}$ and let $D \subset \R$ be a \pN set.
\begin{enumerate}
\item
Consider the following subsets of $\R^2$:
\[
\begin{aligned}
X &:= (0,1) \times \mset 0\\
Y &:= \set{\pair{x, y}: 0 < x < 1 \et \exists t \in D\ \Pa{y = \exp(-t x)}}.
\end{aligned}
\]
Notice that $Y$ is an open subset of $\partial X$ and
$w(X,Y)$ is empty. 
Thus, no Verdier stratification of $\R^2$ compatible with $X$ exists.
\item
Let $X := \set{\R \times \mset{1/d}: d \in D}$ and $Y := \R \times \mset 0$.
Then, $Y = \partial X$ and $b(X,Y) = \emptyset$.
Thus, there is no Whitney stratification of $\R^2$ compatible with~$X$.
Moreover, $w(X,Y) = Y$, and therefore (for a \dminimal structure)
a Verdier stratification is not necessarily a Whitney stratification.
It's an easy exercise to modify the example and also make $X$ connected.
\end{enumerate}
\end{examples}

We can generalize Loi's theorem to locally o-minimal structures.

\begin{lemma}[Wing Lemma]\label{lem:wing}
Assume that $\K$ is locally \ominimal.
Let $V \subseteq \K^k$ be a nonempty open definable set, and
$S \subseteq \K^k \times \K^\ell$ be a definable set.
Suppose that $V \subseteq \partial S$.
Then, there exists a nonempty open subset $U $ of~$U$,
$\alpha_0 > 0$, and a definable map
$\bar \rho: U \times (0, \alpha_0) \to S$, of class $\Cp$,
such that $\bar \rho(y,t) = \pair{y, \bar \rho(y)}$ and 
$\norm{\rho(t)} = t$, for all $y \in U$, $t \in (0, \alpha_0)$.
\end{lemma}
\begin{proof}
Given $x \in \K^k$, denote $S_{[x]} := S \cap (\mset x \times \K^\ell)$.
By Theorem~\ref{thm:imin}\enumref{en:i-fr}, after shrinking~$V$,
\wloG we can assume that, for ever $x \in V$, $x \in \partial{S_{[x]}}$.
By DSF, there exists a definable partial function $\rho: V \times (0, +
\infty) \to \K^\ell$ such that $\rho(x,t) \in V$ and $\norm{\rho(x,t)} = t$
(if there exists an element $y \in S_x$ such that $\norm y = t$, and is
undefined otherwise).
Let $D$ be the domain of $\rho$: notice that, for every $x \in V$,
$0 \in \partial(D_x)$.
Let $\alpha(x) := \sup\set{t >0: (0,t) \subseteq D_{x}}$.
Since $\K$ is locally \ominimal,  $\alpha(x) > 0$ for every $x \in V$.
By Theorem~\ref{thm:imin}\enumref{en:i-f}, after shrinking~$V$,
\wloG we can assume that $\alpha$ is
continuous on $V$, and hence, after shrinking $V$ again,
we can assume that there exists $\alpha_1 > 0$ such that
$\alpha(x) > \alpha_1$ for every $x \in V$.

It remains to shrink $V$ and $\alpha_1$ in order to make $\rho$ a $\Cp$ function.

Be Lemma~\ref{lem:imin-C1}, there exists $W \subseteq V \times (0, \alpha_1)$
definable, open, and dense, such that $\rho$ is $\Cp$ on~$W$.
By Theorem~\ref{thm:imin}\enumref{en:i-fr}, \wloG we can assume that $x
\in \partial(W_{[x]})$, for every $x \in V$.
Repeating the reasoning as above for~$W$, we find $\alpha_0 > 0$ and $U
\subseteq V$ open and definable,
such that $U \times (0, \alpha_0) \subseteq W$.
\end{proof}

\begin{lemma}
Assume that $\K$ is locally \ominimal.
Let $X$ be a definable subset of $\K^n$ and $1 \leq p \in \N$
Then, there exists a definable $\Cp$ Verdier stratification of $\K^n$ compatible with~$X$.
\end{lemma}
\begin{proof}
The proof in \cite{Loi} can be generalized to this situation.
By Proposition~\ref{prop:local-stratification}, it suffices to prove that,
if $X$, and $Y$ are definable $\Cone$ submanifolds of
$\K^n$ with $Y \subseteq \partial X$ and $Y \neq \emptyset$, we have that
the set $w(X,Y)$ is definable and nonempty.
The fact that $w(X,Y)$ is definable is clear, and thus it remains to show that
it is nonempty.
The remainder of the proof is as in \cite{Loi}, using Lemma~\ref{lem:wing}.
\end{proof}

\begin{lemma}
Assume that $\K$ is locally \ominimal.
Let $p \geq 2$ and $X, Y \subseteq \K^n$ be definable $\Cp$ manifolds such
that
$Y \subseteq \partial X$.
Let $y \in Y$ such that $w(X, Y; y)$ holds.
Then, $b(X,Y; y)$ also holds.
\end{lemma}
\begin{proof}
% \begin{claim}
% Let $f: (0, \alpha) \to \K$ be definable, with $f(t) \neq 0$ for every $t> 0$
% and $\lim_{t \to 0^+} f(t) = 0$.
% Then, $\lim_{t \to 0^+} \frac{f(t)}{f(t')} = 0$.
% \end{claim}
% (CHECK)

The   proof in \cite{Loi}*{Prop.~1.10} transfers to this situation
(the main ingredients in the proof is Curve Selection, which holds for locally
\ominimal structures).%
\footnote{There is a misprint in the proof of \cite{Loi}*{Prop.~1.10}, where
  it should be $\phi(t) = (b(t), a(t))$ instead of $\phi(t) = (a(t), b(t))$.}
% (CHECK!)
\end{proof}

\begin{corollary}\label{cor:lomin-Whitney}
Assume that $\K$ is locally \ominimal and $p \geq 2$.
Then, for every definable set $X \subseteq \K^n$ there exists
 exists a definable $\Cp$ Whitney stratification of $\K^n$ compatible with~$X$.
\end{corollary}

For \dminimal structures, we can prove that definable sets have an 
a-stratification.

\begin{proposition}\label{prop:Whitney}
(Let $\K$ be \dminimal.)
Let $X$ be a definable subset of $\K^n$ and $1 \leq p \in \N$.
Then, there exists a definable $\Cp$ $a$-stratification of $\K^n$ compatible with~$X$.
\end{proposition}
\begin{proof}
By Proposition~\ref{prop:local-stratification}, it suffices to prove that,
if $X$, and $Y$ are definable $\Ctwo$ submanifolds of
$\K^n$ with $Y \subseteq \partial X$ and $Y \neq \emptyset$, we have that
the set $a(X,Y)$ is definable and nonempty.
The fact that $a(X,Y)$ is definable is clear, and thus it remains to show that
it is nonempty.

By Corollary~\ref{cor:lomin-Whitney}, we can assume that $\K$ is not locally
\ominimal, and therefore the content of \S\ref{sec:sequence-selection} applies.
Assume, for a contradiction, that $a(X,Y)$ is empty.
Since $a$ is a local property and is invariant under definable
$\Ctwo$ diffeomorphisms, \wloG we can assume that $Y$ is an open subset of 
$\K^k \subset \K^k \times \K^\ell$, where $\ell := n-k$.
In this case, $T_y Y = \K^k$ for every $y \in Y$.

Let $N$ be a \pN set.
By Corollary~\ref{cor:uniform-selection}, 
there exists a definable $\Ctwo$ function
$f: Y \times N \to \K^\ell$ such that
\[
\forall y \in Y\ \forall t \in N\ \Pa{f(y,t) \in X_y \et \abs{f(y,t)} < 1/t}.
\]
%and we can also assume that $X = \Graph(f)$.

Write $f_t(y) := f(y,t)$.
By Corollary~\ref{cor:zero-derivative}, we have that, after shrinking~$Y$,
$\lim_{t \to + \infty} D f_t(y) = 0$, uniformly in~$y$.
But since $T_{\pair{y, f(y,t)}}X \supseteq \Graph (D f_t(y))$, we have that
$a(X,Y) = Y$, absurd.
\end{proof}

We give now an example of a submanifold of $\R^3$ (not definable in any \dminimal
expansion of $\R$) which does not admit an a-stratification.
\begin{example}
Let
\[
\begin{aligned}
X &:= \set{\pair{x, r \sin(x/r), r \cos(x/r)}: 0 < x < 1 \in \R, r \in 2^{-\N}}\\
Y &:= (0,1) \times \mset 0 \times \mset 0.
\end{aligned}
\]
Notice that  $Y \subset \partial X$ and $a(X, Y) = \emptyset$.
Then, there is no a-stratification of $\R^3$ 
(into finitely many $\Cone$ manifolds)
compatible with~$X$.
\end{example}

It remains open the quest of finding a property $b'(X,Y;y)$ such
that:
\begin{enumerate}
\item when $\K$ is o-minimal, $b'(X,Y;y)$ coincides with $b(X,Y;y)$ for
sufficiently smooth definably connected manifolds 
\item  when $\K$ is \dminimal, every definable set
admits a stratification satisfying $a$ and $b'$.
\end{enumerate}
%(And similarly we would like to have an analogue version of property~$w$).

% \begin{proposition}
% Let $X$ be a definable subset of $\K^n$ and $1 \leq p \in \N$
% Then, there exists $\Sfam$ a definable $\Cp$ stratification of $\K^n$
% compatible with~$X$, such that, for every $X, Y \in \Sfam$, if $Y
% \subseteq \partial X$, then $b(X,Y) = Y$.
% \end{proposition}
% \begin{proof}
% As usual, it suffices to show that if $Y$ is an open definable subset of $\K^k$ and
% $X \subseteq \K^k \times \K^\ell$ is a definable $\Ctwo$ manifold,
% with $\emptyset \neq Y \subseteq \partial X$, then $b(X,Y) \neq \emptyset$.
% Assume not.
% The case when $\K$ is locally \ominimal is left to the reader (CHECK); thus,
% we can assume that $N$ is a \pN set.
% \begin{equation}\label{eq:a}
% \forall y \in Y\
% \exists d \in N\ 
% \forall e \in N\
% \exists x \in X \in B(y, 1/e)\ \exists z \in Y \cap B(y, 1/e)\
% \Pa{\delta(\K\cdot(x-z), T_x X) > 1/d}.
% \end{equation}
% By definable choice, we can find definable functions $d: Y \to N$, 
% $x: Y \times N  \to X$ and
% $z: Y \times N \to Y$ satisfying \eqref{eq:a}.
% As usual, after shrinking $Y$, we can assume that $d$ is constant.
% TOBEDONE
% \end{proof}

% -------------------------------------------------------------

\section{Cauchy completion} \label{sec:cauchy}

\begin{definizione}
If $\Kb$ is an ordered field, we denote by $\KbC$ the \intro{Cauchy completion}
of~$\Kb$, that is, the maximal linearly ordered group such that $\Kb$ is dense
in~$\KbC$; notice that $\KbC$ is, in a canonical way, a real closed
field (see~\cite{scott}).
If $\K$ is an expansion of an ordered field~$\Kb$, 
we define $\KC \coloneqq \KbC$ \rom(later we will show how to extend the
structure to~$\KC$: see Proposition~\ref{prop:dmin-Cauchy}\rom).

A \intro{cut} $\Lambda := (\Lambda^L, \Lambda^R)$ of $\K$ is a partition of
$\K$ into  two disjoint subsets $\Lambda^L$ and~$\Lambda^R$, 
such that $y < x$ for
every $y \in \Lambda^L$ and $x \in \Lambda^R$.
We use the notations $y < \Lambda$ if $y \in \Lambda^L$, 
and $x > \Lambda$ if $x \in \Lambda^R$.
A cut $(\Lambda^L, \Lambda^R)$ is a \intro{gap} if $\Lambda^L$ 
is nonempty and has no maximum,
and $\Lambda^R$ is nonempty and has no minimum.
A cut $\Lambda$ is \intro{regular} if $\hat \Lambda = 0$, where
$\hat\Lambda := \inf \set {x - y: x, y \in \K \ \&\ y < \Lambda < x }$.  
\end{definizione}

Notice that the Cauchy completion of $\K$ is the disjoint union of $\K$ and
the set of regular gaps of~$\K$, with suitably defined order $<$ and operations $+$ and~$\cdot$.

\begin{remark}
Let $\kappa$ be the cofinality of $\K$.
The Cauchy completion is the set of Cauchy sequences in $\K$ of
length~$\kappa$,
modulo the set of null sequences.
\end{remark}

\begin{lemma}\label{lem:regular-type}
Let $\K$ be definably complete, $\K^* \succ \K$, $b \in \K^* \setminus \K$.
For every $X \subseteq \K^m$ definable, let $X^* \subseteq {\K^*}^m$ be the
interpretation of $X$ in $\K^*$.

Let $E \subset \K$ be closed, $\K$-definable, and with 
$n := \rkCB(E)  < \omega$, let $\Lambda$ be the cut of $\K$ determined by~$b$.
If $\Lambda$ is regular, then $b \notin E^*$.

If moreover $\K$ is d-minimal, and $D \subset \K$ is $\K$-definable and of
dimension~$0$, then $b \notin D^*$.
\end{lemma}
\begin{proof}
%That is, each element of $E^* \cap \KbC$  is $\K$-definable.
We will prove the lemma by induction on~$n := \rkCB(E)$.
\Wlog, $E$~is bounded (because $b$ is $\K$-bounded).
If $n = 0$, then $E = \emptyset$.
If $n = 1$, then $E$ is discrete; thus, $E$ is pseudo-finite; 
let $\delta := \delta(E)$; notice that $0 < \delta \in \K$. 
Since $\Lambda$ is regular, there exist $a', a'' \in \K$
such that $a' < b < a''$ and $a'' - a' < \delta$.
Thus, $(a', a'')^* \cap E^* = \set b$, and therefore $b$ is $\K$-definable, absurd.

If $n > 1$, let $G := \CBd E {n - 1}$; notice that $G$ is closed, discrete and
nonempty; thus, $\rkCB(G) = 1$, and $G$ is pseudo-finite.
By the case $n = 1$, $b \notin G^*$; let $a', a'' \in \K$ such that $a' < b <
a''$ and $G \cap [a', a''] = \emptyset$ ($a'$ and $a''$ exist by the proof of
the case $n = 1$).
Let $F := E \cap [a', a'']$.
Notice that $F$ is \dcompact.
Moreover, $\rkCB(F) < n$.
Therefore, by inductive hypothesis, $b \notin F^*$.
Hence, $b \notin E^*$.

If $\K$ is d-minimal, let $E$ be the closure of $D$ in~$\K$.
By d-minimality, $\rkCB(D) < \omega$; thus, $b \notin D^*$.
\end{proof}

From the proof of the above lemma, we can deduce the following.
\begin{lemma}\label{lem:regular-cut-dimension}
Let $\K$ be d-minimal, and $\K \prec \K^*$.
Let $c \in \K^* \setminus \K$, 
$\Lambda$~be the cut determined by $c$ over~$\K$,
and $X^*\subseteq \K^*$ be $\K$-definable.
If $\Lambda$ is a regular gap and $c \in X^*$, then there exists an open
interval $I$ with end-points in $\K$ such that $c \in I \subseteq X^*$.
\end{lemma}

\begin{lemma}\label{lem:regular-type-image}
Let $\K$ be definably complete, $\K^* \succ \K$, $b \in \K^* \setminus \K$,
$\Lambda$ be the cut of $\K$ determined by~$b$, 
$f: \K \to \K$ be definable,
$c := f^*(b) \in \K^*$, and $\Gamma$ be the cut determined by $c$ over~$\K$.
Assume that $\Lambda$ is a regular gap 
and $f$ is continuous and  strictly monotone.
Then, $\Gamma$~is a regular gap, and $\K \gen c = \K \gen b$,
where $\K \gen b$ is the definable closure of $\K \cup \mset b$.
\end{lemma}
\begin{proof}
% In fact, by Lemma~\ref{lem:uniform-cont}, \wloG we can assume that
% $b \in (0,1)$ and that $f$ is uniformly continuous and either constant or
% strictly increasing on $(0,1)$.
%If $f$ is constant, then $c \in \K$.
Let $0 < \varepsilon \in \K$; by uniform continuity, 
there exists $0 < \delta \in \K$, such that, if
$x, y \in \K$ and $\abs{x - y} < \delta$, 
then $\abs{f(x) - f(y)} < \varepsilon$.
%Since $\K$ is dense in~$\F$, and $b \in \F$, there exists 
Let $y', y'' \in \K$
such that $y' < b < y''$ and $y'' - y' < \delta$ (they exist because $\Lambda$
is a regular cut).
Thus, $f(y') <  c < f(y'')$ and $f(y'') - f(y') < \varepsilon$, 
and therefore $\Gamma$ is  regular.
Moreover, $f$ is invertible, and therefore $b = f^{-1}(c) \in \K \gen c$.
\end{proof}

\begin{lemma}
Let $\K$ be a d-minimal structure, and $\K^* \succ \K$, such that $\K$ is
dense in $\K^*$.
Then, the algebraic closure in $\K^*$ has the Exchange Property (EP) relative
to~$\K$; that is, if $A \subset \K^*$, and $b, c \in \K^*$ satisfy
$c \in \K\gen{A,b} \setminus \K\gen A$, then $b \in \K\gen{A,c}$.
\end{lemma}
\begin{proof}
Let $b$ and $c$ be as in the hypothesis.
Let $\K' := \K \gen A$; \wloG, $\K = \K'$.
Then, since $\K$ has definable Skolem function,
$b = f^*(c)$ for some $K$-definable $f: \K \to \K$.
By Theorem~\ref{thm:imin}-\ref{en:i-monotonicity}
there exists $D \subseteq \K$ nowhere dense and
$K$-definable, such that $f$ is continuous and either constant or strictly
monotone on each subinterval of $\K \setminus D$.
By Lemma~\ref{lem:regular-type}, $b \notin D^*$.
Therefore, by Lemma~\ref{lem:regular-type-image}, either $c \in \K$, 
or  $b \in \K \gen c$.
\end{proof}

Notice that the hypothesis that $\K$ is dense in $\K^*$ in the above lemma is
necessary: \cite{DMS}*{1.17} shows that if $\K^*$ does not satisfy UF and it is
sufficiently saturated, then the algebraic closure in $\K^*$ does not satisfy
(EP).

\cite{LS} proves the following proposition for o-minimal expansions of
Archimedean ordered groups, and \cite{Frecon} gives a proof for o-minimal
structures expanding a field.

\begin{proposition}[Cauchy completion]\label{prop:dmin-Cauchy}
Let $\K$ be a d-minimal structure, expanding the  ordered field~$\Kb$,
and $\Lang = (0, 1, +, \cdot, <, \dotsc)$ be the language of~$\K$.
%, such that $\K$ is a dense subfield of~$\F$.
There exists a unique expansion of the Cauchy completion
$\KbC$ to an $\Lang$-structure~$\KC$, such that
$\K$ is an elementary substructure of~$\KC$.
\end{proposition}

\begin{proof}
Let $\Sfam$ be set of elementary extensions $\K'$ of~$\K$, such that 
$\K$ is dense in~$\K'$; order $\Sfam$ by elementary inclusion.
Let $\K'$ be a maximal element of~$\Sfam$ ($\K'$~exists by Zorn's lemma).
We claim that $\Kb' = \KbC$ (and, thus, $\KbC$~can be expanded to an 
$\Lang$-structure).

Suppose not. \Wlog, $\K' = \K$.
Let $\Lambda$ be a regular gap of~$\K$, $\K^*$~be an elementary extension
of~$\K$, and $b \in \K^*$ fills the gap $\Lambda$.
Let $\K \gen b$ be the definable closure of $\K \cup \set b$ in~$\K^*$;
since $\K$ has definable Skolem functions (by definable choice),
$\K \prec \K \gen b \preceq \K^*$.

We claim that $\K$ is dense in $\K\gen b$, contradicting the maximality
of~$\K$.
In fact, let $c \in \K \gen b$; thus, $c = f(b)$ for some $\K$-definable
$f: \K^* \to \K^*$.
We have to prove that either $c \in \K$, or that $\Gamma$ is a regular gap,
where $\Gamma$ is the cut determined by $c$ over~$\K$.
By Theorem~\ref{thm:imin}\eqref{en:i-monotonicity}
there exists $D \subseteq \K$ nowhere dense and
definable, such that $f$ is continuous and either constant or strictly
monotone on each subinterval of $\K \setminus D$.
By Lemma~\ref{lem:regular-type}, $b \notin D^*$.
Therefore, by Lemma~\ref{lem:regular-type-image}, either $c \in \K$, 
or $\Gamma$~is a regular gap.

It remains to prove that the $\Lang$-structure on $\KbC$ is unique.
Let $\KC_1$ and $\KC_2$ be two expansion of $\KbC$ to elementary extensions
of~$\K$.
Let $\K'$ be a maximal common elementary substructure of $\KC_1$ and $\KC_2$
extending~$\K$.
Assume, for a contradiction, that $\K' \neq \KC$; 
\wloG, we can assume that $\K = \K'$.
Let $b \in \KC \setminus \K$, and let $\phi(x)$ be an $\Lang$-formula with
parameters in~$\K$.
In order to reach a contradiction, we must prove that $\KC_1 \models \phi(b)$
iff $\KC_2 \models \phi(b)$.
\Wlog, $\KC_1 \models \phi(b)$; let $X := \phi(\K)$.
Moreover, for every $Y \subseteq \K^n$ definable, 
let $Y^C_i$ be the interpretation of $Y$ in $\KC_i$, for $i = 1, 2$.
Notice that $X = U \cap D$, where $U := \inter X$ is open and definable, 
and $D := X \setminus \inter X$ is definable, with $\dim D = 0$.
Since $\K$ is d-minimal, $\rkCB(D) < \omega$.
Thus, by Lemma~\ref{lem:regular-type}, $b \notin D^C_1$.
Hence, $b \in U^C_1$; since $\K$ is dense in $\KC_1$, and $U^C_1$ is open,
there exist $y' < y'' \in \K$ such that $b \in (y', y'')^C_1 \subseteq U^C_1$.
Since $\K \preceq \KC_1$, $(y', y'') \subseteq U$,
and since $\K \preceq \KC_2$,
$b \in (y', y'')^C_2 \subseteq U^C_2 \subseteq X^C_2$.
\end{proof}

\begin{corollary}[of the proof]\label{cor:dmin-regular-cut}
Let $\K$ be a d-minimal structure, $\K^* \succ \K$,
$b \in \K^* \setminus \K$, and $\Lambda$ be the cut of $\K$ determined by~$b$.
If $\Lambda$ is a regular gap, then it uniquely determines the type of $b$
over~$\K$; moreover, $\K$ is dense in $\K \gen b$;
besides, for every $c \in \K \gen b \setminus \K$,
we have $\K \gen c = \K \gen b$. 
\end{corollary}

\begin{corollary}
An Archimedean locally o-minimal structure is o-minimal
\end{corollary}
\begin{proof}
If $\K$ is Archimedean, then $\KbC = \Rbar$; thus, $\K$ has an elementary
extension $\hat \Real$ that is an expansion of~$\Real$.
Therefore, $\hat \Real$ is o-minimal, and thus $\K$ is o-minimal.
\end{proof}

\subsection{Polish structures and theories}

\begin{lemma}\label{lem:Polish}
Let $\F$ be an ordered field.
Assume that $\F$ contains a countable dense subset (not necessarily
definable) and that $\F$ is Cauchy complete.
Then:
\begin{enumerate}
\item 
$\F$ has cofinality~$\omega$.
\item 
$\F$ is a Polish space (\ie, a Cauchy complete separable metric space).
\item\label{en:Polish-3}
$\card{\F} = 2^{\aleph_0}$.
\item 
If $X \subset \F^n$ is perfect, nonempty, and a $\Gd$ (in the topological
sense), then $\card{X} = 2^{\aleph_0}$.
\item
If $X \subset \F^n$ is closed and $\card{X} < 2^{\aleph_0}$, 
then $\isol(X)$, the set of isolated points of~$X$, is dense in~$X$.
\item
If $X \subseteq \F^n$ is a nonempty $\Gd$ (in the topological sense),
then it is a Baire space (again, in the topological sense), and it is even
strong Choquet (see~\cite{kechris}).
\end{enumerate}
\end{lemma}
\begin{proof}
\begin{enumerate}
\item is obvious.

\item requires us to define a metric.
If $\F$ is Archimedean, then $\F$ is homeomorphic to the reals, and we are
done.
Otherwise, let $v$ be the natural valuation on $\F$ induced by the ordering,
and $G$ be the value group of~$\F$.
\begin{claim}
The topology induced by $v$ on $\F$ is the same as the order topology.
\end{claim}
Notice that the claim is false if $\F$ is  Archimedean.

\begin{claim}
$G$ is countable.
\end{claim}
Thus, there exists a coinitial order-reversing embedding  $\iota$ of $(G, >)$ 
in $(\Real_+, >)$ %as ordered sets 
(notice that $\iota$ ignores the group structure).
For every $x, y \in \F$, define
\[
d(x, y) :=
\begin{cases}
\iota(v(x - y)) & \text{if } x \neq y;\\
0 & \text{otherwise}.
\end{cases}\]
\begin{claim}
$(\F, d)$ is a metric space.
\end{claim}
Actually, $(\F, d)$ satisfies the ultra-metric inequality.
\begin{claim}
$(\F,d)$ is homeomorphic to $(\F, v)$ (and hence to $(\F, <)$).
\end{claim}
\begin{claim}
If $(a_n)_{n \in \Nat}$ is a Cauchy sequence in $(\F, d)$,
then $(a_n)_{n \in \Nat}$ is a Cauchy sequence in $(\F, v)$.
\end{claim}
Hence, assume that $(a_n)_{n \in \Nat}$ is a Cauchy sequences in $(\F, d)$.
Since $(\F, v)$ is Cauchy complete (by assumption), $a_n \to a$ for some 
$a \in  \F$, according to the topology induced by~$v$.
However, $v$ and $d$ induce the same topology, and therefore $a_n \to a$
according also to~$d$; thus, $(\F, d)$ is a complete metric space.
Finally, $(\F, d)$ is separable by assumption.

\item %$\card{\F} \leq 2^{\aleph_0}$ is easy.
%The opposite inequality 
Follows from (2) and the fact that if $Y$ is a
nonempty perfect Polish space, then $\card Y = 2^{\aleph_0}$
(see \cite{kechris}*{6.2}).

\item
Notice that a $\Gd$ nonempty subset of $\F^n$ is a Polish
space (see \cite{kechris}*{3.11}). 
Thus, $X$~itself is a nonempty perfect Polish space,
and therefore $\card X = 2^{\aleph_0}$.
%The conclusion follows from~\cite[6.2]{kechris}.

\item Assume, for a contradiction, that $\isol(X)$ is not dense in~$X$;
let $B \subseteq \K^n$ be a closed box, such that $\inter{B} \cap X \neq
\emptyset$, and $Y := X \cap B$ contains no isolated points.
Hence, $Y$~satisfies the hypothesis of~(4), absurd.

\item Every Polish space is Baire and strong Choquet 
(see \cite{kechris}*{8.17}).
\qedhere
\end{enumerate}
\end{proof}

\begin{definizione}
%Let $\F$ be a structure expanding an ordered field.
%We say that $\F$ is a Polish structure (or a Polish model) if, with the
%topology induced by the ordering, $\F$~is a Polish space.
Let $T$ be a complete theory expanding the theory of ordered fields,
in a language~$\Lang$, expanding the language $\Langf$ of ordered fields.
We say that $T$ is a \intro{Polish theory} if, for every finite
language~$\Lang'$,  
such that $\Langf \subseteq \Lang' \subseteq \Lang$, the restriction of
$T$ to $\Lang'$ has a model which is separable and Cauchy complete.
%Polish model.
If $T$ is not complete, we say that $T$ is a Polish theory if every completion
of $T$ is Polish.
\end{definizione}

\begin{lemma}
Let $T$ be a definably complete theory (expanding RCF).
If $T$ is Polish, $\K \models T$ and $X \subseteq \K^n$ is definable
and $\Gd$ (in the definable sense), then $X$ is definably Baire.
\end{lemma}
\begin{proof}
\Wlog, $\K$ is Cauchy complete and separable.
Hence, by Lemma~\ref{lem:Polish}, $X$ is topologically Baire, 
and \emph{a fortiori} definably Baire.
\end{proof}

\begin{proposition}\label{prop:dmin-polish}
A \dminimal theory $T$ is Polish.
In particular, if $\K$ is \dminimal and $X \subseteq \K^n$ is definable,
then $X$ is definably Baire.
\end{proposition}
\begin{proof}
\Wlog, the language of $T$ is countable.
Let $\F'$ be a countable model of~$T$, and $\F$ be the Cauchy completion
of~$\F'$. 
Therefore, $\F$ is a model of~$T$ and a Polish space.
%, and, by lemma~\ref{lem:Polish}, it is a Polish space.
\end{proof}

The above proposition may be a step in the proof of the following conjecture.

\begin{conjecture}
Let $X \subseteq \K^n$ be definable and closed.
Then, $X$ is definably Baire.
\end{conjecture}

% ------------------

\section[Dense pairs]{Dense pairs of d-minimal structures}\label{sec:dense}

\subsection{The Z-closure}
\label{subsec:zclosure}

\begin{definizione}
Let $A \subseteq \K$ and $c \in \K$.
We define the \Zclosure of $A$ inside~$\K$
\[
\zclK(A) := \bigcup\set{C \subset \K: C \text{ nowhere dense and definable with parameters from } A}.
\]
If $\K$ is clear from the context, we drop the superscript~$\K$.
$A \subseteq \K$ is \Zclosed in $\K$ if $\zclK(A) = A$.
\end{definizione}
The notion above is most interesting when $A$ is an elementary substructure 
of~$\K$.

\begin{remark}
$\dcl \subseteq \zcl$. Moreover, if $\K$ is o-minimal, then $\zcl = \dcl$.
\end{remark}

\begin{remark}
If $\K$ is o-minimal and $A \subseteq \K$, then $A$ is \Zclosed in $\K$ if
and only if $A$ is an elementary substructure of~$\K$.
\end{remark}

% \begin{remark}
% For every $A \subseteq \K$, $\dcl(A) \subseteq \zcl(A)$, and
% $\zcl(\dcl(A)) = \zcl(A) = \dcl(\zcl(A))$.
% \end{remark}

\begin{remark}
If $\K$ has \DSF, then $\zclK(A) \preceq \K$.
\end{remark}

\begin{lemma}\label{lem:Z-interior-generator}
If $A \subseteq \K$ has nonempty interior, then $\dcl(A) = \K$, and therefore
$\zcl(A) = \K$.
\end{lemma}
\begin{proof}
Since $A \subseteq \dcl A = \dcl(\dcl(A))$, \wloG $A = \dcl A$.
Let $\varepsilon > 0$ and $a \in \K$
such that $B(a; \varepsilon) \subseteq A$.
Hence, $(-\varepsilon, \varepsilon) \subseteq A$.
Thus, $(1 /\varepsilon, +\infty) \subseteq A$.
Let $b \in \K$; we want to prove that $b \in A$; \wloG, $b > 0$.
Let $a := (1/\varepsilon)$ and $a' := b + 1/\varepsilon$;
notice that $a$ and $a'$ are in~$A$, and therefore $b = a' - a \in A$.
\end{proof}

\begin{remark}
Let $\K \prec \F$ be a dense substructure.
If $\K$ is \dminimal, then $\K$ is \Zclosed in~$\F$.
\end{remark}
\begin{proof}
By Lemma~\ref{lem:regular-type}.
\end{proof}

\begin{remark}
Given $A \subseteq \K$, $\dcl(A)$ does not depend on~$\K$:
that is, if $\K \preceq \K'$, then the definable closure of $A$ inside $\K$
and the definable closure of $A$ inside $\K'$ are the same set.
Instead,  $\zclK(A)$ may depend on~$\K$:
for instance, if an infinite nowhere dense subset of $\K$ is $A$-definable
and if $\K'$ is a $\kappa$-saturated elementary extension of~$\K$,
then $\card{\zcl^{\K'}(A)} \geq \kappa$.
If $\K$ is \dminimal but not o-minimal, then $\zcl$ does depend on~$\K$:
for instance there exists some $\K' \succ \K$ such that
$\zcl^{\K'}(\emptyset) \neq \zclK(\emptyset)$.
\end{remark}

\begin{remark}\label{rem:zcl-extension}
If $A \subseteq \K \preceq \K'$, then $\zclK(A) = \zcl^{\K'}(A) \cap \K$.
\end{remark}

\begin{definizione}
Let $f: X \app Y$ be a definable application
(\ie, a multi-valued partial function), with graph~$F$.
Assume that $\K$ is \iminimal.
For every $x \in X$, let $f(x) := \set{y\in Y: \pair{x,y} \in F}\subseteq Y$.
Such an application $f$ is a \intro{\Zapplication{}} if, for every $x \in X$, 
$\dim\Pa{f(x)} = 0$ (thus, the domain of $f$ is all~$X$);
it is a partial \Zapplication if for every $x \in X$, 
$\dim\Pa{f(x)} \leq 0$.
\end{definizione}

\begin{remark}
Let $A \subseteq \K$, and $b \in \K$.
Then, $b \in \zcl A$ iff there exists an $\emptyset$-definable
\Zapplication $f: \K^n \app \K$ and $\av \in A$,
such that $b \in f(\av)$.
Moreover, if $\cv \in \K^n$, then $b \in \zcl(A\cv)$ iff there exists
an $A$-definable \Zapplication $f: \K^n \to \K$, such that $b \in f(\cv)$.
\end{remark}
\begin{proof}
The ``if'' direction is clear: $f(\av)$ is nowhere dense.
For the converse, let $Z \subset \K$ be nowhere dense and $A$-definable,
such that $b \in Z$.
Let $\phi(x,\av)$ be the formula defining~$Z$.
Let $\psi(x,\y)$ be the formula ``($\psi(\x,\y)$ and $\psi(\K, \y)$ is 
nowhere dense) or ($x = 0$ and $\phi(\K, \y)$ is somewhere-dense)''.
Then, $\psi$ defines a \Zapplication $f: \K^n \app \K$,
and $b \in f(\av)$.

The ``moreover'' part is clear.
\end{proof}

\begin{proviso}
For the rest of this section, $\K$ is \dminimal,
and $\monster \succeq \K$ is ``the'' monster model.
%$\K$ is \iminimal with \DSF.
\end{proviso}
%We will show that, under the above condition, $\zclK$ is a matroid
%(\aka combinatorial pregeometry); we will write $\zcl$ for $\zclK$.

In \cite{fornasiero-matroids}*{\S9}, we defined \dminimal topological structures.
We will now prove that $\K$ is such a structure.

\begin{proposition}
$\K$ is a \dminimal topological structure, in the sense of 
\cite{fornasiero-matroids}*{Definition~9.1}.
\end{proposition}
\begin{proof}
The only nontrivially true conditions in
\cite{fornasiero-matroids}*{Definition~9.1} are:
\begin{sentences}
\item[(4)] For ever $X \subseteq \K^n$ definable and discrete, $\Pi^n_1(X)$ has empty interior;%\\[1ex]
%\smallskip
%\mbox{}\hspace{-\parindent}and
%\smallskip
\item[(5)] Given $X \subseteq \K^2$ and $U \subseteq \Pi^2_1(X)$ definable sets, if $U$ is open and nonempty, and $X_a$ has nonempty interior for every $a \in U$, then $X$ has nonempty interior.
\end{sentences}
Notice that (4) is immediate from Fact~\ref{fact:enumerable}, and
(5) follows from the fact that $\K$ is \ipminimal and Kuratowski-Ulam Theorem.
\end{proof}

%Thus, by \cite[Theorem~9.8]{fornasiero-matroids}, $\Zcl$ is an \intro{existential matroid}.

\begin{corollary}
%If $\K$ is \dminimal, then 
$\zcl$ is an existential matroid.
%\Ie, let $\av \in \monster^n$, $\bv \in \monster^m$, $c \in \monster$,
%and $\phi(x,\y,\z)$ be \xnarrow.
%Assume that, for every conjugate $c'$ of $c$ over $\av$, 
%$\monster \models \phi(c', \av, \bv)$.
%Then, $c$ (and all its conjugates over~$\av$) is in $\zclM(\av)$.
\end{corollary}
\begin{proof}
% Let $\Xi(c/\av)$  be the set of conjugates of $c$ over~$\av$.
% If, for a contradiction, $c \notin \zclM(\av)$,
% then, by the above lemma, $\Xi(c/A)$ has nonempty interior.
% However, $\Xi(c/A) \subseteq \phi(\monster, \av, \bv)$, and the latter
% is nowhere dense, absurd.
It is \cite{fornasiero-matroids}*{Theorem~9.8}.
%The corollary is also a direct consequence of 
%\cite[Lemma~3.22]{fornasiero-matroids}.
\end{proof}

Notice that if $\K$ is o-minimal, then $\zcl = \dcl$, and therefore
$\zind = \tind$, where $\zind$ is the independence relation induced by~$\zcl$.
The converse is also true 
(remember the assumption that $\K$ is \dminimal):
\begin{lemma}
\Tfae:
\begin{enumerate}
\item $\zclM = \dcl$;
\item $\zind = \tind$ \rom(in the monster model\rom);
\item $\K$ is o-minimal.
\end{enumerate}
\end{lemma}
\begin{proof}
($3 \Rightarrow 1$) is clear.
If (1) holds, then $\zind = \mind$ (where $\mind$ was defined in~\cite{adler}); 
moreover, since $\zind$ is symmetric,
$\mind$ is also symmetric, and therefore $\mind = \tind$ 
(see \cite{adler}*{Theorem~2.39}).

Assume that (2) holds.
Let $a \in \zcl(B)$.
Then, $a \zind_B a$, therefore $a \tind_B a$, and thus
$a \in \dcl B$: hence, (1) also holds.
We have to prove that $\monster$ is o-minimal.
Assume, for a contradiction, that $A \subset \monster$ is definable with 
parameters~$\bv$, infinite and with empty interior.
Then, $A \subseteq \zclM(\bv) = \dcl(\bv)$.
However, since $A$ is infinite, $\card{A} \geq \kappa$, and therefore
$\card{\dcl(\bv)} \geq \kappa > \card T$, which is impossible.
\end{proof}

\begin{remark}\label{rem:Zrk}
The dimension induced by $\zcl$ and the geometric notion of dimension
coincide.
That is, if $X \subseteq \monster^n$ is definable, then
$\dim X = \max \set{\RK^{Z}(\x): \x \in X}$ (where $\RK^Z$ is the rank function induced by~$\zcl$).
\end{remark}

% --------------------------------------

Contrast the situation of $\zind$ to the notion of $M$-dividing independence
(defined in~\cite{adler}), where, 
$A \mind_B C$ iff, 
for every $\cv \subset \dcl(B C)$, 
\[
\dcl(A B \cv) \cap \dcl(B C) = \dcl(B \cv).
\]

\begin{lemma}\label{lem:asymmetry}
Assume that $T$ is \dminimal, but \emph{not} o-minimal.
Then, $\mind$~is not symmetric
\rom(and therefore $\dcl$ does not have the Exchange Property\rom).
% However, $\mind$~does satisfy the existence and extension axioms, 
% and therefore coincides with~$\tind$, the \th-forking relation.
% Hence, $T$~is not rosy,
% and in particular $\tind$ is not symmetric.
\end{lemma}
We don't know whether $T$ might be rosy or not.
\begin{proof}
%The fact that $\mind$ satisfy existence and extension and coincides with 
%$\tind$ is immediate from \cite[Lemma~3.22]{fornasiero-matroids}.
%
Let $\K \prec \monster$ such that $\K$ is not Cauchy complete.
%and $\card{\K} < \kappa$.
By expanding the language by $\card{\K}$ new constants, \wloG we can
assume that $\K$ is the prime model of~$T$.
Let $\pi$ be a regular gap of~$\K$; choose $c_0$ and $c_1$
such that $c_1 > \K$ and $c_0 \models \pi$.
Let $\cv \coloneqq \pair{c_0,c_1}$.
Since $\K$ is \dminimal but not o-minimal, there exists an infinite
pseudo-finite set $X$ definable over~$\K$.
By compactness, there exists $b \in X \setminus \dcl(\cv)$.
\Wlog, we can assume that $b > 0$.
Let $\delta \coloneqq \delta(X)$; since $X$ is pseudo-finite and definable
over~$\K$, $0 < \delta \in \K$, and $b$ is bounded by some element of~$\K$.
Let $a_0 \in B(b; \delta/4) \setminus \zcl(b \cv)$
($a_0$ exists because $\zcl$ is an existential matroid).
Define $a \coloneqq a_0/c_1 + c_0$.

% \begin{claim}
% There exist $a_0, b \in \monster$, and a finite set $C$ such that:
% \begin{enumerate}
% \item $c_0, c_1 \in C$;
% \item $b \in \zcl(C)$;
% \item $b \in \dcl(C a_0) \setminus \dcl(C)$;
% \item $a_0 \notin \zcl(Cb)$;
% \item $a_0 > 0$ and $a_0$ is infinitesimal \wrt~$\K$.
% \end{enumerate}
% \end{claim}
% Let $X$ be a definable (with parameters) subset of $\monster$ which is
% discrete and infinite (such a set exists by Corollary~\ref{cor:i-o-min},
% because $T$ is \iminimal but not o-minimal), and let $C$ be any finite set
% containing $c_0$, $c_1$, and the parameters of~$X$. 
% By saturation, there exists $b \in X \setminus \dcl(C)$.
% Let $I$ be an open interval containing~$b$, such that $I \cap X = \set b$.
% By Lemma~\ref{lem:Z-interior-generator}, 
% since $\Zrk(\monster /\K) \geq \kappa$, 
% there exists $a_1 \in I \setminus \zcl(Cb)$.
% \Wlog, $a_1 > 0$.
% If $a_1$ is infinitesimal (\wrt~$\K$), let $a_0 := a_1$;
% if $a_1$ is finite but not infinitesimal, let $a_0 := a_1 / c_1$;
% if $a_1$ is infinite, let $a_0 := 1/ a_1$.
% Notice that $\zcl(C b) = \zcl(C)$.
% \begin{claim}
% With $b$ and $C$ as in the above Claim, there exists $a \in \monster$ such
% that:
Hence, we have the following properties:
\begin{enumerate}
\item $a \models \pi$ (and therefore $a$ is in a Cauchy completion of~$\K$);
\item $a \notin \zcl(\cv b)$;
\item $b \in \dcl( \cv a) \setminus \dcl(\cv)$.
\end{enumerate}
%\end{claim}
% Let $a := c_0 + a_0$.
% Since $a_0$ is infinitesimal, $a \models \pi$.
% Since $a_0 \notin \zcl(C) = \zcl(C b$ and $c_0 \in C$, the second point
% follows. 
% Since $\dcl(C a) = \dcl(C a_0) \ni b$, the third points follows.
\begin{claim}
$a \notind[M]_{\cv} b$.
\end{claim}
$b \in \dcl(\cv a) \cap \dcl(\cv b)$, but $b \notin \dcl(\cv)$.
\begin{claim}
$a \notind[M] \cv b$.
\end{claim}
Follows immediately from the previous claim and transitivity for $\mind$.
\begin{claim}
$\cv b \mind a$.
\end{claim}
In fact, let $A := \dcl(a) = \K \gen a$, and $A' \subseteq A$.
Define $Y := \dcl(\cv A' b) \cap A$; we have to prove that $Y = A'$.
Since $a$ satisfies a regular gap over~$\K$, $\dcl(\emptyset) = \K$ is dense
in~$A$; therefore, $\dcl^A$ satisfies EP.
Hence, either $A' = \K$, or $A' = A$.
If $A' = A$, the conclusion is obvious.
If $A' = \K$, then $a \notin Y$, because $Y \subset \zcl(\cv b)$, 
and $a \notin \zcl(\cv b)$; therefore, since $\dcl^A$ satisfies EP, $Y = \K$.
\end{proof}

We do not know if the above lemma extends to \iminimal theories with \DSF,
or to \ipminimal theories.

% ------------------------------------

\subsection{Dense pairs}
\label{subsec:dense}

Dense pairs of o-minimal structures were studied in~\cite{vdd-dense}.
Dense pairs of \dminimal topological structures were studied in~\cite{fornasiero-matroids}.
%we can apply the results there to our situation.

\begin{proviso}
Remember that $\K$ is \dminimal.
Let
%$\K$ is an \ipminimal structure with \DSF, and 
$T \coloneqq \Th(\K)$.
\end{proviso}

We have seen that the \Zclosure is an existential matroid on~$\K$.
Moreover, $A \subseteq \K$ is topologically dense iff it s dense \wrt to the
matroid $\zcl$, that is iff $X$ intersects every definable subset of $\K$ of
dimension~$1$.

\begin{definizione}
Let $\Td$ be the theory of pairs $\Am \prec \Bm \models T$, 
such that $\Am$ is dense in~$\Bm$. %and $\zcl(\Am) = \Am$.
More generally, for every $n \in \Nat$,
let $T^{n d}$ be the theory of tuples
$\Am_0 \prec \Am_1 \prec \dots \Am_n \models T$,
such that each $\Am_{i+1}$ is a proper elementary extension of $\Am_i$,
and $\Am_0$ is dense in $\Am_n$.
\end{definizione}
Notice that $T^{1 d} = \Td$.
We can apply the results in \cite{fornasiero-matroids} to~$T$, and obtain
the following results.

\begin{thm}
For every $n \in \Nat$,
$T^{n d}$ is consistent and complete.
Besides, $\Am_n$ is the open core of 
$\pair{\Am_0 \prec \Am_1 \prec \dots \Am_n} \models T^{n d}$.
Any model of $T^{nd}$ is definably complete.
%and $\zcl^{\Bm}(\Am) = \Am$. 
%If moreover $T$ is \dminimal, then $\Td$ is the theory of 
%pairs $\Am \prec \Bm \models T$, such that $\Am$
%is dense in $\Bm$ (the fact that $\Am$ is $\zcl$-closed in $\Bm$ follows).
\end{thm}
\begin{proof}
By \cite{fornasiero-matroids}*{Theorems~13.5 and~13.11}.
\end{proof}
%Similar results can be shown for dense tuples of models of $T$
%\cite[\S13]{fornasiero-matroids}.
More results can be proved for~$\Td$ (\eg near model-completeness: see
\cite{fornasiero-matroids}*{Theorem~8.5}).
%If moreover $T$ is \dminimal, then also the results in~\cite[\S9]{fornasiero-matroids} apply to~$T$.

We will give some additional results and conjectures that are more specific to
our situation.

\begin{thm}[\cite{vdd-dense}*{Theorem~2}]
Let $\pair{\Bm, \Am} \models \Td$.
Given a set $Y \subset \Am^n$, \tfae:
\begin{enumerate}
\item 
$Y$ is definable in $\pair{\Bm, \Am}$;
\item 
$Y = Z \cap \Am^n$ for some set $Z \subseteq \Bm^n$ that is definable in~$\Bm$.
%If moreover $T$ is \ipminimal, then the above two conditions are equivalent
\item 
$Y$ is definable in the structure $\pair{\Am, (\Am \cap (0,b))_{b \in \Bm}}$.
\end{enumerate}
\end{thm}

\begin{proof}
$(1 \Rightarrow 2)$ and $(2 \Rightarrow 3)$
are as in \cite{vdd-dense}*{Theorem~2}.
$(3 \Rightarrow 1)$ and $(2 \Rightarrow 1)$ are obvious.
\end{proof}

\begin{question}
Let $T'$ be a complete theory with locally o-minimal open core.
Is there an existential matroid on~$T'$?
\cite{DMS}*{6.2} and \cite{fornasiero-matroids}*{\S8.4}
prove that if $T$ is equal to either
$T^d$ or $T^g$ (see \cite{DMS} for the definition of $T^g$) 
for some o-minimal theory~$T$, then $T'$ admits such a matroid
(in the case of $T^g$, the matroid is~$\acl$).
\end{question}

%\begin{example}
%We show that $\zcl$ does not satisfy Existence on~$\monster \coloneqq \pair{\Bm, \Am}$, a monster model of~$\Td$. 
%Write $x' \elem^2_C x$ if the $\Ltwo$-type of $x$ and $x'$ over $C$ are the
%same, and let $\Xi^2(x / C) := \set{ x' \in \monster: x' \elem^2_C x}$.
%
%Choose $a_1, a_2 \in \Am$ and $b \in \Bm \setminus \Am$, such that
%$a_1$ and $a_2$ are \Zindependent over~$b$.
%Let $c := a_1 b + a_2$, and $f: \Am^2 \to \Bm$ be the definable function
%$f(x_1, x_2) := x_1 \cdot c + x_2$.
%By hypothesis, $\zrk(a_1 a_2/ b) = 2$, and therefore
%$\zrk(a_1, a_2 / b c) \geq 1$.
%Thus, either $a_1 \notin \zcl(b c)$, or $a_2 \notin \zcl(b c)$;
%\wloG, $a_1 \notin \zcl(b c)$.
%However, $f$ is injective, and therefore 
%$a_1$ and $a_2$ are $\Ltwo$-definable over $b c$,
%hence, $\Xi^2(a_1/ b c) = \set{a_1} \subseteq \zcl(b c a_1)$.
%If $\zcl$ did satisfy existence, then $a_1 \in \zcl(b c)$, absurd.
%\end{example}

\subsection{The small closure}
\begin{proviso}
Remember that $T$ is a \dminimal complete theory.  
Let $\Cm := \pair{\Bm, \Am} \models \Td$.
\end{proviso}

We have seen that $\Bm$ is the open core of~$\Cm$.
Hence, since every $\Fs$ subset of $\Cm^n$ is definable in the open core
of~$\Cm$, every such set is constructible.
%We will prove some additional results about this topic.
$\scl$ is the \intro{small closure} on $\Cm$ and $\sdim$ is the
corresponding dimension function, as defined in
\cite{fornasiero-matroids}*{\S8.4}.
I recall that a $\Cm$-definable set $X$ is called \intro{small} if $\sdim(X)
\leq 0$.

\begin{remark}\label{rem:s-dim}
Let $X \subseteq \Cm^n$ be $\Bm$-definable.
Then, $\sdim(X) = \dim(X)$.
\end{remark}
\begin{proof}
\cite{fornasiero-matroids}*{Lemma~8.31}.
\end{proof}

\begin{lemma}\label{lem:cl-U}
Let $(X_t)_{t \in \Cm}$ be a definable increasing family of subsets
of~$\Cm^n$, and $X := \bigcup_t X_t$.
Let $d \leq n$ and assume that, for every $t \in \Cm$,
$\sdim(X_b) \leq d$.
Then, $\sdim(X) \leq d$.
\end{lemma}
\begin{proof}
\cite{fornasiero-matroids}*{Lemma~3.71}, applied to $\scl$.
\end{proof}

% \begin{corollary}
% $\Cm$ is definably Baire.
% \end{corollary}
% \begin{proof}
% By Remark~\ref{rem:s-dim},
% if $X \subseteq S$ is definable in
% $\Cm$ and nowhere dense, then $\sdim(X) = 0$.
% The conclusion then follows from Lemma~\ref{lem:cl-U}.
% \end{proof}

\begin{proposition}\label{prop:small-polish}
\begin{enumerate}
\item $\Td$ is a Polish theory.
\item\label{en:polish-pair} 
Assume that $T$ is countable.
Then, there exists $\pair{\Bm', \Am'} \models \Td$ such that:
\begin{itemize}
\item $\Am'$ is countable;
\item $\Bm'$ is a separable complete metric space;
\item for every $X \subseteq \Bm'^n$ definable in $\pair{\Bm', \Am'}$,
$\sdim(X) \leq 0$ iff $X$ is countable.
\end{itemize}
\end{enumerate}
\end{proposition}
\begin{proof}
It suffices to prove (2): therefore, we can assume that $T$ is countable.
%\Wlog, the language of $T$ is countable.
Let $\Am'$ be a countable model of $T$ and $\Bm'$ be its Cauchy completion.
Notice that $\Bm' \neq \Am'$ and therefore $\pair{\Bm', \Am'}$
is a Cauchy complete and separable model of~$\Td$ (and therefore $\Bm'$ is a
Polish space).
Hence, by Lemma~\ref{lem:Polish}\eqref{en:Polish-3}, $\card{\Bm'} = 2^{\aleph_0}$.
%It is easy to see that $\card{\Bm'} = 2^{\aleph_0}$ (one of the usual proofs
%that $\card{\Real} = 2^{\aleph_0}$ will work).

\begin{claim}\label{cl:small}
Let $Y \subseteq \Bm'$ be definable in $\Bm'$ and of dimension~$0$.
Then, $Y$~is countable.
\end{claim}
In fact, by \dminimality, $Y$ is a finite union of discrete sets, 
and in a Polish space every discrete subset is countable.

\begin{claim}\label{cl:big}
Let $Y \subseteq \Bm'^n$ be definable in $\Bm'$ and of dimension at least~$1$.
Then, $\card Y = 2^{\aleph_0}$.
\end{claim}
In fact, after a permutation of coordinates, $\Pi^n_1(Y)$ will contain an open
interval. 

Let $X \subseteq \Bm'^n$ be definable in  $\pair{\Bm', \Am'}$.

If $\sdim(X) = 0$, then, by \cite{fornasiero-matroids}*{Lemma~8.33},
there exists a $\Bm'$-definable 
\Zapplication $f: \Bm'^m \app \Bm'^n$ such that $X \subseteq f(\Am'^m)$.
By the Claim~\ref{cl:small}, $f(\Am'^m)$ is countable, and therefore $X$ is countable.

If $\sdim(X) > 0$, then, after a permutation of coordinates,
$\sdim(Z) = 1$, where $Z \coloneqq \Pi^n_1(X)$.
By \cite{fornasiero-matroids}*{Proposition~8.36},
there exists $Y \subseteq \Bm'$ such that $Y$ is $\Bm'$-definable,
and $\sdim(Z \sdiff Y) \leq 0$.
Thus, $\dim(Z) = 1$.
By the previous case, $Z \sdiff Y$ is countable, and, by Claim~\ref{cl:big},
$\card Z = 2^{\aleph_0}$.
Thus, $\card Y = 2^{\aleph_0}$.
\end{proof}

% --------------------

\section{Open covers}

In this section $C$ will be a topological space and $\Vfam = (V_i: i \in I)$ a
family of open subsets of~$C$ (indexed, possibly with repetitions, by some set~$I$), 
covering $C$ and with each element in $\Vfam$ nonempty.
We say that: 
$\Vfam$ is pointwise finite (resp., pointwise countable) if, for every
$x$ in $C$, the set $\set{i \in I: x \in V_i}$ is finite
% is contained in at most finitely 
%many sets in~$\Vfam$ (counted with multiplicity); 
(resp., countable); 
$\Vfam$ is pointwise uniformly finite if 
there exists $k \in \Nat$ such that, for  every $x$ in~$C$,
the set $\set{i \in I: x \in V_i}$ has cardinality at most~$k$;
%s contained in at most $k$ many sets in~$\Vfam$;
$\Vfam$ is locally finite if every $x \in C$ has an open
neighborhood~$U$, such that the set
$\set{i \in I: U \cap V_i \neq \emptyset}$ is finite
%U$ meets finitely many sets in~$\Vfam$ 
(and similarly for uniformly locally finite);
finally, $\Vfam$ is finite (resp., infinite) 
if $I$ is finite (resp., infinite).
If $C$ is a definable set, we say that $\Vfam$ is pseudo-finite if $\Vfam$ is a definable family,
and the index set $I$ is pseudo-finite; the definitions of pseudo-enumerable,
pointwise pseudo-finite, etc. are analogous.

We will now give some properties of definable covers of definable sets.
The main results are propositions~\ref{prop:dmin-cover}
and~\ref{prop:pair-cover}, which also show some nice applications of
propositions~\ref{prop:dmin-polish} and~\ref{prop:small-polish}.

First, some observations in the ``topological'' setting.

\begin{enumerate}
\item
$C$ is compact iff every open cover $\Vfam$ of $C$ has a finite subcover.  
\item
$C$ is compact iff, for every $\Vfam$ locally finite open cover of~$C$,  
$\Vfam$ is finite.  
\item
There exists $\Vfam$ open cover of~$\Real$, such that $\Vfam$ is infinite but
uniformly locally finite 
(take \eg $\Vfam \coloneqq \set{(n-1, n+1): n \in \Nat}$).
\item 
There exists $\Vfam$ open cover of $[0, 1]$, such that $\Vfam$ is infinite but
uniformly pointwise finite 
(take $\Vfam \coloneqq \set{[0,1]} \cup \set{(\frac 1{n+1}, \frac 1 n): 0 < n \in \Nat}$).
\end{enumerate}

\begin{fact}\label{fact:polish-cover}
Assume that $C$ is separable.
Let $Q \subseteq C$ be dense and countable.
%contained in a Polish space and
Assume that every point of $Q$ is contained in at most countably many elements
of~$\Vfam$.
%$\Vfam$ is pointwise countable.
Then, $\Vfam$~is \rom(at most\rom) countable.
In particular, if $C$ is a subset of a Polish space and $\Vfam$ is pointwise
countable, then $\Vfam$ is countable.
\end{fact}

Assume now that $C$ is definable, and every $V \in \Vfam$ is also definable.

\begin{enumerate}
\setcounter{enumi}{4}
\item 
If $\K$ is nonarchimedean and $C$ is the interval $[0,1]$ in $\K$, then
there exists $\Vfam$ definable  and with no finite subcover
(take $\varepsilon > 0$ infinitesimal, and let
$\Vfam \coloneqq \set{(x - \varepsilon, x + \varepsilon): x \in [0,1]}$).
\item
If $C$ is the interval $[0,1]$ in $\K$, then
there exists $\Vfam$ covering $C = [0,1]$, such
that $\Vfam$ is uniformly pointwise finite, but $\Vfam$ is infinite 
(use the same cover as in (4)).
%(take \eg
%$\Vfam \coloneqq \set{[0,1]} \cup \set{(1/n+1, 1/n): 0 < n \in \Nat}$).
\end{enumerate}

The following remark is a definable version of Heine-Borel Theorem.
\begin{remark}\label{rem:subcover}
Assume that $\Vfam$ is at most pseudo-enumerable and $C$ is \dcompact.
Then, $\Vfam$~has a pseudo-finite subcover.
\end{remark}
\begin{proof}
%If $I$ (and hence~$\Vfam$) is not already pseudo-finite, then, \wloG, 
\Wlog, $I$~is a definable, closed, and discrete  subset of $\K_{\geq 0}$.
If, for a contradiction, $\Vfam$ has no pseudo-finite subcover, then
for every $i \in I$, the set $C \setminus \bigcup_{j \leq i} V_j$ is
closed in~$C$ and nonempty.
Since $C$ is \dminimal, $C \setminus \bigcup_{j \in I} V_j$ is also nonempty, absurd.
\end{proof}

\begin{lemma}\label{lem:compact-cover}
Assume that $\Vfam$ is definable and locally pseudo-finite.
If $C$ is \dcompact, then $\Vfam$ is pseudo-finite.
If $C$ is an open subset of $\K^n$, then $\Vfam$ is at most pseudo-enumerable.
\end{lemma}
\begin{proof}
Assume that $C$ is \dcompact.
For every definable subset $D \subseteq C$, let $P(D)$ be the property
``$D$ intersects only pseudo-finitely many sets in $\Vfam$''.
Then, $P$ is definable, monotone and additive 
(see \cite{fornasiero-lomin}*{Definitions~4.1 and~5.6}).
Hence, by \cite{fornasiero-lomin}*{Lemma~5.7}, $P(C)$ holds.

If $C$ is open in $\K^n$, then $C = \bigcup_{t \in \K} C_t$, where $(c_t: t
\in \K)$ is an increasing definable family of \dcompact sets.
Each $C_t$ intersects only pseudo-finitely many sets in~$\Vfam$,
hence $\Vfam$ is at most pseudo-enumerable.
\end{proof}

\begin{remark}
Let $\Vfam$ be definable.
If $\K$ is locally o-minimal
and $\Vfam$ is pointwise pseudo-finite, then $\Vfam$ is pseudo-finite.
In particular, if $\K$ is o-minimal and $\Vfam$ is  pointwise finite, then $\Vfam$ is finite.
\end{remark}
\begin{proof}
See Proposition~\ref{prop:dmin-cover}(1) (since pseudo-enumerable and pseudo-finite
coincide when $\K$ is locally o-minimal).
\end{proof}

Notice that in the above remark we did \emph{not} assume that $C$ is
\dcompact: therefore, in the o-minimal and locally o-minimal cases, what
\emph{prima facie} would seems a  property of \dcompact sets,
is true instead for every definable set.
For \dminimal non locally o-minimal structures instead, a converse
of Lemma~\ref{lem:compact-cover} holds: 
see  Proposition~\ref{prop:dmin-cover}(2).
%
%As the next lemma shows, when $\K$ is \dminimal but not locally o-minimal,
%is it possible to distinguish between \dcompact sets and non \dcompact ones.
%
The assumption in Proposition~\ref{prop:dmin-cover}(2) is equivalent to the
fact that the open core of $\K$ is not locally o-minimal (see
\cite{fornasiero-lomin}*{Theorem~A}): 
in particular, it holds if $\K$ is \dminimal but not locally o-minimal.

\begin{proposition}\label{prop:dmin-cover}
Assume that $\Vfam$ is definable.
\begin{enumerate}\item 
If $\K$ is \dminimal and 
$\Vfam$ is pointwise at most pseudo-enumerable, then $\Vfam$ is
at most pseudo-enumerable.
\item
Assume that there exists $N$ a definable, closed, discrete, and unbounded subset of~$\K_{\geq 0}$.
Then, \tfae:
{\pointlessenum
\begin{enumerate}%[i\textup )]
\item\label{en:cover-1} 
$C$ is \dcompact, 
\item\label{en:cover-2} 
for every $\Vfam$ definable open cover of $C$, if
$\Vfam$ is locally pseudo-finite, then it is pseudo-finite.
\end{enumerate}
}
\end{enumerate}
\end{proposition}
\begin{proof}
1) Since $\K$ is \dminimal, ``being pseudo-enumerable'' is equivalent to
``having dimension~$0$'', which is a first-order property.
Therefore, by Proposition~\ref{prop:dmin-polish}, \wloG $\K$ is a Polish
space, and hence $C$ is separable.
%and therefore $C$ is also a Polish space.
By Fact~\ref{fact:polish-cover}, $\Vfam$~is at most countable, and therefore
at most pseudo-enumerable.

2) The fact that $(\ref{en:cover-1} \Rightarrow \ref{en:cover-2})$ is
Lemma~\ref{lem:compact-cover}.
For the converse, assume that $C$ is not \dcompact.
Hence, \wloG $C$ is unbounded.

For every $n \in N$, let 
$p(n) \coloneqq \max\set{m \in N \cup \set{- \infty}:  m < n}$, 
and $s(n) \coloneqq \min \set{m \in N : m > n}$,
and $V_n \coloneqq \set{x \in C: p(n) < \abs x < s(n)}$.
Let $M \coloneqq \set{n \in N: V_n \neq \emptyset}$ (notice that $M$ is
unbounded, because $C$ is unbounded).
Then, $(V_n: n \in M)$ is a definable locally pseudo-finite open cover of $C$
which is not pseudo-finite.
\end{proof}

\begin{proposition}\label{prop:pair-cover}
Let $\K \coloneqq \pair{\Bm, \Am}$ be a dense pair of \dminimal structures
\rom(with $\Afam \prec \Bfam$\rom).
Assume that $\Vfam$ is definable and 
pointwise small
%such that every point of $\Am$ is in  
%``small-many'' elements of $\Bm$
\rom(that is, for every $x \in C$, $\sdim(\set{i \in I: x \in V_i}) = 0$\rom).
Then, $\Vfam$~is small.
\end{proposition}
\begin{proof}
By Proposition~\ref{prop:small-polish}\eqref{en:polish-pair}, \wloG $\Bfam$ is
Polish, $\Afam$ is a dense countable subset, and every small definable set is countable.
Therefore, by Fact~\ref{fact:polish-cover}, $\Vfam$ is countable and
definable; thus, again by
Proposition~\ref{prop:small-polish}\eqref{en:polish-pair}, $\Vfam$ is small.
\end{proof}

\begin{enumerate}
\setcounter{enumi}{6}
\item 
Let $\K \coloneqq \pair{\Bm, \Am}$ be as in the above proposition,
%dense pair of o-minimal structures
and $C$ be the interval $[0,1]$ in~$\K$.
Assume that $\K$ is nonarchimedean.
Then, there exists $\Vfam$ definable open cover of~$C$, 
such that $\Vfam$ is small but has no finite subcovers
(take $\varepsilon > 0$ infinitesimal, and let
$\Vfam \coloneqq \set{(x-\varepsilon, x + \varepsilon): x \in \Afam \cap [0, 1]}$).
\end{enumerate}

\bibliographystyle{alpha}	% uses file "filename.bst"
\bibliography{tame}		% expects file "filename.bib"

\end{document}